\documentclass[10pt,a4paper,twoside,fleqn]{article}
\usepackage{amsfonts,amssymb,amsmath,amsthm,latexsym,mathrsfs}
\usepackage{a4wide,fancyhdr}
\usepackage{graphicx,epic,eepic}
\usepackage{titlesec}
\usepackage{pgfplots}
\pgfplotsset{compat=1.16}
\usepackage{tikz}
\usetikzlibrary{quotes,angles}
\usepackage{subfigure}

\titleformat{\section}{\large\bfseries}{\thesection}{1em}{}

\numberwithin{equation}{section}

\newtheorem{theorem}{Theorem}[section]
\newtheorem{lemma}[theorem]{Lemma}
\newtheorem{proposition}[theorem]{Proposition}

\newtheorem{corollary}[theorem]{Corollary}

\theoremstyle{definition}

\newcommand{\mres}{\mathbin{\vrule height 1.6ex depth 0pt width
0.13ex\vrule height 0.13ex depth 0pt width 1.3ex}}

\def\XXint#1#2#3{{\setbox0=\hbox{$#1{#2#3}{\int}$ } 
\vcenter{\hbox{$#2#3$ }}\kern-.6\wd0}}

\begin{document}

\begin{center}
{\bf\large Some isoperimetric inequalities in the plane with radial power weights}
\end{center}

\begin{center} 
I McGillivray\\
School of Mathematics\\
University of Bristol
\end{center}

\begin{abstract}
We consider the punctured plane with volume density $|x|^\alpha$ and perimeter density $|x|^\beta$. We show that centred balls are uniquely isoperimetric for indices $(\alpha,\beta)$ which satisfy the conditions $\alpha-\beta+1>0$, $\alpha\leq 2\beta$ and $\alpha(\beta+1)\leq\beta^2$ except in the case $\alpha=\beta=0$ which corresponds to the classical isoperimetric inequality. 
\end{abstract}

\smallskip

\noindent Key words: weighted isoperimetric inequality, radial densities.

\smallskip

\noindent Mathematics Subject Classification 2020: 49K20, 49Q20, 51M16, 53A10.

\section{Introduction}

\smallskip

\noindent For $\alpha\in\mathbb{R}$ the weighted volume measure $V_\alpha:=|x|^\alpha\mathscr{L}^2$ is defined on the $\mathscr{L}^2$-measurable sets in $\mathbb{R}^2_0:=\mathbb{R}^2\setminus\{0\}$. For $\beta\in\mathbb{R}$ the weighted perimeter of a set of locally finite perimeter $E$ in $\mathbb{R}^2_0$ is defined by
\begin{equation}\label{weighted_perimeter}
P_\beta(E):=\int_{\mathbb{R}^2_0}|x|^\beta\,d|D\chi_E|\in[0,\infty].
\end{equation}
We study minimisers for the weighted isoperimetric problem
\begin{equation}\label{isoperimetric_problem}
I(v):=\inf\Big\{ P_\beta(E):E\text{ is a set with locally finite perimeter in }\mathbb{R}^2_0\text{ and }V_\alpha(E)=v\Big\}
\end{equation}
for $v>0$. Let us introduce the parameter set $\mathcal{P}$ given by
\begin{align}
\mathcal{P}&:=\Big\{(\alpha,\beta)\in\mathbb{R}^2:\alpha-\beta+1>0,\alpha\leq 2\beta\text{ and }\alpha(\beta+1)\leq\beta^2\Big\}.
\end{align}
Our main result is the following.

\smallskip

\begin{theorem}\label{main_theorem}
Let $(\alpha,\beta)\in\mathcal{P}\setminus\{(0,0)\}$ and let $v>0$. Then the centred ball is a unique minimiser for the problem (\ref{isoperimetric_problem}) up to equivalence.
\end{theorem}

\smallskip

\noindent This result may be formulated alternatively. For $(\alpha,\beta)\in\mathcal{P}\setminus\{(0,0)\}$,
\begin{align}
P_\beta(E)^{\alpha+2}\geq(2\pi)^{\alpha-\beta+1}(\alpha+2)^{\beta+1}V_\alpha(E)^{\beta+1}\label{weighted_isoperimetric_inequality}
\end{align}
for any set of locally finite perimeter $E$ in $\mathbb{R}^2_0$ with finite weighted volume and equality holds if and only if $E$ is equivalent to a centred ball.

\smallskip

\noindent This result has been formulated in \cite{Diazetal2010} Conjecture 4.22 and \cite{Alvinoetal2016} Conjecture 5.1. It arises from stability considerations as in \cite{Morganetal2016} Conjecture 3.12 (but see also  \cite{Morgan2016}). The result has been proved amongst convex competitors in \cite{DiGiosiaetal2016} and \cite{DiGiosiaetal2019} Theorem 3.1. It is proved in part in \cite{Alvinoetal2016} Theorem 1.1; in particular, for indices $(\alpha,\beta)$ satisfying $\alpha-\beta+1\geq 0$ as well as $\alpha\leq 0\leq\beta\leq 1/3$ or a more complicated condition with $\beta\geq 1/3$. In \cite{ChibaHoriuchi2016} Theorem 1.3 the result is proved in the r\'egime $\alpha-\beta+1\geq 0$ and $\alpha\leq 2\beta\leq 0$. Apart from its intrinsic interest another motivation to study this inequality is its application in the theory of Sobolev spaces. In \cite{Alvinoetal2016} the weighted isoperimetric inequality is used to obtain the best constant in a  corresponding Caffarelli-Kohn-Nirenberg inequality (see \cite{Alvinoetal2016} Theorem 8.3 or \cite{ChibaHoriuchi2016} for example). 

\smallskip

\begin{figure}[h]
\centering
\begin{tikzpicture}[>=stealth]
    \begin{axis}[
        xmin=-2.5,xmax=3.5,
        ymin=-1.5,ymax=3.5,
        axis x line=middle,
        axis y line=middle,
        axis line style=->,
        xlabel={$\alpha$},
        ylabel={$\beta$},
        ]
        \addplot[black, domain=-2:0, smooth]{x/2};
        \addplot[black, dashed, domain=0:3.5, smooth]{x/2};

        \addplot[dotted, domain=0:3.5, smooth]{x};
        \draw[black] (-2,-1) -- (-2,3.5);

        \addplot[black, domain=0:3.5, smooth]{(x+sqrt(x^2+4*x))/2};
    \end{axis}
\end{tikzpicture}
\caption{In the region between the solid lines and curve centred balls are isoperimetric (see Theorem \ref{main_theorem} and \cite{Diazetal2010}, \cite{Howe2016}). The solid curve is given by $\alpha(\beta+1)-\beta^2=0$ with $\beta\geq 0$. Balls tangent to the origin are isoperimetric on the dotted line (see \cite{Dahlbergal2010} Theorem 3.16). In the region between the dotted and dashed lines pinched circles through the origin are isoperimetric (see \cite{Diazetal2010} Proposition 4.21). For the region between the solid curve and dotted line see \cite{Diazetal2010} Conjecture 4.22.} 
\label{fig:Summary_of_results}
\end{figure}
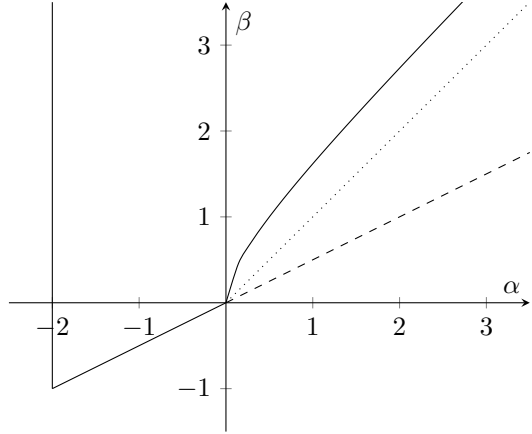

\noindent Let us mention some related results. Centred balls are isoperimetric in case $\alpha-\beta+1\leq 0$ and $\alpha+2>0$ as shown in \cite{Diazetal2010} Proposition 4.21 (see also \cite{Howe2016} Example 3.5 (4)). In \cite{Dahlbergal2010} Theorem 3.16 it is shown that discs tangent to the origin are isoperimetric in the case $\alpha=\beta>0$. In contrast, isoperimetric minimisers do not exist if $-2<\alpha<0$ and $\alpha=\beta$ according to \cite{Carrolletal2007} Proposition 4.2. If $\beta\leq\alpha\leq 2\beta$ and $\beta\geq 0$ then pinched circles through the origin are isoperimetric. This is contained in \cite{Diazetal2010} Proposition 4.21. The case $\beta\geq 1$ and $\alpha=0$ is treated in \cite{Bettaetal1999} Theorem 2.1. A similar isoperimetric problem but with additional constraints is discussed in the case $\alpha=0$ and $\beta<0$ in \cite{Csato2017} Theorem 3. We refer to Figure \ref{fig:Summary_of_results}.

\section{Outline of the argument}

\smallskip

\noindent Let us give the bones of the argument.  

\smallskip

\noindent Select a minimiser $E$ of (\ref{isoperimetric_problem}). Then $E$ is essentially bounded. Existence and boundedness is guaranteed in Corollary \ref{existence_and_boundedness_for_indices_in_curly_Q}. Replace $E$ by its spherical cap symmetral $E^{sc}$. This is also a minimiser because the perimeter and volume densities are radial.
By choosing a suitable $\mathscr{L}^2$-version of the set  $E^{sc}$ we can arrange that its topological boundary is well-behaved in a measure-theoretic sense. Call this version $\widetilde{E}$. It has analytic boundary $M$ in $\mathbb{R}^2_0$. So we may choose $\widetilde{E}$ to be open, bounded, spherical cap symmetric and with analytic boundary relative to $\mathbb{R}^2_0$. The altered minimiser $\widetilde{E}$ is related to the original minimiser $E$ via the relation $L_{\widetilde{E}}=L_E$ a.e. on $(0,\infty)$. Here, $L_E(t)$ stands for the $\mathscr{H}^1$-measure of the trace of $E$ on the centred circle with radius $t>0$ (and likewise for $\widetilde{E})$. These facts are contained in Theorem \ref{C1_property_of_reduced_boundary}.

\smallskip

\begin{figure}[h]
\centering
\begin{tikzpicture}[>=stealth]
\draw[->] (-0.2,0) -- (3,0);
\draw[->] (0,-0.2) -- (0,3);
\node[inner sep=1.2pt,fill,circle] at (1,1) (point) {};
\draw [thick] (1,1) to[out=135,in=20] (0.2,1.2);
\draw [thick] (1,1) to[out=-45,in=120] (1.4,0.4);
\draw[dashed] (0,0) to (1,1);
\draw[->] (1,1) to (1.7,1.7) node[right]{$x$}; 
\draw[->] (1,1) to (0.3,1.7) node[above]{$t$}; 
\coordinate (first) at (1,1);
\coordinate (second) at (1.5,1.5);
\coordinate (third) at (0.5,1.5);
\pic [draw, ->, "$\sigma$", angle eccentricity=1.5] {angle=second--first--third};  
\end{tikzpicture}
\caption{The inner unit normal along the boundary of $E$ is denoted $n$. The tangent vector $t$ is chosen in such a way that the pair $\{t,n\}$ is positively oriented. The angle $\sigma$ stands for the angle between the position vector $x$ and the tangent vector $t$ measured in an anticlockwise direction.} 
\label{fig:angle_sigma}
\end{figure}
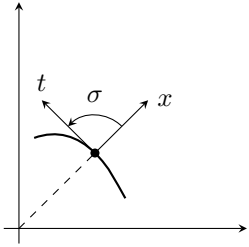

\noindent Let $n$ stand for the inner unit normal vector field along $M$. Choose a tangent vector $t$ in such a way that the pair $\{t,n\}$ is positively oriented. The angle $\sigma$ stands for the angle between the position vector $x$ and the tangent vector $t$ measured in an anticlockwise direction (as in Figure \ref{fig:angle_sigma}). Define the set $\Omega$ to be the collection of radii for which the centred circle with this radius meets the boundary curve $M$ transversally. This set is open in $(0,\infty)$. The heart of the demonstration is to show that $\Omega$ is empty. In other words centred circles meet the boundary $M$ tangentially if at all. This result is contained in Theorem \ref{Omega_is_empty}. We shall say more about this shortly. With this result in hand we may deduce that $\widetilde{E}$ is a finite union of centred annuli. This leads to the fact that $\widetilde{E}$ is a centred ball. The proof of uniqueness exploits the a.e. relation $L_{\widetilde{E}}=L_E$ mentioned above. 

\smallskip

\noindent Let us say a little more about how to show that $\Omega$ is empty. Suppose for a contradiction that $\Omega\neq\emptyset$. The function $u$ at a point $\tau$ in $\Omega$ is defined to be the sine of the angle $\sigma$ introduced above. This function is continuously differentiable on $\Omega$ and there exists a real constant $\lambda$  such that
\begin{align}
u^\prime+\frac{\beta+1}{\tau}u+\lambda\tau^{\alpha-\beta}&=0\label{ode_for_u_1}
\end{align}
on $\Omega$. This last is a reformulation of the constant generalised (mean) curvature condition; that is, constancy of the expression
\[
\tau^{\beta-\alpha}\Big\{k+\frac{\beta}{\tau}u\Big\}
\]
along $M$. The curvature $k$ of $M$ as a function of the radial variable $\tau$ is given by the expression $k=\frac{1}{\tau}(\tau u)^\prime$. This means that the curve $M$ can be reconstructed from the function $u$ up to rotation around the origin. In fact the angular variable admits an explicit integral representation. Suppose that $(a,b)$ is a maximal connected component of  $\Omega$. Then the angular variable of a point with polar coordinates $(\tau,\theta)$ on $M$ is given by
\[
\theta_u=\theta_u(b)+\int_\cdot^b\frac{u}{\sqrt{1-u^2}}\frac{dt}{t}
\]
on $(a,b)$. The angle integral is a hyper-elliptic integral in case the indices $\alpha$ and $\beta$ are integral.

\smallskip

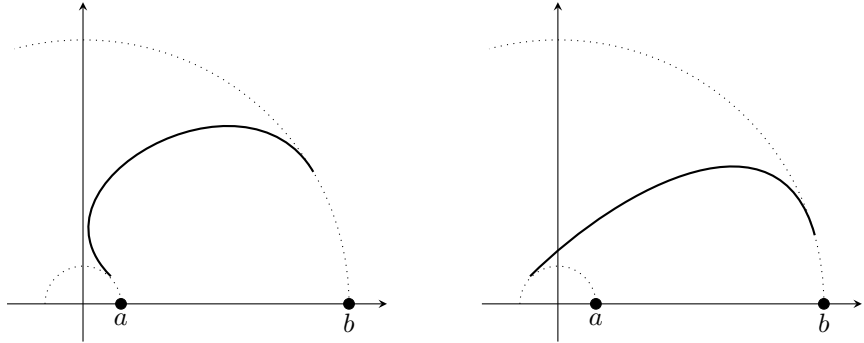
\begin{figure}[h]
\centering
\begin{subfigure}
\centering
\begin{tikzpicture}[>=stealth][scale=0.33]
\draw[->] (-1,0) -- (4,0);
\draw[->] (0,-0.5) -- (0,4);
\draw[dotted] (-0.5,0) arc (180 : 0 : 0.5);
\draw[dotted] (3.5,0) arc (0 : 105 : 3.5);
\draw[thick] (0.366,0.366) to[out=135,in=120,,looseness=1.5] (3.03,1.75);

\coordinate (A) at   (0.5,0);
\coordinate (B) at   (3.5,0);

\filldraw (A) circle (2pt) node[below] {$a$};
\filldraw (B) circle (2pt) node[below] {$b$};

\end{tikzpicture}
\end{subfigure}
\begin{subfigure}
\centering
\hspace{1cm}
\begin{tikzpicture}[>=stealth][scale=0.33]
\draw[->] (-1,0) -- (4,0);
\draw[->] (0,-0.5) -- (0,4);
\draw[dotted] (-0.5,0) arc (180 : 0 : 0.5);
\draw[dotted] (3.5,0) arc (0 : 105 : 3.5);
\draw[thick] (-0.366,0.366) to[out=45,in=105,looseness=1.2] (3.38,0.91);

\coordinate (A) at   (0.5,0);
\coordinate (B) at   (3.5,0);

\filldraw (A) circle (2pt) node[below] {$a$};
\filldraw (B) circle (2pt) node[below] {$b$};
\end{tikzpicture}
\end{subfigure}
\label{fig:boundary_conditions}
\caption{The diagram on the left gives an illustration of the boundary  condition $-u(a)=u(b)=1$; that on the right illustrates the boundary condition $u(a)=u(b)=1$.} 
\end{figure}

\noindent Suppose that the left-hand end-point $a$ is positive. The cosine of $\sigma$ vanishes at the end-points $a$ and $b$. So $u$ takes values $\pm 1$ there. Suppose first of all that $-u(a)=u(b)=1$ (see Figure \ref{fig:boundary_conditions}). Then $\theta_u(a)\in(0,\pi]$ by Corollary \ref{integral_of_solution_of_first_order_ode} and in fact in $(0,\pi)$. As the sine of $\sigma$ is $-1$ at this point the curvature is $-1/a$; but this conflicts with the curvature obtained from the generalised (mean) curvature equation. Now suppose that $u(a)=u(b)=1$. Then $\theta_u(a)>\pi$ (Corollary \ref{integral_of_w}); again a contradiction. We turn to the case $a=0$. Then the weighted perimeter of $E$ relative to the closed centred ball with radius $b$ exceeds that of the centred ball with the same weighted volume. This result is contained in (\ref{relation_between_isoperimetric_ratios}). This involves a delicate beta function inequality in Theorem \ref{beta_function_inequality}. We refer to the proof of Theorem \ref{Omega_is_empty} for a discussion of the remaining boundary conditions. We have mentioned two inequalities involving the angle integral. Let us supply a little more detail.

\smallskip

\noindent Consider once more the boundary condition $-u(a)=u(b)=1$. The angle integral can be expressed in terms of the distribution function $\mu_{\pm u}$ of $\pm u$ with respect to the measure $\mu$ on $(0,\infty)$ with infinitessimal element $(1/\tau)\,d\tau$; to be more explicit $\mu_u(t)=\mu(\{u>t\})$. In detail,
\[
\theta_u(b)-\theta_u(a)=-\int_a^b\frac{u}{\sqrt{1-u^2}}\frac{d\tau}{\tau}
=-\int_0^1\Big\{\mu_u-\mu_{-u}\Big\}\frac{dt}{(1-t^2)^{3/2}}
\]
after a change in the order of integration. We aim to show that the right-hand side is strictly negative. This entails that $\theta_u(a)\in(0,\pi]$ as mentioned in the last paragraph. The difference between the distribution functions can be expressed in terms of their derivatives
\[
-\Big\{\mu_u(t)-\mu_{-u}(t)\Big\}=\int_{(t,1]}\Big[\mu_u^\prime-\mu_{-u}^\prime\Big]\,ds.
\]
By the coarea formula
\[
\mu_u^\prime(t)-\mu_{-u}^\prime(t)=-\int_{\{u=t\}}\frac{1}{u^\prime}\frac{d\mathscr{H}^0}{\tau}
+\int_{\{u=-t\}}\frac{1}{u^\prime}\frac{d\mathscr{H}^0}{\tau}.
\]
Let us imagine for a moment that $t=1$ so that $\{u=\pm t\}$ coincides with the point $b$ respectively $a$. Then this becomes 
\[
\mu_u^\prime(1)-\mu_{-u}^\prime(1)=-\frac{1}{bu^\prime(b)}+\frac{1}{au^\prime(a)}
\]
As can be seen from the ODE in (\ref{ode_for_u_1}) the derivative of $u$ can be expressed in terms of $\lambda$, $u$ and $\tau$. We thus require
\[
-\frac{1}{bu^\prime(b)}+\frac{1}{au^\prime(a)}=\frac{1}{\beta+1+\lambda b^\gamma}+\frac{1}{\beta+1-\lambda a^\gamma}<0
\] 
where $\gamma=\alpha-\beta+1$. The multiplier $\lambda$ depends explicitly on the end-points $a$ and $b$ as well as the parameters $\alpha$ and $\beta$. This results in an algebraic inequality. The necessary inequality is established in Proposition \ref{inequality_for_hat_m}. Now suppose that $t\in(0,1)$; that is, $u$ takes the values $-t$ and $+t$ at the end-points of the interval $\{|u|<t\}$. The scaled function $u/t$ has boundary values $\pm 1$ and satisfies an ODE similar to (\ref{ode_for_u}) but with a scaled $\lambda$. Strict negativity of $\mu^\prime_u(t)-\mu^\prime_{-u}(t)$ follows as just described.

\smallskip

\noindent Now let us turn to the boundary condition $u(a)=u(b)=1$ with again $0<a<b<\infty$. In this case we consider the reciprocal $w:=1/u$ of $u$. This satisfies the Riccati equation
\[
w^\prime+\lambda\tau^{\alpha-\beta}w^2=\frac{\beta+1}{\tau}w
\]
as can be seen from (\ref{ode_for_u_1}) but with boundary conditions $w(a)=w(b)=1$. Likewise the angle integral is
\[
\theta_u(b)-\theta_u(a)=-\int_a^b\frac{1}{\sqrt{w^2-1}}\frac{dt}{t}.
\]
Expressed in terms of the distribution function we arrive at a formulation that differs from the above
\[
\int_a^b\frac{1}{\sqrt{w^2-1}}\frac{d\tau}{\tau}=\pi-\int_1^\infty\Big[\mu_w-\mu_{w_0}\Big]\frac{t}{(t^2-1)^{3/2}}\,dt.
\]
We establish the differential inequality $-\mu_w^\prime>\sigma(t,\mu_w)$ where equality $-\mu_{w_0}^\prime=\sigma(t,\mu_{w_0})$ pertains for the distribution function associated to the case $\alpha=\beta=0$ and corresponds to a semicircle through the points $(-a,0)$ and $(b,0)$. This differential inequality entails that  $\mu_w<\mu_{w_0}$ and the angle integral in the last displayed equation exceeds $\pi$. This results in the contradiction that $\theta_u(a)>\pi$ as mentioned above.

\smallskip

\noindent Let us give some indication as how to obtain the last-mentioned differential inequality. Again by the coarea formula
\[
-\mu_w^\prime(t)+\mu_{w_0}^\prime(t)=\int_{\{w=t\}}\frac{1}{|w^\prime|}\frac{d\mathscr{H}^0}{\tau}
-\int_{\{w_0=t\}}\frac{1}{|w_0^\prime|}\frac{d\mathscr{H}^0}{\tau}.
\]
Let us imagine that $t=1$. Then this becomes 
\[
-\mu_w^\prime(1)+\mu_{w_0}^\prime(1)
=-\frac{1}{bw^\prime(b)}+\frac{1}{aw^\prime(a)}
+\frac{1}{bw_0^\prime(b)}-\frac{1}{aw_0^\prime(a)}
\]
As before the multiplier $\lambda$ depends explicitly on the end-points $a$ and $b$ as well as the parameters $\alpha$ and $\beta$. This results in an algebraic inequality. The necessary inequality is established in Theorem \ref{main_algebraic_inequality} and is somewhat involved. Now suppose that $t\in(0,1)$; that is, $w$ takes the value $t$ at the end-points of the interval $\{w>t\}$. The scaled function $w/t$ has boundary value $1$ on an interval contained inside $(a,b)$ and satisfies a Riccati equation there but with a scaled $\lambda$. The above argument leads to the inequality $\mu_w(t)<\mu_{w_0}(t)$.

\smallskip

\noindent Finally, we mention that Sections \ref{A_tangential_Volpert_Theorem}-\ref{Further_preliminary_results} contain some preliminary material that is needed in the proofs of existence and boundedness contained in Section \ref{Existence_and_boundedness}.

\section{A tangential Vol'pert Theorem}\label{A_tangential_Volpert_Theorem}

\smallskip

\noindent The main result of this Section is a tangential Vol'pert theorem given in Theorem \ref{tangentialVolperttheorem}. This result appears in \cite{CagnettiPeruginiStoger2020} Theorem 3.7. This last is couched in the language of integer multiplicity rectifiable currents. We prefer to work here in the framework of sets of locally finite perimeter. The proof is an adaptation of the approach in \cite{Barchiesietal2013} Theorem 2.4.

\smallskip

\noindent{\em Variation on the sphere.} The focus of this subsection is a De Giorgi rectifiability theorem on the circle (cf. \cite{Ambrosio2000} Theorem 3.59). We use the notation $\mathbb{S}^{1}_\tau$ to stand for the sphere in $\mathbb{R}^2$ with centre at the origin and radius $\tau>0$. We often drop the superscript; $\mathbb{S}$ stands for the unit sphere. The induced Riemannian metric on $\mathbb{S}$ will be denoted by $\langle\cdot,\cdot\rangle_{\mathbb{S}}$ and the tangential gradient and divergence by $\nabla_\mathbb{S}$ and $\mathrm{div}_\mathbb{S}$ (see \cite{Ambrosio2000} Section 7.3 for example). Let $\Omega$ be an open set in $\mathbb{S}$. For $u\in L^1(\Omega,\mathscr{H}^{1})$ define its variation on the sphere by
\begin{equation}\label{variation_on_the_sphere}
V_{\mathbb{S}}(u,\Omega):=
\sup\left\{
\int_{\mathbb{S}}u\,\mathrm{div}_{\mathbb{S}}X\,d\mathscr{H}^{1}
:X\in C^1_c(\Omega,\mathbb{R}^2),\|X\|_\infty\leq 1\right\}\in[0,\infty].
\end{equation}
We say that $u$ has bounded variation on the sphere if $u\in L^1(\Omega,\mathscr{H}^{1})$ and $V_{\mathbb{S}}(u,\Omega)<\infty$ and write $u\in\mathrm{BV}_{\mathbb{S}}(\Omega)$. 

\smallskip

\noindent Fix $x\in\mathbb{R}^2$.  The projection
\[
(\pi v)(x)=\pi(x)v:=
\left\{
\begin{array}{ll}
v - |x|^{-2}\,\langle v,\, x\rangle\,x & \text{for }x\neq 0,\\
0                                                                  & \text{for }x=0,\\
\end{array}
\right. 
\]
maps the vector $v\in\mathbb{R}^2$ onto the subspace $x^\perp$ perpendicular to $x$. Note that $\pi$ is symmetric in the sense that
$
\langle(\pi v)(x), w\rangle = \langle v,(\pi w)(x)\rangle
$
for any $v,\,w\in\mathbb{R}^2$. 

\smallskip

\begin{lemma}\label{properties_of_variation_on_sphere}
Let $\Omega$ be an open set in $\mathbb{S}$. Then
\begin{itemize}
\item[(i)] the functional $u\mapsto V_\mathbb{S}(u,\Omega)$ is lower semicontinuous with respect to convergence in $L^1(\Omega,\mathscr{H}^{1})$.
\end{itemize}
Assume that $u\in\mathrm{BV}_{\mathbb{S}}(\Omega)$. Then
\begin{itemize}
\item[(ii)]
there exists a unique $\mathbb{R}^2$-valued finite Radon measure $D_{\mathbb{S}}u$ on $\Omega$ such that
\[
\int_\Omega u\,\mathrm{div}_{\mathbb{S}}\,\pi X\,d\mathscr{H}^{1}
=
-\int_\Omega\langle X,\,dD_{\mathbb{S}}u\rangle_{\mathbb{S}}
\text{ for any }X\in C^1_c(\Omega,\mathbb{R}^2)
\]
and $|D_\mathbb{S} u|(\Omega)=V_\mathbb{S}(u,\Omega)$;
\item[(iii)]
$D_\mathbb{S} u=\nabla_\mathbb{S} u\,\mathscr{H}^{1}$ 
and 
$|D_\mathbb{S}u|=|\nabla_\mathbb{S} u|\,\mathscr{H}^{1}$ 
for $u\in C^1(\Omega)\cap\mathrm{BV}_\mathbb{S}(\Omega)$.
\end{itemize}
\end{lemma}

\smallskip

\noindent{\em Proof.}
{\em (i)} follows from \cite{Ambrosio2000} Remark 3.5 while {\em (ii)} is a consequence of Riesz's Theorem (see \cite{Ambrosio2000} Theorem 1.54).
\qed

\smallskip

\noindent Recall that the set $E$ has locally finite perimeter in $\mathbb{R}^2_0$ if $\chi_E$ belongs to $\mathrm{BV}_{\mathrm{loc}}(\mathbb{R}^2_0)$. The reduced boundary $\mathscr{F}E$ of $E$ is defined by
\[
\mathscr{F}E:=\Big\{x\in\mathrm{supp}|D\chi_E|:\nu_E(x):=\lim_{\rho\downarrow 0}\frac{D\chi_E(B(x,\rho))}{|D\chi_E|(B(x,\rho))}\text{ exists in }\mathbb{R}^2\text{ and }|\nu^E(x)|=1\Big\}
\]
as in for example \cite{Ambrosio2000} Definition 3.54. Note that $\mathscr{F}E$ is a Borel set and the map $\nu_E:\mathscr{F}E\rightarrow\mathbb{S}$ is Borel (cf. \cite{Ambrosio2000} Theorem 2.22). 

\smallskip

\begin{lemma}\label{on_cone_set}
Let $E$ be a set of locally finite perimeter in $\mathbb{R}^2_0$. Assume that $E$ is scale invariant; that is, $\lambda E=E$ for each $\lambda>0$. Then
\begin{itemize}
\item[(i)] $\mathscr{F}E$ is scale-invariant;
\item[(ii)] $\nu_E(x)=\nu_E(\lambda x)$ for $x\in\mathscr{F}E$ and $\lambda>0$;
\item[(iii)] $\nu_E(x)\perp x$ for each $x\in\mathscr{F}E$.
\end{itemize}
\end{lemma}

\smallskip

\noindent{\em Proof.}
{\em (i)} A scaling argument gives
\[
D\chi_E(B(\lambda x,\lambda r))=\lambda D\chi_E(B(x,r))
\text{ and }
|D\chi_E|(B(\lambda x,\lambda r))=\lambda|D\chi_E|(B(x,r))
\]
for $x\in\mathbb{R}^2_0$ and $0<r<|x|$. This leads to the first claim. The considerations in {\em (i)} entail {\em (ii)}. {\em (iii)} Let $\phi\in C^1_c((0,\infty))$ and $\psi\in C^1(\mathbb{S})$ and define a $C^1_c$ vector field on $\mathbb{R}^2_0$ by $X(x):=\phi(\tau)\psi(\omega)\omega$ for $x=\tau\omega\in\mathbb{R}^2_0$. Then
\[
-\int_{\mathscr{F}E}\phi\psi\langle\omega,\nu_E\rangle\,d\mathscr{H}^1
=\int_E\mathrm{div}\,X\,dx=\int_E\frac{1}{\tau}(\tau\phi)^\prime\psi\,dx
=\int_\mathbb{S}\int_0^\infty\chi_E(\tau\phi)^\prime\,d\tau\psi\,d\mathscr{H}^1
=0
\]
upon converting to polar coordinates. The function $\mathbb{S}\rightarrow\mathbb{R};\omega\mapsto\langle\nu_E(\omega),\omega\rangle$ is Borel. By Lusin's Theorem (cf. \cite{Ambrosio2000} Theorem 1.45) we derive
\[
\int_{\mathscr{F}E\cap\{\langle\cdot,\nu_E\rangle>0\}}\phi\langle\omega,\nu_E\rangle\,d\mathscr{H}^1=0
\]
from which we infer that the set $\mathscr{F}E\cap\{\langle\cdot,\nu_E\rangle>0\}$ is an $\mathscr{H}^1$-null set. We arrive at a similar conclusion with the inequality reversed. The result follows.
\qed

\bigskip

\noindent For a Borel set $E$ in $\mathbb{R}^2_0$ and $\tau>0$ the $\tau$-section of $E$ is defined by $E_\tau:=\{x\in E: |x|=\tau\}$. For a real-valued Borel function $u$ on $\mathbb{R}^2_0$ define $u_\tau$ to be the restriction of $u$ to $\mathbb{S}_\tau$.

\smallskip

\begin{proposition}\label{tangential_integral_equality}
Let $E$ be a set of locally finite perimeter in $\mathbb{R}^2_0$ and let $g:\mathbb{R}^2_0\rightarrow[0,\infty]$ be a Borel function. Then
\[
\int_{\mathscr{F}E}g|\nu_E^\omega|\,d\mathscr{H}^{1}
=\int_0^\infty 
\int_{(\mathscr{F}E)_\tau}g\,d\mathscr{H}^{0}\,d\tau
\]
where $\nu_E^\omega=\pi\nu_E$ stands for the tangential component of the inner unit normal vector $\nu_E$.
\end{proposition}

\smallskip

\noindent{\em Proof.}
Define $f:\mathbb{R}^2_0\rightarrow(0,\infty)$ via $x\mapsto |x|$. Let $x\in\mathbb{R}^2_0$ and $M$ be a line through $x$ perpendicular to the unit vector $v$. Then 
\begin{align}
C_1d^Mf_x&=(\pi v)(x).\label{formula_for_C_1_of_differential}
\end{align}
The result follows by the coarea formula \cite{Ambrosio2000} Theorem 2.93 and (2.72).
\qed

\medskip

\noindent An $\mathscr{H}^{1}$-measurable set $E$ in $\Omega$ is said to have finite perimeter in $\Omega$ if $\chi_E\in\mathrm{BV}_{\mathbb{S}}(\Omega)$. The reduced boundary $\mathscr{F}_{\mathbb{S}}E$ of $E$ consists of the set of points $x\in\mathrm{supp}|D_{\mathbb{S}}\chi_E|\cap\Omega$ such that
\[
\nu_E^{\mathbb{S}}(x):=
\lim_{\rho\downarrow 0}
\frac{D_{\mathbb{S}}\chi_E(B_\rho(x))}{|D_{\mathbb{S}}\chi_E|(B_\rho(x))}
\]
exists in $\mathbb{R}^2$ and has unit length. We may now state and prove a De Giorgi rectifiability theorem on the circle (cf. \cite{Ambrosio2000} Theorem 3.59).

\smallskip

\begin{theorem}\label{deGiorgivariationsphere}
Let $\Omega$ be an open set in $\mathbb{S}$ and $E$ a set of finite perimeter in $\Omega$. Then
\[
D_{\mathbb{S}}\chi_E = \nu_E^{\mathbb{S}}\mathscr{H}^{0}\mres\mathscr{F}_{\mathbb{S}}E
\text{ and }
|D_{\mathbb{S}}\chi_E| = \mathscr{H}^{0}\mres\mathscr{F}_{\mathbb{S}}E.
\]
\end{theorem}

\smallskip

\noindent{\em Proof.}
Let
\[
\widetilde{\Omega}:=\{x\in\mathbb{R}^2_0:x/|x|\in\Omega\}
\]
be the cone over $\Omega$ and define $\widetilde{E}$ likewise. Let $X\in C^1_c(\widetilde{\Omega},\mathbb{R}^2)$. We first establish the identity
\begin{align}
\int_{\widetilde{\Omega}}\chi_{\widetilde{E}}\,\mathrm{div}\,X\,dx
&=\int_0^\infty\int_{\Omega}\chi_{E}\,\mathrm{div}(\pi\Big[X\circ I^\tau\Big])\,d\mathscr{H}^1\,d\tau\label{integral_identity}
\end{align}
where $I^\tau(x)=\tau x$ is the homothety with scale factor $\tau$. It is apparent from this identity that $\widetilde{E}$ is a set of locally finite perimeter in $\widetilde{\Omega}$ because $E$ is a set of finite perimeter in $\Omega$. We decompose $X=\pi X+(I-\pi)X$ and note that $(I-\pi)X=\phi\omega$ where $\phi(x)=\langle X,\omega\rangle$ and $\omega=x/|x|$. Converting to polar coordinates,
\begin{align}
\int_{\widetilde{\Omega}}\chi_{\widetilde{E}}\,\mathrm{div}\,X\,dx
&=\int_{\widetilde{\Omega}}\chi_{\widetilde{E}}\,\mathrm{div}\,\pi X\,dx
+\int_{\widetilde{\Omega}}\chi_{\widetilde{E}}\,\mathrm{div}(I-\pi)X\,dx\nonumber\\
&=\int_0^\infty\int_{\widetilde{\Omega}_\tau}\chi_{\widetilde{E}}\mathrm{div}(\pi X)\,d\mathscr{H}^1\,d\tau
+\int_{\Omega}\chi_E\int_0^\infty\partial_\tau(\tau\phi)\,d\tau\,d\mathscr{H}^1\nonumber\\
&=\int_0^\infty\int_{\widetilde{\Omega}_\tau}\chi_{\widetilde{E}}\,\mathrm{div}(\pi X)\,d\mathscr{H}^1\,d\tau
=\int_0^\infty\int_{\Omega}\chi_{E}\,\mathrm{div}(\pi\Big[X\circ I^\tau\Big])\,d\mathscr{H}^1\,d\tau
\nonumber
\end{align}
using the fact that $\mathrm{div}(I-\pi)X=\tau^{-1}\partial_\tau(\tau\phi)$ and $\mathrm{div}(\pi X\circ I^\tau)(x)=\tau\,\mathrm{div}(\pi X)(\tau x)$. This establishes (\ref{integral_identity}). 

\smallskip

\noindent Let $\phi\in C_c^1((0,\infty))$. Let $Y$ be an $\mathbb{R}^2$-valued vector field defined on $\Omega$ which is of class $C^1$ with compact support and which is tangential to $\mathbb{S}$. Define $X\in C^1_c(\widetilde{\Omega},\mathbb{R}^2)$ by $X(x):=\phi(\tau) Y(\omega)$ for $x=\tau\omega\in\widetilde{\Omega}$. The set $\widetilde{E}$ is scale-invariant. So $|\pi\nu_{\widetilde{E}}|=1$ on $\mathscr{F}\widetilde{E}$ by Lemma \ref{on_cone_set} which implies in turn that $|\nu_{\widetilde{E}}^\omega|=1$. By Proposition \ref{tangential_integral_equality} the left-hand side of (\ref{integral_identity}) may be written
\begin{align}
\int_{\widetilde{\Omega}}\chi_{\widetilde{E}}\,\mathrm{div}\,X\,dx
&=-\int_{\mathscr{F}\widetilde{E}}\langle X,\nu_{\widetilde{E}}\rangle\,d\mathscr{H}^1\nonumber
=-\int_0^\infty \int_{(\mathscr{F}\widetilde{E})_\tau}\langle X,\nu_{\widetilde{E}}\rangle\,d\mathscr{H}^{0}\,d\tau\nonumber\\
&=-\int_0^\infty\phi\int_{(\mathscr{F}\widetilde{E})_\tau}\langle Y,\nu_{\widetilde{E}}\rangle\,d\mathscr{H}^0\,d\tau.
\label{LHS_of_integral_identity}
\end{align}
The right-hand side of (\ref{integral_identity}) may be expressed
\begin{align}
\int_0^\infty\int_{\Omega}\chi_{E}\,\mathrm{div}(\pi\Big[X\circ I^\tau\Big])\,d\mathscr{H}^1\,d\tau
&=
-\int_0^\infty\int_\Omega\langle X\circ I^\tau,dD_{\mathbb{S}}\chi_E\rangle_{\mathbb{S}}\,d\tau\nonumber\\
&=-\int_0^\infty\phi\int_\Omega\langle Y,dD_{\mathbb{S}}\chi_E\rangle_{\mathbb{S}}\,d\tau
\label{RHS_of_integral_identity}
\end{align}
by Lemma \ref{properties_of_variation_on_sphere}. On combining (\ref{LHS_of_integral_identity}) and (\ref{RHS_of_integral_identity}) using (\ref{integral_identity}) we may equate
\[
\int_0^\infty\phi\int_{(\mathscr{F}\widetilde{E})_\tau}\langle Y,\nu_{\widetilde{E}}\rangle\,d\mathscr{H}^0\,d\tau
=\int_0^\infty\phi\int_\Omega\langle Y,dD_{\mathbb{S}}\chi_E\rangle_{\mathbb{S}}\,d\tau
\] 
which leads in turn to the identity
\[
\int_{(\mathscr{F}\widetilde{E})_1}\langle Y,\nu_{\widetilde{E}}\rangle\,d\mathscr{H}^0
=
\int_\Omega\langle Y,dD_{\mathbb{S}}\chi_E\rangle_{\mathbb{S}}.
\]
From the uniqueness property in Riesz's Theorem (cf. \cite{Ambrosio2000} Theorem 1.54) we deduce that
$D_\mathbb{S}\chi_E=\nu_{\widetilde{E}}\mathscr{H}^{0}\mres(\mathscr{F}\widetilde{E})_1$. We derive that $\mathscr{F}E=(\mathscr{F}\widetilde{E})_1$ and $\nu_E^{\mathbb{S}}=\nu_{\widetilde{E}}$ there. This concludes the proof. 
\qed
 
\bigskip

\noindent{\em Spherical variation.}  Let $\Omega$ be an open set in $\mathbb{R}^2_0$. The spherical variation of $u\in L^1_{\mathrm{loc}}(\Omega)$ on $A\subset\subset\Omega$ is defined by
\[
V_\omega(u,A):=\sup\Big\{\int_A u\,\text{div}(\pi X)\,dx:X\in C^1_c(A,\mathbb{R}^2)\text{ and }\|X\|_\infty\leq 1\Big\}\leq\infty.
\]

\smallskip

\begin{lemma}\label{properties_of_spherical_variation}
Let $\Omega$ be an open set in $\mathbb{R}^2_0$. Then
\begin{itemize}
\item[(i)] the functional $u\mapsto V_\omega(u,A)$ is lower semicontinuous with respect to convergence in $L^1_{\mathrm{loc}}(\Omega)$ for each $A\subset\subset\Omega$.
\end{itemize}
Assume that $u\in\mathrm{BV}_{\mathrm{loc}}(\Omega)$. Then
\begin{itemize}
\item[(ii)]
there exists a unique $\mathbb{R}^2$-valued Radon measure $D_\omega u$ such that
\[
\int_\Omega u\,\mathrm{div}(\pi X)\,dx=-\int_\Omega\langle X,\,dD_\omega u\rangle\text{ for any }
X\in C^1_c(\Omega,\mathbb{R}^2)
\]
and $|D_\omega u|(\Omega)=V_\omega(u,\Omega)$;
\item[(iii)]
$D_\omega u=(\pi\rho)\,|Du|$ and $|D_\omega u|=|\pi\rho|\,|Du|$ where $Du=\rho\,|Du|$ is the polar decomposition of $D u$ with $\rho$ a unique $\mathbb{S}$-valued function in $L^1(\Omega,\,|Du|)^2$;
\item[(iv)]
$D_\omega u=\nabla_\omega u\,\mathscr{L}^2$ 
and 
$|D_\omega u|=|\nabla_\omega u|\,\mathscr{L}^2$ 
for $u\in C^1(\Omega)\cap\mathrm{BV}_{\mathrm{loc}}(\Omega)$ where $\nabla_\omega = \pi\nabla$.
\end{itemize}
\end{lemma}

\smallskip

\noindent{\em Proof.}
{\em (i)} and {\em (ii)} follow from \cite{Ambrosio2000} Remark 3.5 and Riesz's Theorem (cf. \cite{Ambrosio2000} Theorem 1.54). The fact that $u$ belongs to the class $\mathrm{BV}_{\mathrm{loc}}(\Omega)$ leads to {\em (iii)} and {\em (iv)}. 
\qed

\smallskip

\begin{lemma}\label{spherical_gradient_of_mollified_function}
Let $\Omega$ be an open set in $\mathbb{R}^2_0$. Let $u\in BV(\Omega)$ and $U\subset\subset\Omega$ open. Let $(\varphi_\varepsilon)_{\varepsilon>0}$ be a family of mollifiers on $\mathbb{R}^2$ as in \cite{Ambrosio2000} Section 2.1. Then $\nabla_\omega(u\ast\varphi_\varepsilon)=\pi\,(Du\ast\varphi_\varepsilon)$ on $U$ for  small $\varepsilon>0$.
\end{lemma}

\smallskip

\noindent{\em Proof.}
Let  $0<\varepsilon<\mathrm{dist}(U,\partial\Omega)$. Let $X\in C^1_c(\Omega,\mathbb{R}^2)$ with $\mathrm{supp}\,X\subset U$. We use the fact that
$\mathrm{div}(\varphi_\varepsilon\ast\pi X)=\varphi_\varepsilon\ast\mathrm{div}(\pi X)$ on $\Omega$. Note that the restriction of $u\ast\varphi_\varepsilon$ to $U$ belongs to $\mathrm{BV}(U)$. By symmetry of $\pi$, Lemma \ref{properties_of_spherical_variation} and the convolution identity \cite{Ambrosio2000} (2.3),
\begin{align}
-\int_U\langle\nabla_\omega(u\ast\varphi_\varepsilon),X\rangle\,dx
&=
\int_U(u\ast\varphi_\varepsilon)\mathrm{div}\,\pi X\,dx
=
\int_\Omega u(\mathrm{div}\,\pi X\ast\varphi_\varepsilon)\,dx
\nonumber\\
&
=
\int_\Omega u\,\mathrm{div}(\varphi_\varepsilon*\pi X)\,dx
=
-\int_\Omega\langle\varphi_\varepsilon*\pi X,\,dDu\rangle
\nonumber\\
&
=
-\int_U\langle\pi(\varphi_\varepsilon*Du),X\rangle\,dx.\nonumber
\end{align}
The result now follows by the fundamental lemma of the calculus of variations and smoothness of the mollified functions. 
\qed

\smallskip

\begin{proposition}\label{convergence_of_convolution}
Let $\Omega$ be an open set in $\mathbb{R}^2_0$. Let $u\in\mathrm{BV}(\Omega)$ and $U\subset\subset\Omega$ open with the property that $|D u|(\partial U)=0$. Then
\[
|D_\omega u|(U)=
\lim_{\varepsilon\downarrow 0}|D_\omega(u*\varphi_\varepsilon)|(U)
\]
where $(\varphi_\varepsilon)_{\varepsilon>0}$ is a family of mollifiers.
\end{proposition}

\smallskip

\noindent{\em Proof.}
By lower semicontinuity of the spherical variation (Lemma \ref{properties_of_spherical_variation}),
\[
|D_\omega u|(U)\leq\liminf_{\varepsilon\downarrow 0}|D_\omega u_\varepsilon|(U).
\]
The bulk of the proof is dedicated to the reverse inequality. Let  $0<\varepsilon<\mathrm{dist}(U,\partial\Omega)$. By Lemma \ref{properties_of_spherical_variation}, Lemma \ref{spherical_gradient_of_mollified_function} and Fubini's theorem,
\[
|D_\omega u_\varepsilon|(U) = \int_U|\nabla_\omega (u*\varphi_\varepsilon)|\,dx
= \int_U|\pi(Du*\varphi_\varepsilon)|\,dx
\]
\[
\leq
\int_U\int_{\Omega}\varphi_\varepsilon(x-y)|\pi(x)\rho(y)|\,d|Du|(y)\,dx
=
\int_{\Omega}\int_U\varphi_\varepsilon(x-y)|\pi(x)\rho(y)|\,dx\,d|Du|(y)
\]
for small $\varepsilon>0$. Denoting the integrand by $\phi_\varepsilon$ we may write the right-hand side of the last identity as
\[
\int_{\Omega}\phi_\varepsilon\,d|Du|
=
\int_U\phi_\varepsilon\,d|Du|
+
\int_{\Omega\setminus\overline{U}}\phi_\varepsilon\,d|Du|
\]
in virtue of the fact that $|D u|(\partial U)=0$ by assumption.
The mapping $x\mapsto\pi(x)v$ is continuous on $\mathbb{R}^2_0$ for any $v\in\mathbb{R}^2$. So for $y\in\Omega$,
\[
\lim_{\varepsilon\downarrow 0}\phi_\varepsilon(y)
=
\left\{
\begin{array}{lcl}
|\pi(y)\rho(y)| & \text{ if } & y\in U,\\
0              & \text{ if } & y\in\Omega\setminus\overline{U}.
\end{array}
\right.
\]
By the dominated convergence theorem,
\[
\limsup_{\varepsilon\downarrow 0}|D_\omega u_\varepsilon|(U)
\leq\int_U|\pi\rho|\,d|Du|
=|D_\omega u|(U)
\]
by Lemma \ref{properties_of_spherical_variation}. That is, the reverse inequality does indeed hold.
\qed

\smallskip

\noindent The next proposition is the spherical variation counterpart to \cite{Ambrosio2000} Proposition 3.103.

\begin{proposition}\label{tangential_variation}
 Let $\Omega$ be an open set in $\mathbb{R}^2_0$ and $u\in\mathrm{BV}_{\mathrm{loc}}(\Omega)$. Then $u_\tau\in\mathrm{BV}_{\mathbb{S}_\tau}(\Omega_\tau)$ for a.e. $\tau>0$ and 
\[
V_\omega(u,\Omega) = \int_0^\infty V_{\mathbb{S}_\tau}(u_\tau,\Omega_\tau)\,d\tau.
\]
\end{proposition}

\smallskip

\noindent{\em Proof.}
First note that equality holds for $u\in C^1(\Omega)\cap\mathrm{BV}(\Omega)$. Indeed by Lemma \ref{properties_of_spherical_variation} {\em (iv)}, switching to polar coordinates and Lemma \ref{properties_of_variation_on_sphere} {\em (iii)},
\begin{align}
|D_\omega u|(\Omega)
&=\int_\Omega|\nabla_\omega u|\,dx
=\int_0^\infty\int_{\Omega_\tau}|\nabla_\omega u|\,d\mathscr{H}^{1}\,d\tau\nonumber\\
&=\int_0^\infty\int_{\Omega_\tau}|\nabla_{\mathbb{S}_\tau} u_\tau|\,d\mathscr{H}^{1}\,d\tau
=\int_0^\infty V_{\mathbb{S}_\tau}(u_\tau,\Omega_\tau)\,d\tau\label{disintegration_of_spherical_variation}
\end{align}
as required.

\smallskip

\noindent By converting to polar coordinates it can be seen that $u_\tau\in L^1(\Omega_\tau,\mathscr{H}^{1})$ for a.e. $\tau>0$. For small $t>0$ put
\[
\Omega^t:=\Big\{x\in\Omega:\,|x|<1/t\text{ and }\mathrm{dist}(x,\partial\Omega)>t\Big\}.
\]
It holds that $|D u|(\partial\Omega^t)=0$ for all but countably many $t>0$ as $|D u|$ is a Radon measure on $\Omega$. For $\varepsilon>0$ small,
\[
\int_{\Omega^t}|u*\varphi_\varepsilon- u|\,dx
=
\int_0^\infty\|(u*\varphi_\varepsilon)_\tau - u_\tau\|_{L^1((\Omega^t)_\tau)}\,d\tau
\]
upon converting to polar coordinates and the left-hand side converges to zero in the limit $\varepsilon\downarrow 0$ by a standard property of mollification. So there exists a subsequence $\varepsilon_h\downarrow 0$ such that with $v_h=u*\varphi_{\varepsilon_h}\in C^\infty_c(\Omega)$,
\[
(v_h)_\tau\rightarrow u_\tau\text{ in }L^1((\Omega^t)_\tau,\,\mathscr{H}^{1})\text{ for a.e. }\tau>0.
\]
By lower semicontinuity Lemma \ref{variation_on_the_sphere}, Fatou's lemma, (\ref{disintegration_of_spherical_variation}) and Proposition \ref{convergence_of_convolution},
\begin{align}
\int_0^\infty V_{\mathbb{S}_\tau}(u_\tau,(\Omega^t)_\tau)\,d\tau
& \leq
\int_0^\infty\left\{\liminf_{h\rightarrow\infty}V_{\mathbb{S}_\tau}((v_h)_\tau,(\Omega^t)_\tau)\right\}
\,d\tau\nonumber\\
& \leq
\liminf_{h\rightarrow\infty}\int_0^\infty V_{\mathbb{S}_\tau}((v_h)_\tau,(\Omega^t)_\tau)
\,d\tau\nonumber\\
&=
\liminf_{h\rightarrow\infty}|D_\omega v_h|(\Omega^t)
=
|D_\omega u|(\Omega^t)
\leq 
V_\omega(u,\Omega).\nonumber
\end{align}
By the monotone convergence theorem,
\[
\int_0^\infty V_{\mathbb{S}_\tau}(u_\tau,\Omega_\tau)\,d\tau
\leq 
V_\omega(u,\Omega)
\]
so that $u_\tau\in\mathrm{BV}_{\mathbb{S}_\tau}(\Omega_\tau)$ for a.e. $\tau>0$.

\smallskip

\noindent As for the reverse inequality, let $X\in C^1_c(\Omega,\mathbb{R}^2)$ with $\|X\|_\infty\leq 1$. In polar coordinates
\[
\int_\Omega u\,\mathrm{div}(\pi X)\,dx
= 
\int_{p(\Omega)}\int_{\Omega_\tau}u_\tau\mathrm{div}_{\mathbb{S}_\tau}(\pi X)\,d\mathscr{H}^{1}
\,d\tau
\leq
\int_{p(\Omega)}V_{\mathbb{S}_\tau}(u_\tau,\Omega_\tau)\,d\tau
\]
recalling (\ref{variation_on_the_sphere}). Taking the supremum yields the reverse inequality.
\qed

\bigskip

\noindent Let $E$ be a set of locally finite perimeter in $\mathbb{R}^2_0$. In line with Lemma \ref{properties_of_spherical_variation} {\em (iii)},
\begin{equation}\label{spherical_variation_measure}
D_\omega\chi_E=\nu_E^\omega\,\mathscr{H}^{1}\mres\mathscr{F}E
\text{ and }
|D_\omega\chi_E|=|\nu_E^\omega|\,\mathscr{H}^{1}\mres\mathscr{F}E,
\end{equation}
where we write $\nu_E^\omega=\pi\nu_E$ as before. With the foregoing preparation we arrive at the main theorem of this Section. It is a tangential counterpart to \cite{Barchiesietal2013} Theorem 2.4.

\smallskip

\begin{theorem}\label{tangentialVolperttheorem}
Let $E$ be a set of locally finite perimeter in $\mathbb{R}^2_0$. Then for a.e. $\tau>0$,
\begin{itemize}
\item[(i)]
$E_\tau$ is a set of finite perimeter in $\mathbb{S}_\tau$;
\item[(ii)]
$\mathscr{H}^0(\mathscr{F}_{\mathbb{S}_\tau}E_\tau\Delta(\mathscr{F}E)_\tau)=0$.
\end{itemize}
\end{theorem}

\smallskip

\noindent{\em Proof.}
{\em (i)} This follows from Proposition \ref{tangential_variation}. {\em (ii)} By this last as well as Lemma \ref{properties_of_spherical_variation} and Lemma \ref{properties_of_variation_on_sphere}, 
\begin{equation}\label{E_r_and_finite_perimeter}
|D_\omega\chi_E|(\Omega)
=V_\omega(\chi_E,\Omega)
=\int_0^\infty V_{\mathbb{S}_\tau}((\chi_E)_\tau,\Omega_\tau)\,d\tau
=\int_0^\infty |D_{\mathbb{S}_\tau}\chi_{E_\tau}|(\Omega_\tau)\,d\tau 
\end{equation}
for any relatively compact open set $\Omega$ contained in $\mathbb{R}^2_0$. Let $A$ stand for the open centred annulus with radii $0<r<R<\infty$. Put
\[
\mathscr{M}:=\Big\{B\subset A:\,B\text{ is a Borel set and }
|D_\omega\chi_E|(B)
=\int_0^\infty |D_{\mathbb{S}_\tau}\chi_{E_\tau}|(B_\tau)\,d\tau \Big\}.
\]
By (\ref{E_r_and_finite_perimeter}) the collection $\mathscr{M}$ contains all open sets in $A$; it also has the properties
\begin{itemize}
\item[(a)]
$(B_h)_{h\in\mathbb{N}}\subset\mathscr{M},\,B_h\uparrow B\Rightarrow B\in\mathscr{M}$;
\item[(b)]
$B,\,B^\prime,\,B\cup B^\prime\in\mathscr{M}\Rightarrow B\cap B^\prime\in\mathscr{M}$;
\item[(c)]
$B\in\mathscr{M}\Rightarrow A\setminus B\in\mathscr{M}$.
\end{itemize}
$\mathscr{M}$ therefore coincides with the $\sigma$-algebra of Borel sets in $A$ by \cite{Ambrosio2000} Remark 1.9. This entails that the identity in (\ref{E_r_and_finite_perimeter}) holds for each relatively compact Borel set $B$ in $\mathbb{R}^2_0$. 

\smallskip

\noindent Let $B$ be a relatively compact Borel set in $\mathbb{R}^2_0$. By Proposition \ref{tangential_integral_equality}, (\ref{spherical_variation_measure}), the last observation, and Theorem \ref{deGiorgivariationsphere},
\begin{align}
\int_0^\infty\int_{(\mathscr{F}E)_\tau}\chi_B\,d\mathscr{H}^{0}\,d\tau
&=
\int_{\mathscr{F}E}\chi_B|\nu_E^\omega|\,d\mathscr{H}^{1}
=
|D_\omega\chi_E|(B)\nonumber\\
&=\int_0^\infty|D_{\mathbb{S}_\tau}\chi_{E_\tau}|(B_\tau)\,d\tau
=
\int_0^\infty\int_{\mathscr{F}_{\mathbb{S}_\tau}E_\tau}\chi_B\,d\mathscr{H}^{0}\,d\tau.\nonumber
\end{align}
Let $\mathscr{G}$ be a countable base for the $\sigma$-algebra of Borel sets in $\mathbb{S}$. Fix $G\in\mathscr{G}$ and any Borel set $I\subset(0,\infty)$. This last identity holds with $B=\{\tau\omega:\tau\in I\text{ and }\omega\in G\}$. We derive that 
\begin{equation}\label{integral_identity_on_shell}
\int_{(\mathscr{F}E)_\tau}\chi_G\circ I^{1/\tau}\,d\mathscr{H}^{0}
=
\int_{\mathscr{F}_{\mathbb{S}_\tau}E_\tau}\chi_G\circ I^{1/\tau}\,d\mathscr{H}^{0}
\end{equation}
for a.e. $\tau>0$. Identity (\ref{integral_identity_on_shell}) therefore holds for every Borel set $G$ in $\mathbb{S}$ for a.e. $\tau$ in $(0,\infty)$. Item {\em (ii)} now follows. 
\qed

\section{A radial Vol'pert Theorem}\label{A_radial_Volpert_Theorem}

\smallskip

\noindent The goal of this Section is to prove a radial Vol'pert theorem in Theorem \ref{radialVolperttheorem}; there are similarities with the previous Section. For $x\in\mathbb{R}^2_0$ and $v\in\mathbb{R}^2$ set $(\tilde{\pi} v)(x):=\langle v,\omega\rangle\omega$ where $\omega=x/|x|$. That is, $(\tilde{\pi} v)(x)$ stands for the radial component of $v$ at the point $x$. Let $\Omega$ be a relatively compact open set in $\mathbb{R}^2_0$. The radial variation of $u\in L^1(\Omega)$ is defined by
\[
V_{\mathrm{rad}}(u,\Omega):=\sup\Big\{\int_{\Omega}u\,\mathrm{div}(r^{-1}\tilde{\pi} X)\,dx:X\in C^1_c(\Omega,\mathbb{R}^2)\text{ and }\|X\|_\infty\leq 1\Big\}.
\]

\smallskip

\begin{lemma}\label{properties_of_radial_variation}
Let $\Omega$ be a relatively compact open set in $\mathbb{R}^2_0$. Then
\begin{itemize}
\item[(i)] the functional $u\mapsto V_\mathrm{rad}(u,\Omega)$ is lower semicontinuous with respect to convergence in $L^1(\Omega)$.
\end{itemize}
Let $u\in\mathrm{BV}(\Omega)$. Then
\begin{itemize}
\item[(ii)]
there exists a unique $\mathbb{R}^2$-valued Radon measure $D_\mathrm{rad} u$ on $\Omega$ such that
\[
\int_{\Omega}u\,\mathrm{div}(r^{-1}\tilde{\pi} X)\,dx
=
-\int_\Omega\langle X,dD_{\mathrm{rad}}u\rangle
\]
for any $X\in C^1_c(\Omega,\mathbb{R}^2)$ and $|D_\mathrm{rad} u|(\Omega)=V_\mathrm{rad}(u,\Omega)$;
\item[(iii)]
$D_{\mathrm{rad}}u=\frac{1}{r}(\tilde{\pi}\varrho)|Du|$ and $|D_{\mathrm{rad}}u|=(1/r)|\langle\varrho,\omega\rangle||Du|$ where $Du=\rho\,|Du|$ is the polar decomposition of $D u$ with $\rho$ a unique $\mathbb{S}$-valued function in $L^1(\Omega,\,|Du|)^2$;
\item[(iv)] $D_{\mathrm{rad}}u=\frac{1}{r}(\partial_r u)\mathscr{L}^2$ and $|D_{\mathrm{rad}}u|=\frac{1}{r}|\partial_r u|\mathscr{L}^2$ for $u\in C^1(\Omega)\cap\mathrm{BV}(\Omega)$.
\end{itemize}
\end{lemma}

\smallskip

\noindent{\em Proof.}
{\em (i)} follows as in \cite{Ambrosio2000} Remark 3.5. {\em (ii)} For $X\in C^1_c(\Omega,\mathbb{R}^2)$ the vector field $r^{-1}\tilde{\pi} X$ belongs to the same class with uniform norm bounded by $\|X\|_\infty/r_{\mathrm{min}}$ where $r_\mathrm{min}:=\inf\{|x|:x\in\Omega\}\in(0,\infty)$. This means that $V_{\mathrm{rad}}(u,\Omega)$ is finite as $u\in\mathrm{BV}(\Omega)$. The assertion now follows from Riesz's Theorem (cf. \cite{Ambrosio2000} Theorem 1.54). The fact that $u$ belongs to the class $\mathrm{BV}(\Omega)$ leads to {\em (iii)} and {\em (iv)}. 
\qed

\smallskip

\begin{proposition}\label{convergence_of_convolution_of_radial_variation}
Let $\Omega$ be a relatively compact open set in $\mathbb{R}^2_0$ and $u\in\mathrm{BV}(\Omega)$. Let $U\subset\subset\Omega$ open with the property that $|D u|(\partial U)=0$. Then
\[
|D_\mathrm{rad} u|(U)=
\lim_{\varepsilon\downarrow 0}|D_\mathrm{rad}(u*\varphi_\varepsilon)|(U)
\]
where $(\varphi_\varepsilon)_{\varepsilon>0}$ is a family of mollifiers.
\end{proposition}

\smallskip

\noindent{\em Proof.}
By the lower semicontinuity of the radial variation in Lemma \ref{properties_of_radial_variation}, 
\[
|D_\mathrm{rad} u|(U)\leq\liminf_{\varepsilon\downarrow 0}|D_\mathrm{rad}(u*\varphi_\varepsilon)|(U).
\]
As for the reverse inequality let us first choose $0<\varepsilon<\mathrm{dist}(U,\partial\Omega)$. By Lemma \ref{properties_of_radial_variation}, 
\begin{align}
|D_\mathrm{rad}(u*\varphi_\varepsilon)|(U)&=\int_U\frac{1}{|x|}|\langle\nabla(u*\varphi_\varepsilon),\omega\rangle|\,dx.\nonumber
\end{align}
As $u\in\mathrm{BV}(\Omega)$,
\begin{align}
\langle\nabla(u*\varphi_\varepsilon)(x),\omega\rangle
&=\int_\Omega\langle(\nabla\varphi_\varepsilon)(x-y),\frac{x}{|x|}\rangle u(y)\,dy
=-\int_\Omega\mathrm{div}_y\Big(\varphi_\varepsilon(x-\cdot)\frac{x}{|x|}\Big)u(y)\,dy
\nonumber\\
&=\int_\Omega\langle\varphi_\varepsilon(x-\cdot)\frac{x}{|x|},dDu\rangle
=\int_\Omega\langle\varphi_\varepsilon(x-\cdot)\frac{x}{|x|},\varrho\rangle\,d|Du|\nonumber
\end{align}
and so
\begin{align}
|\langle\nabla(u*\varphi_\varepsilon)(x),\omega\rangle|
&\leq
\int_\Omega\varphi_\varepsilon(x-\cdot)|\langle\frac{x}{|x|},\varrho\rangle|\,d|Du|
\nonumber
\end{align}
Inserting this estimate into the above equality and using Tonelli's Theorem leads to 
\begin{align}
|D_\mathrm{rad}(u*\varphi_\varepsilon)|(U)
&\leq\int_U\frac{1}{|x|}\int_\Omega\varphi_\varepsilon(x-\cdot)|\langle\frac{x}{|x|},\varrho\rangle|\,d|Du|\,dx\nonumber\\
&=\int_\Omega\int_U\frac{1}{|x|}\varphi_\varepsilon(x-\cdot)|\langle\frac{x}{|x|},\varrho\rangle|\,dx\,d|Du|
=\int_\Omega\phi_\varepsilon\,d|Du|
\nonumber
\end{align}
in an obvious notation. We may write the last expression as
\[
\int_{\Omega}\phi_\varepsilon\,d|Du|
=
\int_U\phi_\varepsilon\,d|Du|
+
\int_{\Omega\setminus\overline{U}}\phi_\varepsilon\,d|Du|
\]
in virtue of the fact that $|D u|(\partial U)=0$ by assumption. For $y\in\Omega\setminus\partial U$,
\[
\lim_{\varepsilon\downarrow 0}\phi_\varepsilon(y)
=
\left\{
\begin{array}{ll}
\frac{1}{|y|}|\langle\frac{y}{|y|},\varrho(y)\rangle| & \text{ if }y\in U,\\
0              & \text{ if }y\in\Omega\setminus\overline{U}.
\end{array}
\right.
\]
By the dominated convergence theorem,
\[
\limsup_{\varepsilon\downarrow 0}|D_\mathrm{rad}(u*\varphi_\varepsilon)|(U)
\leq\int_U\frac{1}{r}|\langle\omega,\varrho\rangle|\,d|Du|
=|D_\mathrm{rad} u|(U).
\]
This proves the result. 
\qed

\smallskip

\noindent For a Borel set $E$ in $\mathbb{R}^2_0$ and $\omega\in\mathbb{S}$ the $\omega$-section of $E$ is defined to be the intersection of $E$ with the open ray $(0,\infty)\omega$ in direction $\omega$. The notation $u_\omega$ refers to the restriction of the function $u$ to $(0,\infty)\omega$.

\smallskip

\begin{proposition}\label{expression_for_radial_variation}
Let $\Omega$ be a relatively compact open set in $\mathbb{R}^2_0$ and $u\in\mathrm{BV}(\Omega)$. Then $u_\omega\in\mathrm{BV}(\Omega_\omega)$ for a.e. $\omega\in\mathbb{S}$ and
\begin{align}
V_{\mathrm{rad}}(u,\Omega)&=\int_{\mathbb{S}}V(u_\omega,\Omega_\omega)\,d\mathscr{H}^1.
\label{V_rad_equality}
\end{align}
\end{proposition}

\smallskip

\noindent{\em Proof.}
Let $u\in C^1_c(\Omega)$. By Lemma \ref{properties_of_radial_variation},
\begin{align}
V_{\mathrm{rad}}(u,\Omega)&=|D_{\mathrm{rad}}u|(\Omega)
=\int_{\Omega}\frac{1}{r}|\partial_r u|\,dx
=\int_{\mathbb{S}}\int_{\Omega_\omega}|\partial_\tau u|\,d\tau\,d\mathscr{H}^1
=\int_{\mathbb{S}}V(u_\omega,\Omega_\omega)\,d\mathscr{H}^1\label{V_rad_identity}
\end{align}
upon switching to polar coordinates.

\smallskip

\noindent Now let $u\in\mathrm{BV}(\Omega)$.  It can be seen that $u_\omega\in L^1(\Omega_\omega,\mathscr{H}^1)$ for a.e. $\omega\in\mathbb{S}$ by converting to polar coordinates. For small $t>0$ put
\[
\Omega^t:=\{x\in\Omega:\mathrm{dist}(x,\partial\Omega)>t \}.
\]
It holds that $|D u|(\partial\Omega^t)=0$ for all but countably many $t>0$ as $|D u|$ is a finite measure on $\Omega$. For $\varepsilon>0$ small,
\[
\int_{\Omega^t}|u*\varphi_\varepsilon- u|\,\frac{dx}{|x|}
=
\int_\mathbb{S}\|(u*\varphi_\varepsilon)_\omega - u_\omega\|_{L^1((\Omega^t)_\omega)}\,d\mathscr{H}^1
\]
and the left-hand side converges to zero in the limit $\varepsilon\downarrow 0$ by a standard property of mollification. So there exists a subsequence $\varepsilon_h\downarrow 0$ such that with $v_h=u*\varphi_{\varepsilon_h}\in C^\infty_c(\Omega)$,
\[
(v_h)_\omega\rightarrow u_\omega\text{ in }L^1((\Omega^t)_\omega,\,\mathscr{H}^1)\text{ for a.e. }\omega\in\mathbb{S}.
\]
By lower semicontinuity, Fatou's lemma, (\ref{V_rad_identity}) and Proposition \ref{convergence_of_convolution_of_radial_variation},
\begin{align}
\int_{\mathbb{S}}V(u_\omega,(\Omega^t)_\omega)\,d\mathscr{H}^1
& \leq
\int_{\mathbb{S}}\left\{\liminf_{h\rightarrow\infty}V((v_h)_\omega,(\Omega^t)_\omega)\right\}\,d\mathscr{H}^1 
\leq
\liminf_{h\rightarrow\infty}\int_{\mathbb{S}}V((v_h)_\omega,\,(\Omega^t)_\omega)\,d\mathscr{H}^1\nonumber\\
&=
\liminf_{h\rightarrow\infty}V_{\mathrm{rad}}(v_h,\Omega^t)
=V_{\mathrm{rad}}(u,\Omega^t)
\leq 
V_{\mathrm{rad}}(u,\Omega).\nonumber
\end{align}
By the monotone convergence theorem,
\[
\int_{\mathbb{S}}V(u_\omega,\Omega_\omega)\,d\mathscr{H}^1
\leq
V_{\mathrm{rad}}(u,\Omega)
\]
and this permits us to draw the conclusion that $u_\omega\in\mathrm{BV}(\Omega_\omega)$ for a.e. $\omega\in\mathbb{S}$.

\smallskip

\noindent For the reverse inequality, let $X\in C^1_c(\Omega,\mathbb{R}^2)$ with $\|X\|_\infty\leq 1$. In polar coordinates
\[
\int_\Omega u\,\mathrm{div}(r^{-1}\tilde{\pi} X)\,dx
= 
\int_{\mathbb{S}}\int_{\Omega_\omega}u_\omega r\,\mathrm{div}(r^{-1}\tilde{\pi} X)\,dr\,d\mathscr{H}^1
\leq
\int_{\mathbb{S}}V(u_\omega,\Omega_\omega)\,d\mathscr{H}^1
\]
as $\mathrm{div}(r^{-1}\tilde{\pi} X)=r^{-1}\partial_r\langle X,\omega\rangle$ which completes the proof after taking the supremum.
\qed

\smallskip

\begin{proposition}\label{integral_equality}
Let $E$ be a set of locally finite perimeter in $\mathbb{R}^2_0$ and let $g:\mathbb{R}^2_0\rightarrow[0,\infty]$ be a Borel function. Then
\[
\int_{\mathscr{F}E}g(1/r)|\langle\nu_E,\omega\rangle|\,d\mathscr{H}^1=\int_{\mathbb{S}}\int_{(\mathscr{F}E)_\omega}g\,d\mathscr{H}^0\,d\mathscr{H}^1.
\]
\end{proposition}

\smallskip

\noindent{\em Proof.}
Define $f:\mathbb{R}^2_0\rightarrow\mathbb{S}$ via $x\mapsto x/|x|$. Let $x\in\mathbb{R}^2_0$ and $M$ be a line through $x$ in direction $v\in\mathbb{S}$. Then 
\[
J_1d^M f_x=(1/r)|\langle\omega,v^\perp\rangle|.
\]
Now appeal to the generalised area formula \cite{Ambrosio2000} Theorem 2.91.
\qed

\smallskip

\noindent The following theorem is a Vol'pert-type theorem; see \cite{Barchiesietal2013} Theorem 2.4.

\smallskip

\begin{theorem}\label{radialVolperttheorem}
Let $E$ be a set of locally finite perimeter in $\mathbb{R}^2_0$. Then for a.e. $\omega\in\mathbb{S}$,
\begin{itemize}
\item[(i)]
$E_\omega$ is a set of locally finite perimeter in $(\mathbb{R}^2)_\omega$;
\item[(ii)]
$\mathscr{H}^0(\mathscr{F}E_\omega\Delta(\mathscr{F}E)_\omega)=0$.
\end{itemize}
\end{theorem}

\smallskip

\noindent{\em Proof.}
{\em (i)} This follows from Proposition \ref{expression_for_radial_variation}. {\em (ii)} By Proposition \ref{expression_for_radial_variation} and Lemma \ref{properties_of_radial_variation}, 
\begin{equation}\label{total_radial_variation_identity}
|D_\mathrm{rad}\chi_E|(\Omega)
=V_\mathrm{rad}(\chi_E,\Omega)
=\int_{\mathbb{S}}V((\chi_E)_\omega,\Omega_\omega)\,d\mathscr{H}^1
=\int_{\mathbb{S}}|D(\chi_E)_\omega|(\Omega_\omega)\,d\mathscr{H}^1
\end{equation}
for any relatively compact open set $\Omega\subset\mathbb{R}^2_0$. Let $A$ stand for the open centred annulus with radii $0<r<R<\infty$. Put
\[
\mathscr{M}:=\Big\{B\subset A:\,B\text{ is a Borel set and }
|D_\mathrm{rad}\chi_E|(B)
=\int_{\mathbb{S}}|D(\chi_E)_\omega|(B_\omega)\,d\mathscr{H}^1\Big\}.
\]
By (\ref{total_radial_variation_identity}) the collection $\mathscr{M}$ contains all open sets in $A$; it also has the properties
\begin{itemize}
\item[(a)]
$(B_h)_{h\in\mathbb{N}}\subset\mathscr{M},\,B_h\uparrow B\Rightarrow B\in\mathscr{M}$;
\item[(b)]
$B,\,B^\prime,\,B\cup B^\prime\in\mathscr{M}\Rightarrow B\cap B^\prime\in\mathscr{M}$;
\item[(c)]
$B\in\mathscr{M}\Rightarrow A\setminus B\in\mathscr{M}$.
\end{itemize}
$\mathscr{M}$ therefore coincides with the $\sigma$-algebra of Borel sets in $A$ by \cite{Ambrosio2000} Remark 1.9. This entails that the identity holds for each relatively compact Borel set $B$ in $\mathbb{R}^2_0$. For such a Borel set $B$ in $\mathbb{R}^2_0$,
\begin{align}
\int_{\mathbb{S}}\int_{(\mathscr{F}E)_\omega}\chi_B\,d\mathscr{H}^0\,d\mathscr{H}^1&=
\int_{\mathscr{F}E}(1/r)\chi_B|\langle\nu_E,\omega\rangle|\,d\mathscr{H}^1
=|D_{\mathrm{rad}}\chi_E|(B)
\nonumber\\
&=\int_{\mathbb{S}}|D(\chi_E)_\omega|(B_\omega)\,d\mathscr{H}^1
=\int_{\mathbb{S}}\int_{\mathscr{F}E_\omega}\chi_B\,d\mathscr{H}^0\,d\mathscr{H}^1\nonumber
\end{align}
by Proposition \ref{integral_equality} and the above identity. This leads to {\em (ii)}.
\qed

\section{Further preliminary results}\label{Further_preliminary_results}

\smallskip

\noindent In this Section we establish two miscellaneous results for use in the next Section. The first is a standard result which characterises the image of a set of locally finite perimeter under a proper $C^1$ diffeomorphism. The second result in Lemma \ref{on_toplogical_boundary} describes a weak regularity property of the boundary of an isoperimetric minimiser. This is used in the proof of the boundedness theorem Theorem \ref{boundedness_theorem}. It is not as deep as the regularity results treated in \cite{Tamanini1984} or \cite{Maggi2012} Part III for example and we include its proof for completeness. 

\smallskip

\begin{lemma}\label{diffeomorhism_invariance_of_boundary}
Let $\Omega$ and $\Omega^\prime$ be open sets in $\mathbb{R}^2$ and $\Phi:\Omega\rightarrow\Omega^\prime$ a proper $C^1$ diffeomeorphism. Let $E$ be a set of locally finite perimeter in $\Omega$ and put $F:=\Phi(E)$. Then 
\begin{itemize}
\item[(i)] $F$ is a set of locally finite perimeter in $\Omega^\prime$;
\item[(ii)] $(\partial^* F)\cap\Omega^\prime=\Phi(\partial^* E\cap\Omega)$;
\item[(iii)] $\mathscr{F}F\cap\Omega^\prime$ is $\mathscr{H}^1$-equivalent to
$\Phi(\mathscr{F}E\cap\Omega)$.
\end{itemize}
\end{lemma}

\smallskip

\noindent{\em Proof.}
Part {\em (i)} follows from \cite{Ambrosio2000} Theorem 3.16. {\em (ii)} We first claim that for each $x\in E^0\cap\Omega$ the point $y:=\Phi(x)$ belongs to $F^0\cap\Omega^\prime$. Choose open sets $\Omega_1^\prime$ and $\Omega_2^\prime$ in $\Omega^\prime$ with $\Omega_1^\prime\subset\subset\Omega_2^\prime\subset\subset\Omega^\prime$.  Their inverse images in $\Omega$ are denoted $\Omega_1$ and $\Omega_2$. Let $x\in\Omega_1$. Let $K$ stand for the Lipschitz constant of the restriction of $\Phi^{-1}$ to $\Omega_2^\prime$. Then $B(y,r)\subset\Phi(B(x,Kr))$ for $r>0$ small.
By \cite{Ambrosio2000} Theorem 2.53, the fact that $\Phi$ is a bijection, and \cite{Ambrosio2000} Proposition 2.49,
\begin{align}
|F\cap B(y,r)|&=\mathscr{H}^2(F\cap B(y,r))
\leq\mathscr{H}^2(F\cap\Phi(B(x,Kr))
=\mathscr{H}^2(\Phi(E\cap B(x,Kr))\nonumber\\
&\leq
L^2\mathscr{H}^2(E\cap B(x,Kr))
=L^2|E\cap B(x,Kr)|\nonumber
\end{align}
where $L$ stands for the Lipschitz constant of $\Phi$ restricted to the set $\Phi^{-1}(\Omega_2^\prime)$ noting that this latter set is compactly embedded in $\Omega$ as the mapping $\Phi$ is proper. Noting that $B(x,r)\subset\Phi^{-1}(B(y,Lr))$ a similar argument shows that $|B(x,r/L)|\leq K^2|B(y,r)|$ for $r$ small. This means that
\[
\frac{|F\cap B(y,r)|}{|B(y,r)|}
\leq(KL)^2\frac{|E\cap B(x,Kr)|}{|B(x,r/L)|}
=(KL)^4\frac{|E\cap B(x,Kr)|}{|B(x,Kr)|}.
\]
This proves the claim. This leads to the identity $\Phi(E^0)=\Phi(E)^0$. The observation that $E^1=(\Omega\setminus E)^0$ establishes {\em (ii)}. Item {\em (iii)} follows from Federer's Theorem \cite{Ambrosio2000} Theorem 3.61 and {\em (ii)}. 
\qed

\medskip
 


\noindent Let $\Omega$ be an open set in $\mathbb{R}^2$. Given positive locally Lipschitz densities $f$ and $g$ on $\Omega$ we shall refer to the variational problem
\begin{equation}\label{fg_isoperimetric_problem}
I(v):=\inf\Big\{P_g(E):E\text{ is a set of locally finite perimeter in }\Omega\text{ and }V_f(E)=v\Big\}
\end{equation}
for $v>0$. In the next result we establish mild regularity of the boundary of an isoperimetric minimiser of (\ref{fg_isoperimetric_problem}).

\smallskip

\noindent Let $E$ be a set of locally finite perimeter in $\Omega$. According to \cite{Maggi2012} Proposition 12.19 there exists a Borel set $F$ equivalent to $E$ with the property that the topological boundary $\partial F\cap U$ in $\Omega$ satisfies the condition 
\begin{align}
\Big\{x\in\Omega:0<|F\cap B(x,r)\cap\Omega|<\pi r^2\text{ for all }r>0\Big\}&=\partial F\cap\Omega.\label{topological_boundary}
\end{align}

\smallskip

\begin{lemma}\label{on_toplogical_boundary}
Let $\Omega$ be an open set in $\mathbb{R}^2$ and assume that $f$ and $g$ are positive locally Lipschitz densities on $\Omega$. Let $v>0$ and suppose that the set $E$ is a minimiser of (\ref{fg_isoperimetric_problem}). Assume that $E$ satisfies the condition in (\ref{topological_boundary}). Then for each $x\in\partial E\cap\Omega$ there exists $\delta>0$ and a constant $0<c<1/2$ such that
\[
c\leq\frac{|E\cap B(x,r)|}{\pi r^2}\leq 1-c
\text{ for each }0<r<\delta;
\]
and in particular $\mathscr{H}^1(\partial E\cap\Omega\setminus\mathscr{F}E)=0$.
\end{lemma}

\smallskip 

\noindent{\em Proof.}
Let $x\in\Omega$ and $t>0$ such that $B(x,t)\subset\subset\Omega$. An adaptation of the argument in \cite{McGillivray2018_1} Proposition 3.5 (and \cite{Ambrosio2000} Corollary 3.89) leads to the formulae
\begin{align}
P_g(E)&=P_g(E,B(x,t))+\int_{\partial B(x,t)}g|\chi^{B(x,t)}_E-\chi^{\Omega\setminus\overline{B}(x,t)}|\,d\mathscr{H}^1+P_g(E,\Omega\setminus\overline{B}(x,t))\label{g_perimeter_decomposed_1}
\end{align}
and
\begin{align}
P_g(E\cap B(x,t))&=P_g(E,B(x,t))+\int_{\partial B(x,t)}g\chi^{B(x,t)}_E\,d\mathscr{H}^1.\label{g_perimeter_decomposed_2}
\end{align}

\smallskip

\noindent Let $x\in\partial E\cap\Omega$. By \cite{McGillivray2018_1} Proposition 3.6 there exist constants $C>0$ and $\delta\in(0,d(x,\partial\Omega))$ with the following property. For any $0<r<\delta$,
\begin{equation}\label{Morgan_lemma_assertion_1}
P_g(E)-P_g(F)\leq C\big|V_f(E)-V_f(F)\big|
\end{equation}
where $F$ is any set with locally finite perimeter in $\Omega$ such that $E\Delta F\subset\subset B(x,r)$. Fix $0<r<\delta$ and choose $s$ with $0<r<s<\delta$. Put $F:=E\setminus B(x,r)$. By (\ref{Morgan_lemma_assertion_1}) and (\ref{g_perimeter_decomposed_1}),
\begin{align}
P_g(E,B(x,s))&\leq P_g(E\setminus B(x,r),B(x,s))+CV_f(E\cap B(x,r))\nonumber\\
&=P_g(E,B(x,s)\setminus\overline{B}(x,r))+\int_{\partial B(x,r)}g\chi^{\Omega\setminus\overline{B}(x,r)}_E\,d\mathscr{H}^1+CV_f(E\cap B(x,r)).\nonumber
\end{align}
Upon letting $s\downarrow r$ we obtain
\[
P_g(E,B(x,r))\leq\int_{\partial B(x,r)}g\chi^{\Omega\setminus\overline{B}(x,r)}_E\,d\mathscr{H}^1
+CV_f(E\cap B(x,r)).
\]
It follows from this last that
\[
P_g(E\cap B(x,r))\leq\int_{\partial B(x,r)}g\Big\{
\chi^{B(x,r)}_E+\chi^{\Omega\setminus\overline{B}(x,r)}_E\Big\}
+CV_f(E\cap B(x,r))
\]
with the help of (\ref{g_perimeter_decomposed_2}). From Federer's Theorem (\cite{Ambrosio2000} Theorem 3.60) and the tangential Vol'pert-type theorem (Theorem \ref{tangentialVolperttheorem}),
\[
\chi_E^{B(x,t)}=\chi_{E^1}\text{ }\mathscr{H}^1\text{-a.e. on }\partial B(x,t)
\]
for a.e. $t\in(0,\delta)$ and likewise for $\chi_E^{\Omega\setminus\overline{B}(x,t)}$. So in fact we may write
\begin{align}
P_g(E\cap B(x,t))&\leq 2\mathscr{H}^1_g((E^1)_t)+CV_f(E\cap B(x,t))\label{inequality_for_P_g}
\end{align}
for a.e. $t$ in $(0,\delta)$.

\smallskip

\noindent Put $m(t):=|E\cap B(x,t)|$ for $t\in(0,\delta)$. The function $m$ is absolutely continuous on $(0,\delta)$ and $m^\prime=\mathscr{H}^1(E^1\cap\partial B(x,t))$ a.e. on $(0,\delta)$. Let $0<c<1$. As both $f$ and $g$ are positive and continuous on $\Omega$ we may choose $0<\delta_1<\delta$ such that 
\[
(1-c)f(x)\leq f\leq(1+c)f(x)
\text{ and }
(1-c)g(x)\leq g\leq(1+c)g(x)
\text{ on }\overline{B}(x,\delta_1).
\]
From (\ref{inequality_for_P_g}) and the classical isoperimetric inequality we obtain
\[
\sqrt{4\pi}(1-c)g(x)\leq 2(1+c)g(x)\frac{m^\prime}{\sqrt{m}}+C(1+c)f(x)\sqrt{m}
\]
a.e. on $(0,\delta_1)$. This yields the estimate
\[
\frac{1-c}{2(1+c)}-\frac{Cf(x)}{8g(x)}t\leq\sqrt{\frac{m}{\pi t^2}}
\]
on $(0,\delta_1)$. This leads to the lower bound in the statement of the Lemma. The upper bound follows from the fact that the complement $\Omega\setminus E$ is a relative isoperimetric minimiser.

\smallskip

\noindent The inequality entails that $\partial E\cap\Omega$ is contained in the essential boundary $\partial^* E\cap\Omega$. This implies the last assertion of the Lemma by Federer's Theorem (see for example \cite{Ambrosio2000} Theorem 3.60).
\qed

\section{Existence and boundedness}\label{Existence_and_boundedness}

\smallskip

\noindent In this Section we turn to the topic of existence and boundedness of isoperimetric minimisers. We seek to adapt the existence result contained in \cite{MorganPratelli2011} Theorem 3.3 in the case of Euclidean space with density to the two-weighted case with radial volume and perimeter densities $f$ respectively $g$;  and we supply a variation of the boundedness result in \cite{MorganPratelli2011} Theorem 5.9. We mention \cite{DePhilippisetal2017} Theorem 1.2 in passing but we do not pursue this approach here. Our results are contained in Theorem \ref{existence_theorem} and Theorem \ref{boundedness_theorem}. They extend \cite{DiGiosiaetal2016} and \cite{DiGiosiaetal2019} Theorem 2.3 and Theorem 2.4. We require that the density relative to the metric $\psi$ diverges at infinity as in condition (A.4) below. In the case of radial power weights as considered in this paper this condition fails if $\alpha=2\beta$ with $\beta<0$. We then have recourse to \cite{MorganRitore2002} Theorem 1.1 (and \cite{BrayMorgan2001} Corollary 2.4) through a reparametrisation of the problem. 

\smallskip

\noindent We begin with this latter result. Let us introduce the index set
\begin{align}
\mathcal{P}^-&:=\Big\{(\alpha,\beta)\in\mathbb{R}^2:-2<\alpha\leq 0\text{ and }\alpha=2\beta\Big\}.\nonumber
\end{align}
For indices in $\mathcal{P}^-$ existence holds as we show shortly. A crucial tool in the proof is the mapping (\ref{Morgan_map}). This also plays a r\^ole in preparing the proof of Theorem \ref{existence_theorem} so we choose to introduce this material at an early stage. Towards the end of this Section we show that centred balls are uniquely isoperimetric in Theorem \ref{isoperimetry_on_P_-}. 

\smallskip

\begin{theorem}\label{isoperimetry_on_P_-}
Let $(\alpha,\beta)\in\mathcal{P}^-\setminus\{0\}$. Then centred balls are minimisers for the problem (\ref{isoperimetric_problem}) up to equivalence.
\end{theorem}

\smallskip

\noindent  We begin with a remark. Choose a pair $(\alpha,\beta)\in\mathcal{P}^-\setminus\{0\}$. The angle $\gamma\in(0,\pi/2)$ is characterised by the relation $\sin\gamma=\beta+1$. Define a (punctured) cone 
\begin{align}
C&=C_\beta:=\Big\{(x,z)\in\mathbb{R}^2_0\times\mathbb{R}:z=(\cot\gamma)|x|\Big\}\label{cone}
\end{align}
in $\mathbb{R}^3$. The mapping $\phi:\mathbb{R}^2_0\rightarrow C;x\mapsto|x|^\beta(x,(\cot\gamma)|x|)$ has first fundamental form
\[
ds^2=|x|^{2\beta}(dx_1^2+dx_2^2)
\]
and is a conformal parametrisation of $C$. In short, the punctured plane equipped with the conformal metric $ds=|z|^\beta|dz|$ is conformally equivalent to the punctured cone $C$.

\smallskip

\noindent Given an index $-1<\beta<0$ we make use of the bijection
\begin{align}
&\Phi:\mathbb{R}^2_0\rightarrow\mathbb{R}^2_0;x\mapsto \frac{1}{\beta+1}|x|^\beta x.\label{Morgan_map}
\end{align} 
This is the composition of the map $\phi$ above with the map denoted $f$ in the proof of \cite{MorganRitore2002} Theorem 1.1. Fix $x\in\mathbb{R}^2_0$ and consider the vector $v=v_\omega\omega+v_{\omega^\perp}\omega^\perp$ where $\omega$ and $\omega^\perp$ are unit vectors, $\omega$ lies in direction $x$ and the pair $\{\omega,\omega^\perp\}$ is positively oriented. A calculation gives
\[
d\Phi_x v=|x|^\beta\Big\{v_\omega\omega+\frac{|x|}{\beta+1}v_{\omega^\perp}\omega^\perp\Big\}
\]
from which it follows that
\begin{align}
J_2d\Phi_x&=\frac{|x|^{2\beta}}{\beta+1}\label{J_2_of_Phi}
\end{align}
and
\begin{align}
J_1d^{x+M}\Phi_x
&=\frac{|x|^{\beta}}{\beta+1}\Big\{(\beta+1)^2v_\omega^2+|x|^2 v_{\omega^\perp}^2\Big\}^{1/2}
\leq\frac{|x|^{\beta}}{\beta+1}\label{J_1_of_Phi}
\end{align}
where $M:=\mathrm{span}\{v\}$.

\smallskip

\noindent{\em Proof of Theorem \ref{isoperimetry_on_P_-}.}
Let $E$ be a set of locally finite perimeter in $\mathbb{R}^2_0$. By (\ref{J_2_of_Phi}) and the area formula (\cite{Ambrosio2000} Theorem 2.71 for example), $V_\alpha(E)=(\beta+1)V_0(\Phi(E))$. By (\ref{J_1_of_Phi}) and the generalised area formula (\cite{Ambrosio2000} Theorem 2.91 for example), $P_\beta(E)\geq(\beta+1)P_0(\Phi(E))$ with the help of Lemma \ref{diffeomorhism_invariance_of_boundary}. It follows that the centred ball is a minimiser for the problem (\ref{isoperimetric_problem}) as the inverse of $\Phi$ preserves centred circles. 
\qed

\bigskip

\noindent We now derive an existence and a boundedness theorem for the isoperimetric problem in the two-weighted situation (namely, Theorems \ref{existence_theorem} and \ref{boundedness_theorem}). Suppose that $\mathtt{f}$ and $\mathtt{g}$ are continuous positive functions defined on $(0,\infty)$. Put $\zeta:=\mathtt{f}/\mathtt{g}$ and $\psi:=\mathtt{g}^2/\mathtt{f}$ on $(0,\infty)$. Assume that
\begin{itemize}
\item[(A.1)] $\int_0^1 t\mathtt{f}\,dt<\infty$;
\item[(A.2)] $\int_1^{\infty} t\mathtt{f}\,dt=\infty$;
\item[(A.3)] $t^\nu\mathtt{g}$ is non-decreasing for some $\nu\in[0,1)$;
\item[(A.4)] $\psi$ diverges as $t\rightarrow\infty$.
\end{itemize}
Define continuous positive densities $f$ and $g$ on $\mathbb{R}^2_0$ via
\begin{align}
&f:\mathbb{R}^2_0\rightarrow(0,\infty);
x\mapsto\mathtt{f}(|x|)\label{definition_of_density}
\end{align}
and similarly for $g$. The weighted $f$-volume $V_f:=f\mathscr{L}^2$ is defined on the $\mathscr{L}^2$-measurable sets in $\mathbb{R}^2_0$. The weighted $g$-perimeter of a set $E$ of locally finite perimeter in $\mathbb{R}^2_0$ is defined by
\begin{equation}\label{weighted_perimeter}
P_g(E):=\int_{\mathbb{R}^2_0}g\,d|D\chi_E|\in[0,\infty].
\end{equation}
We refer to the isoperimetric problem (\ref{fg_isoperimetric_problem}) for $v>0$. We require a number of preparatory results before we reach the existence theorem. 

\smallskip

\begin{lemma}\label{lower_semicontinuity_result}
Let $\Omega$ be an open set in $\mathbb{R}^2$  and $g$ a positive continuous function on $\Omega$. Let $(u_h)_{h\in\mathbb{N}}$ be a sequence of functions in $\mathrm{BV}_{\mathrm{loc}}(\Omega)$ which converge to $u\in L^1_{\mathrm{loc}}(\Omega)$ in $L^1_{\mathrm{loc}}(\Omega)$. Assume that 
\[
\liminf_{h\rightarrow\infty}\int_\Omega g\,d|Du_h|<\infty.
\]
Then $u\in\mathrm{BV}_{\mathrm{loc}}(\Omega)$ and
\begin{align}
\int_\Omega g\,d|Du|&\leq\liminf_{h\rightarrow\infty}\int_\Omega g\,d|Du_h|.\label{lower_semicontinuity}
\end{align}
\end{lemma}

\smallskip

\noindent{\em Proof.}
Let $A$ be a relatively compact open set in $\Omega$. As $g$ is positive and continuous on $\Omega$ it is bounded away from zero on $\overline{A}$ by a positive constant $c$. By lower semi-continuity of the variation (see \cite{Ambrosio2000} Remark 3.5),
\begin{align}
V(u,A)&\leq\liminf_{h\rightarrow\infty}V(u_h,A)\leq\frac{1}{c}\liminf_{h\rightarrow\infty}\int_\Omega g\,d|Du_h|.\nonumber
\end{align}
By \cite{Ambrosio2000} Proposition 3.6 the limit function $u$ belongs to $\mathrm{BV}_{\mathrm{loc}}(\Omega)$. Then $gDu$ is an $\mathbb{R}^2$-valued Radon measure on $\Omega$. By \cite{Ambrosio2000} Remark 1.46 and Proposition 1.47,
\[
|gDu|(A)=\int_A g\,d|Du|.
\]
By \cite{Ambrosio2000} Proposition 3.13 the sequence $(u_h)_{h\in\mathbb{N}}$ weakly* converges to $u$ on $A$. In particular, the finite Radon measures $Du_h\rightarrow Du$ weakly* on $A$ and as a consequence $gDu_h\rightarrow gDu$ weakly* on $A$ as $h\rightarrow\infty$. By \cite{Ambrosio2000} Corollary 1.60,
\begin{align}
\int_A g\,d|Du|&=|gDu|(A)\leq\liminf_{h\rightarrow\infty}|gDu_h|(A)=\liminf_{h\rightarrow\infty}\int_A g\,d|Du_h|
\leq\liminf_{h\rightarrow\infty}\int_\Omega g\,d|Du_h|.\nonumber
\end{align}
The inequality in (\ref{lower_semicontinuity}) follows from the monotone convergence theorem.
\qed

\bigskip

 \noindent Let $E$ be a set of locally finite perimeter in $\mathbb{R}^2_0$. We write $E_t:=E\cap\mathbb{S}_t$ for the $t$-section of $E$ for each $t>0$. Define
\begin{align}
L(t)=L_E(t)&:=\mathscr{H}^1((E^1)_t)\label{definition_of_L}\\
p(t)=p_E(t)&:=\mathscr{H}^0((\mathscr{F}E)_t).\label{definition_of_p}
\end{align}
By the co-area formula (\cite{Ambrosio2000} Theorem 2.93) and the De Giorgi structure theorem (\cite{Ambrosio2000} Theorem 3.59) the function $p$ is $\mathscr{L}^1$-measurable on $(0,\infty)$. Note that $L$ does not depend on the $\mathscr{L}^2$-version of $E$. 

\smallskip
 
\noindent Let $B(t)$ stand for the open centred ball in $\mathbb{R}^2$ with radius $t\geq 0$; the closed centred ball is denoted $\overline{B}(t)$. Sometimes we simply write $\overline{B}$ for shortness. Let $E$ be a set of locally finite perimeter in $\mathbb{R}^2_0$. Let us define
\begin{align}
V_f(t)&:=V_f(E\cap\mathbb{R}^2_0\setminus\overline{B}(t));\nonumber\\
P_g(t)&:=\int_{\mathbb{R}^2_0\setminus\overline{B}(t)}g\,d|D\chi_E|;\nonumber
\end{align}
for $t\geq 0$. In this notation, $V_f(0)=V_f(E)$ and $P_g(0)=P_g(E)$.  If $f=g=1$ we add a subscript in the form $V_0$ and $P_0$.

\smallskip

\begin{lemma}\label{lower_bound_for_P_0}
Let $\mathtt{f}$ be a positive continuous function on $(0,\infty)$ which satisfies the conditions (A.1)-(A.2) and define $f$ as in (\ref{definition_of_density}). Let $E$ be a set of locally finite perimeter in $\mathbb{R}^2_0$ with $V_f(E)<\infty$. Then
\[
P_0(E,\mathbb{R}^2_0\setminus\overline{B}(t))\geq L(t)
\]
for each $t>0$.
\end{lemma}
 
 \smallskip
 
 \noindent{\em Proof.} Let us fix $t>0$ and write $\overline{B}$ instead of $\overline{B}(t)$ for shortness. We may assume that $E$ has finite perimeter relative to $\mathbb{R}^2_0\setminus\overline{B}$. Let us first note that $\mathscr{F}E$ is countably $1$-rectifiable and 
 \[
 P_0(E,\mathbb{R}^2_0\setminus\overline{B})=\mathscr{H}^1(\mathscr{F}E\setminus\overline{B})
 \]
 as in \cite{Ambrosio2000} Theorem 3.59. The map $\alpha:\mathbb{R}^2\setminus\overline{B}\rightarrow\mathbb{S}_t$ defined by $x\mapsto(t/|x|)x$ is $1$-Lipschitz. We remark in passing that $\alpha$ admits an extension $\widetilde{\alpha}:\mathbb{R}^2\rightarrow\mathbb{S}_t$ that is $\sqrt{2}$-Lipschitz by \cite{Ambrosio2000} Proposition 2.12. By the generalised area formula (cf. \cite{Ambrosio2000} Theorem 2.91),
\[
\int_{\mathbb{R}^2_0}\mathscr{H}^0(\mathscr{F}E\setminus\overline{B}\cap\alpha^{-1}(\omega))\,\mathscr{H}^1(d\omega)
=
\int_{\mathscr{F}E\setminus\overline{B}}J_1d^{\mathscr{F}E}\alpha_x\,d\mathscr{H}^1
\leq
\mathscr{H}^1(\mathscr{F}E\setminus\overline{B}).
\]

\smallskip

\noindent Put $F:=E^1\setminus\overline{B}$. Consider the set $\Lambda$ of all points $\omega\in\mathbb{S}_t$ which satisfy the properties
\begin{itemize}
\item[(a)] $\omega\in(E^1)_t$;
\item[(b)] $F_\omega$ is a set of locally finite perimeter in $(\mathbb{R}^2)_\omega$;
\item[(c)] $(\mathscr{F}F)_\omega=\mathscr{F}F_\omega$;
\item[(d)] $\omega\in\mathscr{F}F_\omega$.
\end{itemize}
We argue that $\Lambda$ has full $\mathscr{H}^1$-measure in $(E^1)_t$. First note that the set $F$ is a set of locally finite perimeter in $\mathbb{R}^2_0$ by \cite{Ambrosio2000} Corollary 3.89. By Theorem \ref{radialVolperttheorem} the set $F_\omega$ is a set of locally finite perimeter in $(\mathbb{R}^2)_\omega$ and the set $\{\omega:\mathscr{F}F_\omega=(\mathscr{F}F)_\omega\}$ has full measure in $\mathbb{S}_t$ for $\mathscr{H}^1$-a.e. $\omega\in\mathbb{S}_t$. Note that if $\omega\in(E^1)_t$ then $\omega\in F^{1/2}$. By Federer's Theorem (cf. \cite{Ambrosio2000} Theorem 3.61) it follows that $(E^1)_t$ is contained in the reduced boundary of $F$ apart from at most an $\mathscr{H}^1$-null set. These considerations entail that the conditions (a)-(d) hold a.e. in $(E^1)_t$.

\smallskip

\noindent Assume for a moment that the set $H:=\Lambda\cap\{\omega:\mathscr{F}F_\omega\setminus\{\omega\}=\emptyset\}$ has positive $\mathscr{H}^1$-measure in $\mathbb{S}_t$.  For each $\omega$ in this latter set $F_\omega$ is equivalent to the ray $(1,\infty)\omega$ by \cite{Ambrosio2000} Proposition 3.52. We infer that the truncated cone $(1,\infty)H$ is contained within $F$ up to a Lebesgue null set. This truncated cone  and hence $E$ has infinite $V_f$-volume on appealing to property (A.2) which contradicts our hypotheses that $V_f(E)<\infty$. We conclude from this that $\Lambda\cap\{\omega:\mathscr{F}F_\omega\setminus\{\omega\}\neq\emptyset\}$ has full measure in $(E^1)_t$.

\smallskip

\noindent Suppose that $\omega$ belongs to the set $\Lambda$ and $\mathscr{F}F_\omega\setminus\{\omega\}\neq\emptyset$. By (c), $(\mathscr{F}F)_\omega\setminus\{\omega\}\neq\emptyset$ and this implies in turn that $\omega$ belongs to the range of $\alpha$ restricted to $\mathscr{F}E\setminus\overline{B}$. So
\[
\int_{\mathbb{R}^2_0}\mathscr{H}^0(\mathscr{F}E\setminus\overline{B}\cap\alpha^{-1}(\omega))\,\mathscr{H}^1(d\omega)
\geq
\int_{\mathbb{R}^2_0}\chi_{(E^1)_t}\,d\mathscr{H}^1=L(t).
\]
This leads to the result when combined with the above inequality.  
 \qed
   
\smallskip

\begin{lemma}\label{lower_bound_on_P_nu}
Let $\nu\in[0,1)$ and $\mathtt{f}$ be a positive continuous function on $(0,\infty)$ which satisfies the conditions (A.1)-(A.2) and define $f$ as in (\ref{definition_of_density}). Let $E$ be a set of locally finite perimeter in $\mathbb{R}^2_0$ with $V_f(E)<\infty$. Then
\[
 P_{-\nu}(t)\geq(|\cdot|^{-\nu} L)(t)
 \]
for each $t>0$.
\end{lemma}

\smallskip

\noindent{\em Proof.}
Let $t>0$.  By (\ref{J_1_of_Phi}) and the generalised area formula (cf. \cite{Ambrosio2000} Theorem 2.91),
\begin{align}
P_{-\nu}(t)
&\geq(-\nu+1)\mathscr{H}^1(\Phi(\mathscr{F}E\setminus\overline{B}))\nonumber\\
&=(-\nu+1)\mathscr{H}^1(\mathscr{F}\Phi(E)\setminus\Phi(\overline{B}))
=(-\nu+1)P_0(\Phi(E),\mathbb{R}^2_0\setminus\Phi(\overline{B}))\label{equality_for_P_nu_perimeter}
\end{align}
where in the first equality we use the fact that $\Phi(\mathscr{F}E)$ is $\mathscr{H}^1$-a.e. equivalent to $\mathscr{F}\Phi(E)$ according to Lemma \ref{diffeomorhism_invariance_of_boundary}. The push-forward $\Phi_\sharp V_\alpha$ of $V_\alpha$ under $\Phi$ is a non-finite Radon measure on $\mathbb{R}^2_0$ with radial density satisfying (A.1)-(A.2) and $\Phi(E)$ has finite $\Phi_\sharp V_\alpha$-measure. We make use of  Lemma \ref{lower_bound_for_P_0} to continue,
\begin{align}
P_{-\nu}(t)
\geq(-\nu+1)\mathscr{H}^1((\Phi(E)^1)_{\frac{t^{-\nu+1}}{-\nu+1}})
=(|\cdot|^{-\nu} L)(t)
\nonumber
\end{align}
where once more we appeal to Lemma \ref{diffeomorhism_invariance_of_boundary}. This establishes the result.
\qed
  
\smallskip
 
\begin{proposition}\label{lower_bounds_for_P_t}
Suppose that the positive continuous function $\mathtt{g}$ satisfies condition (A.3) and the positive continuous function $\mathtt{f}$ satisfies conditions (A.1)-(A.2). Let $E$ be a set of locally finite perimeter in $\mathbb{R}^2_0$ with $V_f(E)<\infty$. Assume that $P_g(t)<\infty$ for each $t>0$. Then for each $t>0$,
\begin{itemize}
\item[(i)] $P_g(t)\geq\int_t^\infty\mathtt{g}p\,d\tau$;
\item[(ii)] $P_g(t)\geq(\mathtt{g}L)(t)$;
\end{itemize}
and
\begin{itemize}
\item[(iii)] the set $\{L=2\pi\tau\}$ is bounded in $(0,\infty)$.
\end{itemize}
Moreover if $\mathtt{f}$ satisfies condition (A.4) then
 \begin{itemize}
\item[(iv)] there exists $T>0$ such that $P_g(t)^2\geq\psi_-(t)V_f(t)$ for each $t>T$;
\end{itemize}
where $\psi_-(t):=\min_{[t,\infty)}\psi$.
\end{proposition}

\smallskip

\noindent{\em Proof.}
{\em (i)} By Proposition \ref{tangential_integral_equality},
\[
\int_{\mathscr{F}E\setminus\overline{B}}g|\nu_E^\omega|\,d\mathscr{H}^1
=\int_t^\infty\mathtt{g}\mathscr{H}^0((\mathscr{F}E)_\tau)\,d\tau
\]
for each $t>0$ and the statement follows by the definition in (\ref{definition_of_p}).

\smallskip

\noindent{\em (ii)} Choose $\nu\in[0,1)$ as in (A.3). By Lemma \ref{lower_bound_on_P_nu},
\begin{align}
P_g(t)&=\int_{\mathscr{F}E\setminus\overline{B}}(|x|^\nu g)|x|^{-\nu}\,d\mathscr{H}^1
\geq t^\nu\mathtt{g}(t)\int_{\mathscr{F}E\setminus\overline{B}}|x|^{-\nu}\,d\mathscr{H}^1
=t^\nu\mathtt{g}(t)P_{-\nu}(t)
\geq (gL)(t).\nonumber
\end{align}

\smallskip

\noindent{\em (iii)} For $t>s>0$ we may write $t\mathtt{g}(t)=(t^\nu\mathtt{g}(t))t^{1-\nu}$ with $\nu\in[0,1)$ as in (A.3). In light of this condition the function $t\mathtt{g}$ diverges as $t\uparrow\infty$. Suppose that the set $A:=\{L=2\pi t\}$ is unbounded. By {\em (ii)}, $P_g(s)\geq P_g(t)\geq(\mathtt{g}L)(t)=2\pi t\mathtt{g}(t)$ for each $t\in A$. This contradicts the fact that $P_g(s)<\infty$.

\smallskip

\noindent{\em (iv)} For each $s>t>0$,
\[
\infty>P_g(s)\geq P_g(t)\geq(\mathtt{g}L)(t)=\psi(t)(\zeta L)(t)
\]
by {\em (ii)}. It follows that $(\zeta L)(t)\rightarrow 0$ as $t\rightarrow\infty$ because $\psi$ diverges according to (A.4). This means in particular that $\zeta L$ is bounded on $[t,\infty)$ for each $t>0$. Put
\[
M(t):=\sup\Big\{(\zeta L)(\tau):\tau\geq t\Big\}<\infty.
\]
From the above estimate $P_g(t)\geq\psi_-(t)M(t)$. 

\smallskip

\noindent Let $T$ be an upper bound for the set $A$ according to {\em (iii)}. Choose $t>T$. We may assume that $V_f(t)=\int_t^\infty\mathtt{f}L\,d\tau>0$ otherwise the inequality follows immediately. This latter condition entails that $M(t)>0$. This in turn implies that $\zeta L/M(t)\leq 1$ on $[t,\infty)$. Note that $p\geq 1$ on $\{p>0\}$. With the justification supplied by these comments and {\em (i)} we may write
\begin{align}
P_g(t)&\geq\int_t^{\infty}\mathtt{g}p\,d\tau
=\int_t^{\infty}\chi_{\{p>0\}}\mathtt{g}p\,d\tau
\geq\int_t^{\infty}\chi_{\{p>0\}}\mathtt{g}\,d\tau
\geq\frac{1}{M(t)}\int_t^{\infty}\chi_{\{p>0\}}\psi(\zeta^2 L)\,d\tau\nonumber\\
&=\frac{1}{M(t)}\int_t^{\infty}\chi_{\{L>0\}}\mathtt{f}L\,d\tau
=\frac{1}{M(t)}\int_t^{\infty}\mathtt{f}L\,d\tau
=\frac{1}{M(t)}V_f(t).\nonumber
\end{align}
In the first equality in the second line we make use of the the fact that the sets $\{p>0\}$ and $\{L>0\}$ coincide up to a set of measure zero on $(T,\infty)$. Let us explain this. Let $\Lambda$ signify the set of $t>0$ with the property that both $(E^1)_t$ is a set of finite perimeter in $\mathbb{S}_t$ and $(\mathscr{F}E)_t=\mathscr{F}(E^1)_t$. By Theorem \ref{tangentialVolperttheorem} this is a set of full measure in $(0,\infty)$. For $t\in\Lambda$,
\begin{align}
L(t)&=\mathscr{H}^1((E^1)_t)\text{ and }p(t)=\mathscr{H}^0((\mathscr{F}E)_t)=\mathscr{H}^0(\mathscr{F}(E^1)_t).\label{formulae_for_L_and_p}
\end{align}
Moreover,
\[
\Lambda\cap\{0<L\}\cap(T,\infty)=\Lambda\cap\{0<L<2\pi t\}\cap(T,\infty)
\]
and this last coincides with the set $\Lambda\cap\{0<p\}\cap(T,\infty)$.

\smallskip

\noindent As a last step the product of the estimates in the first and second paragraphs of this section leads to the result in {\em (iv)}. 
\qed

\smallskip

\begin{theorem}\label{existence_theorem}
Suppose that $\mathtt{f}$ and $\mathtt{g}$ satisfy (A.1)-(A.4). Then (\ref{fg_isoperimetric_problem}) has a minimiser for each $v>0$.
\end{theorem}
 
 \smallskip
 
\noindent{\em Proof.}
We adapt the proof of \cite{MorganPratelli2011} Theorem 3.3 (or \cite{DiGiosiaetal2016} Theorem 2.3). Let $(E_h)_{h\in\mathbb{N}}$ be a minimising sequence for (\ref{fg_isoperimetric_problem}).  By \cite{Ambrosio2000} Theorem 3.23 we may assume that the sequence $(u_h)_{h\in\mathbb{N}}$ converges to a function $u\in\mathrm{BV}_{\mathrm{loc}}(\Omega)$ in $L^1_{\mathrm{loc}}(\Omega)$ where $\Omega:=\mathbb{R}^2_0$ and $u_h:=\chi_{E_h}$. By choosing an a.e. convergent subsequence we may assume that the limit function $u$ takes the form $u=\chi_E$ for some Borel set $E$ in $\Omega$. Assume for the moment that
\begin{align}
\lim_{R\rightarrow\infty}\sup_{h\in\mathbb{N}}V_f(E_h\setminus B(R))&=0.\label{condition_for_existence}
\end{align}
By \cite{Ambrosio2000} Proposition 1.27 the collection $\mathscr{F}:=\{u_h:h\in\mathbb{N}\}$ is an equi-integrable subset of $L^1(\Omega,V_f)$. By the Dunford-Pettis Theorem (\cite{Ambrosio2000} Theorem 1.38 for example) and uniqueness of the weak limit we may assume that $(u_h)_{h\in\mathbb{N}}$ converges to $u$ weakly in $L^1_{\mathrm{loc}}(\Omega)$. In particular, $V_f(E)=v$. So $E$ is a minimiser for the variational problem (\ref{fg_isoperimetric_problem}) by Proposition \ref{lower_semicontinuity_result}.

\smallskip

\noindent Let us suppose for a contradiction that the condition (\ref{condition_for_existence}) does not hold; in particular, there exists $\varepsilon>0$ such that for each $t>0$ there exists $h\in\mathbb{N}$ with the property that $V_f(E_h\setminus B(t))\geq\varepsilon$. We first remark that there exists $T>0$ with the property that $\{L_h=2\pi t\}\cap(T,\infty)=\emptyset$ for each $h\in\mathbb{N}$ (in an obvious notation). This is a consequence of the fact that the set $\{P_g(E_h):h\in\mathbb{N}\}$ is bounded in $\mathbb{R}$ and the estimate in Proposition \ref{lower_bounds_for_P_t} {\em (ii)}. Let $t>T$ and choose $h$ as above. 
By Proposition \ref{lower_bounds_for_P_t},
\[
P_g(E_h)^2\geq\psi_-(t)\varepsilon.
\]
The right-hand side in this inequality diverges as $t\uparrow\infty$ while the left-hand side remains bounded. 
This gives a contradiction. We conclude that the condition (\ref{condition_for_existence}) holds true.
\qed

\smallskip

\begin{lemma}\label{BV_functions_and_positivity_set}
Let $u$ be a non-negative function in $\mathrm{BV}_{\mathrm{loc}}((0,\infty))$ with $u\not\equiv 0$. Then there exists a good representative $\overline{u}$ of $u$ with the property that the set $\{\overline{u}>0\}$ is open in $(0,\infty)$.
\end{lemma}

\smallskip

\noindent{\em Proof.}
Let $\overline{u}$ be a good representative of $u$ (cf. \cite{Ambrosio2000} Theorem 3.28). Denote by $A$ the set of atoms of $Du$; a countable set as $|Du|$ is a Radon measure. Put
\[
A_1:=\{t\in A:\overline{u}(t-)=0\text{ or }\overline{u}(t+)=0\}
\]
and define $\overline{u}_1:=\chi_{(0,\infty)\setminus A_1}\overline{u}$. Then $\overline{u}_1$ is a good representative of $\overline{u}$ by \cite{Ambrosio2000} Theorem 3.28. Moreover, the set $\{\overline{u}_1>0\}$ is open in $(0,\infty)$.
\qed

\smallskip

\noindent The next is contained in \cite{CagnettiPeruginiStoger2020} Lemma 4.1. We include its short proof.

\smallskip

\begin{lemma}\label{L_is_BV_loc}
Let $E$ be a set of locally finite perimeter in $\mathbb{R}^2_0$. Then $L\in\mathrm{BV}_{\mathrm{loc}}((0,\infty))$.
\end{lemma}

\smallskip

\noindent{\em Proof.}
Let $\phi\in C^1_c((0,\infty))$. Then
\[
\int_0^\infty\phi^\prime(L/t)\,dt=-\int_{\mathscr{F}E}\frac{1}{t}\phi\nu_E^\omega\,d\mathscr{H}^1.
\]
So $L/t$ belongs to $\mathrm{BV}_{\mathrm{loc}}((0,\infty))$ and likewise $L$ by \cite{Ambrosio2000} Proposition 3.2.
\qed

\smallskip

\begin{lemma}\label{on_existence_of_complementary_annulus}
Suppose that the positive continuous functions $\mathtt{f}$ and $\mathtt{g}$ on $(0,\infty)$ satisfy the conditions (A.1)-(A.5). Let $E$ be a set of locally finite perimeter in $\mathbb{R}^2_0$ with $V_f(E)<\infty$. Assume $P_g(t)<\infty$ for each $t>0$. Then for each $t>0$ there exists $a$ and $b$ with $t<a<b<\infty$ such that $|E\cap A(a,b)|=0$.
\end{lemma}

\smallskip

\noindent The notation $A(a,b)$ refers to the open annulus with inner radius $a$ and outer radius $b$.

\smallskip

\noindent{\em Proof.} 
First suppose that $\mathtt{g}\equiv 1$. In other words, $E$ is a set of locally finite perimeter in $\mathbb{R}^2_0$ with finite $V_f$-measure and $P_0(t)<\infty$ for each $t>0$. Let us decompose $(0,\infty)$ into a disjoint union
\[
(0,\infty)=\{L=0\}\cup\{0<L<2\pi\tau\}\cup\{L=2\pi\tau\}.
\]
By Proposition \ref{lower_bounds_for_P_t} {\em (iii)} the set $\{L=2\pi\tau\}$ is bounded and has finite measure. The set $\{0<L<2\pi\tau\}$ is contained in $\{p>0\}$ up to a Lebesgue null set bearing in mind ({\ref{formulae_for_L_and_p}). Moreover the set $\{p>0\}$ has finite measure in $(0,\infty)$ thanks to Proposition \ref{lower_bounds_for_P_t} {\em (i)}. The upshot of this is that the set $\{0<L<2\pi\tau\}$ has finite measure. Accordingly the set $\{L=0\}$ has infinite measure and so does $\{L=0\}\cap(t,\infty)$ for any $t>0$. 

\smallskip

\noindent The function $L$ is a nonnegative function in $\mathrm{BV}_{\mathrm{loc}}((0,\infty))$ by the previous Lemma \ref{L_is_BV_loc}. Choose a good representative $\overline{L}$ of $L$ as in Lemma \ref{BV_functions_and_positivity_set}. Then the set $\{\overline{L}>0\}$ is open. This latter set is a countable disjoint union of open intervals. Its complement in $(0,\infty)$ is a countable union of disjoint closed intervals relative to $(0,\infty)$. The intersection of this complement with $(t,\infty)$ has infinite measure for each $t>0$. So $\{\overline{L}=0\}$ contains a closed interval $[a,b]$ in $(t,\infty)$ with positive length for any choice of $t>0$. It follows that the intersection of the set $E$ with the annulus $A(a,b)$ has measure zero. 

\smallskip

\noindent Now suppose that the functions $\mathtt{f}$ and $\mathtt{g}$ satisfy the conditions (A.1)-(A.4). Let $\nu\in[0,1)$ and $\Phi$ as in (\ref{Morgan_map}) with $\beta$ replaced by $-\nu$. By (A.3) and (\ref{equality_for_P_nu_perimeter}),
\begin{align}
P_g(t)&=\int_{\mathscr{F}E\setminus\overline{B}}g\,d\mathscr{H}^1
\geq t^\nu\mathtt{g}(t)P_{-\nu}(t)
=(-\nu+1) t^\nu\mathtt{g}(t)P_0(\Phi(E),\mathbb{R}^2_0\setminus\Phi(\overline{B}));\nonumber
\end{align}
that is, the set $\Phi(E)$ has the properties mentioned in the first paragraph of this proof. The set $\Phi(E)$ has finite $\Phi_\sharp V_\alpha$-measure and this measure has radial density satisfying (A.1)-(A.2). By the result for the case $\mathtt{g}\equiv 1$ for each $t>0$ there exists $a$ and $b$ with  $\frac{t^{-\nu+1}}{-\nu+1}<a<b<\infty$ such that $|\Phi(E)\cap A(a,b))|=0$. As $\Phi$ preserves centred circles a similar result holds for $E$.
\qed

\smallskip

\begin{theorem}\label{boundedness_theorem}
Suppose that $\mathtt{f}$ and $\mathtt{g}$ satisfy (A.1)-(A.4). Then for each $v>0$ any minimiser of (\ref{fg_isoperimetric_problem}) is bounded.
\end{theorem}

\smallskip

\noindent{\em Proof.}
Let $v>0$ and $E$ be a set of locally finite perimeter in $\mathbb{R}^2_0$ which is a minimiser for (\ref{fg_isoperimetric_problem}) according to Theorem \ref{existence_theorem}. Let us assume that the topological boundary $\partial E\cap\mathbb{R}^2_0$ of $E$ in $\mathbb{R}^2_0$ enjoys the property in (\ref{topological_boundary}). Assume for a contradiction that $E$ is not bounded. This means that $V_f(t)>0$ for each $t>0$. Choose a centred open ball $B_0$ with weighted volume $v$. So $\delta:=V_f(B_0\setminus E)>0$. In fact $\delta=V_f(B_0\setminus E^1)=V_f(B_0\cap E^0)+V_f(B_0\cap\partial^*E)=V_f(B_0\cap E^0)>0$ by Federer's Theorem (\cite{Ambrosio2000} Theorem 3.61). Assume for a moment that the interior of the set $B_0\cap E^{0}$ is empty. Then $B_0\cap E^{0}\subset\partial E\cap\mathbb{R}^2_0$. To see this consider a point $x$ in $B_0\cap E^0$. The complement of $E$ has density $1$ at $x$ as this point belongs to $E^0$. On the other hand if $E\cap B(x,\rho)$ is empty for some $\rho>0$ then $B_0\cap E^0$ has non-empty interior. This inclusion contradicts Lemma \ref{on_toplogical_boundary}. So the set $B_0\cap E^{0}$ has non-empty interior.

\smallskip

\noindent The upshot of the argument above is that there exists an open ball $B_1$ contained in $B_0$ with the properties $0<V_f(B_1)<\delta$, $|B_1\cap E|=0$ and $0\not\in\overline{B}_1$. Put
\[
\eta:=\sup\Big\{\frac{P_g(B)^2}{V_f(B)}:B\text{ is an open ball contained in }B_1\Big\}<\infty.
\]
Bearing in mind property (A.4) as well as Lemma \ref{on_existence_of_complementary_annulus} and our hypothesis on the behaviour of $V_f(\cdot)$ we may choose $t>0$ such that
\begin{itemize}
\item[(a)] $\psi_-(t)>\eta$;
\item[(b)] $V_f(t)<V_f(B_1)$;
\item[(c)] there exists $\delta_1>0$ such that $|A(t-\delta_1,t+\delta_1)\cap E|=0$.
\end{itemize}
Choose an open ball $B_2$ in $B_1$ with $V_f(B_2)=V_f(t)$. Define $E_1:=(E\cap B(t))\cup B_2$. We claim that $P_g(E_1)< P_g(E)$. We first note that
\[
\frac{P_g(B_2)^2}{V_f(t)}=\frac{P_g(B_2)^2}{V_f(B_2)}\leq\eta<\psi_-(t)\leq\frac{P_g(E,\mathbb{R}^2\setminus\overline{B}(t))^2}{V_f(t)}
\]
by Proposition \ref{lower_bounds_for_P_t}. By \cite{Ambrosio2000} Corollary 3.89,
\[
P_g(E_1)= P_g(E\cap B(t))+P_g(B_2)<P_g(E\cap B(t))+P_g(E,\mathbb{R}^2\setminus\overline{B}(t))=P_g(E)
\]
This gives rise to the contradiction that the set $E$ is in fact not a minimiser for (\ref{fg_isoperimetric_problem}). This proves the theorem.
\qed

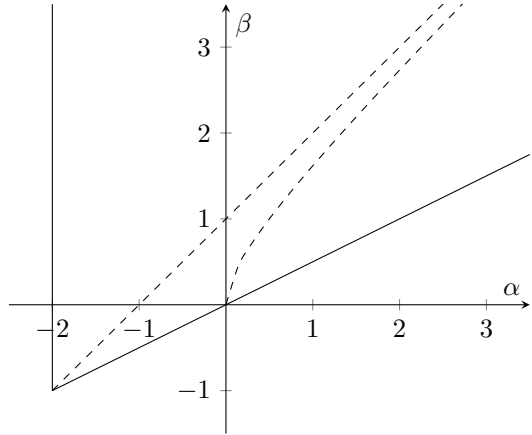
\begin{figure}[h]
\centering
\begin{tikzpicture}[>=stealth]
    \begin{axis}[
        xmin=-2.5,xmax=3.5,
        ymin=-1.5,ymax=3.5,
        axis x line=middle,
        axis y line=middle,
        axis line style=->,
        xlabel={$\alpha$},
        ylabel={$\beta$},
        ]
        \addplot[black, domain=-2:3.5, smooth]{x/2};
        \draw[black] (-2,-1) -- (-2,3.5);

        \addplot[black, dashed, domain=0:3.5, smooth]{(x+sqrt(x^2+4*x))/2};
        \addplot[black, dashed, domain=-2:3.5, smooth]{x+1};
    \end{axis}
\end{tikzpicture}
\caption{The region $\mathcal{Q}$ is enclosed between the solid lines $\alpha+2=0$ and $\alpha-2\beta=0$. The region $\mathcal{P}$ is enclosed between the dashed lines and is contained inside $\mathcal{Q}$.}
\label{fig:Region_Q}
\end{figure}

\smallskip

\noindent The next result is closely related to \cite{MorganRitore2002} Theorem 1.1 (and \cite{BrayMorgan2001} Corollary 2.4) upon parametrising the cone $C$ in (\ref{cone}) using the map $\varphi$.

\smallskip

\begin{theorem}\label{isoperimetry_on_P_-}
Let $(\alpha,\beta)\in\mathcal{P}^-\setminus\{0\}$. Then centred balls are unique minimisers for the problem (\ref{isoperimetric_problem}) up to equivalence.
\end{theorem}

\smallskip

\noindent{\em Proof.}
We prove uniqueness. Let $B$ be the centred open ball with volume $V_f(E)=V_f(B)=v$. Assume that $P_\beta(E)=P_\beta(B)$. In this case equality holds in the inequality (\ref{J_1_of_Phi}) and $\nu_E^\omega=0$ $\mathscr{H}^1$-a.e. on $\mathscr{F}E$. From Proposition \ref{tangential_integral_equality}, $(\mathscr{F}E)_t$ is empty for a.e. $t>0$. By Theorem \ref{tangentialVolperttheorem}, $(E^1)_t$ is equivalent to either the empty set $\emptyset$ or $\mathbb{S}_t$ for a.e. $t>0$. Let $R$ stand for the radius of $B$. Suppose $(E^1)_t$ is equivalent to $\mathbb{S}_t$ for some $t>R$. By Proposition \ref{lower_bounds_for_P_t}, $P_\beta(E)\geq P_\beta(t)\geq(|\cdot|^\beta L)(t)=2\pi t^{\beta+1}>2\pi R^{\beta+1}=P_\beta(B)$ as $\beta+1>0$. This is a contradiction. So in fact $(E^1)_t$ is equivalent to the empty set for a.e. $t>0$. This implies that $E$ is equivalent to $B$.
\qed

\smallskip

\noindent Let us now apply the above existence and boundedness theorems to the case of radial power weights which are the focus of this paper. We introduce the parameter set 
\begin{align}
\mathcal{Q}&:=\Big\{(\alpha,\beta)\in\mathbb{R}^2:-2<\alpha\text{ and }\alpha\leq 2\beta\text{ if }\beta\leq 0
\text{ or }\alpha<2\beta\text{ if }\beta>0\Big\}
\end{align}
illustrated in Figure \ref{fig:Region_Q}. 

\smallskip

\begin{corollary}\label{existence_and_boundedness_for_indices_in_curly_Q}
For each $(\alpha,\beta)\in\mathcal{Q}$ the isoperimetric problem (\ref{isoperimetric_problem}) has a minimiser for each $v>0$ and any such minimiser is bounded.
\end{corollary}

\smallskip

\noindent{\em Proof.}
For $(\alpha,\beta)\in\mathcal{Q}$ with $\alpha<2\beta$ this follows from Theorem \ref{existence_theorem} and Theorem \ref{boundedness_theorem}. Now consider a pair of indices $(\alpha,\beta)$ satisfying $\alpha=2\beta$ and $-1<\beta<0$. Then existence and boundedness is a special case of Theorem \ref{isoperimetry_on_P_-}. The case $(\alpha,\beta)=(0,0)$ corresponds to the classical isoperimetric inequality.
\qed

\smallskip

\noindent The parameter set $\mathcal{P}$ is included in the set $\mathcal{Q}$. This is because
\[
\alpha=\frac{\beta^2}{\beta+1}<\beta<2\beta
\]
whenever $\beta>0$.

\section{Regularity, spherical cap symmetry and curvature}\label{Section_on_regularity_etc}

\smallskip

\noindent The isoperimetric problem (\ref{isoperimetric_problem}) has a bounded minimiser after Corollary \ref{existence_and_boundedness_for_indices_in_curly_Q} so long as the indices $(\alpha,\beta)$ belong to the parameter set $\mathcal{Q}$. In this Section we collate some symmetry, regularity and curvature properties of the minimiser. Let us start with spherical cap symmetry and analyticity. The notion of spherical cap symmetry is discussed in for example \cite{McGillivray2018_1} Section 4. 

\smallskip

\begin{theorem}\label{C1_property_of_reduced_boundary}
Let $(\alpha,\beta)\in\mathcal{Q}$. Let $v>0$ and suppose that $E$ is a bounded minimiser of (\ref{isoperimetric_problem}). Then there exists an $\mathscr{L}^2$-measurable set $\widetilde{E}$ in $\mathbb{R}^2_0$ with the properties
\begin{itemize}
\item[(i)] $\widetilde{E}$ is a minimiser of (\ref{isoperimetric_problem});
\item[(ii)] $L_{\widetilde{E}}=L_E$ a.e. on $(0,\infty)$;
\item[(iii)] $\widetilde{E}$ is open with analytic boundary in $\mathbb{R}^2_0$;
\item[(iv)] $\widetilde{E}=\widetilde{E}^{\mathrm{sc}}$.
\end{itemize}
\end{theorem}

\smallskip

\noindent{\em Proof.} See \cite{McGillivray2018_1} Theorem 3.1 and Theorem 4.8. A regularity result for manifolds with density is contained in \cite{Morgan2003} 3.10. In our case the metric may be degenerate or singular at the origin while the density relative to the metric may vanish at the origin (admittedly a minor change). Let us nevertheless demonstrate the analyticity property in {\em (iii)}. Let $E$ be an isoperimetric set. By \cite{McGillivray2018_1} Theorem 3.7 the boundary of $E$ is almost minimal and hence $M=\partial E$ is of class $C^1$ in $\mathbb{R}^2_0$ by \cite{Tamanini1984} Theorem 1.9. Consider a point $p\in M$ and assume that the normal to $M$ at $p$ is not parallel to the $x_1$-axis. There exists an open interval $G$ not containing the origin and $u\in C^1(G)$ such that $M$ is locally the graph of $u$,
\[
M\cap Q=\{(x_1,u(x_1)):x_1\in G\}
\]
where $Q$ is the vertical strip $G\times\mathbb{R}$ (making this latter assumption for the sake of simplicity).  Let $\zeta\in C^1_c(G)$ and define $X\in C^1_c(\mathbb{R}^2_0,\mathbb{R}^2)$ by $X(x)=\zeta(x_1)\phi(x_2)e_2$ for a suitable cut-off function $\phi\in C^\infty_c(\mathbb{R})$. We may assume that the inner unit normal to $\partial E$ is given locally by
\[
\nu_E=\frac{(-u^\prime,1)}{\sqrt{1+(u^\prime)^2}}
\]
and 
\[
\mathrm{div}^M X
=\mathrm{div}X-\langle\nu_E,D_{\nu_E}X\rangle
=-\langle\nu_E,D_{\nu_E}X\rangle
=\frac{u^\prime\zeta^\prime}{1+(u^\prime)^2}
\]
in a neighbourhood of $M$ in $Q$. The Jacobean determinant of the parametrisation $\varphi:G\rightarrow M;x_1\mapsto(x_1,u(x_1))$ is $J_1d\varphi=\sqrt{1+(u^\prime)^2}$. As in \cite{McGillivray2018_1} Proposition 6.4 (for example) and making use of the area formula \cite{Ambrosio2000} Theorem 2.71,
\begin{align}
0&=\int_{M}\langle\nabla g, X\rangle + g\,\mathrm{div}^{M}X
-\lambda f\langle\nu_E, X\rangle\,d\mathscr{H}^{1}\nonumber\\
&=\int_G\big[\partial_{x_2} g(\cdot,u)\sqrt{1+(u^\prime)^2}-\lambda\,f(\cdot,u) \big]\zeta
+g(\cdot,u)\frac{u^\prime\zeta^\prime}{\sqrt{1+(u^\prime)^2}}
\,dx_1.\label{Euler_Lagrange_identity}
\end{align}
Define
\[
I(u,G):=\int_G j(\cdot,u,u^\prime)\,dx_1
\text{ with }
j(x_1,z,p):=g(x_1,z)\sqrt{1+p^2}-\lambda\,F(x_1,z)
\]
where the function $F$ is chosen such that $\partial_{x_2} F=f$ on $Q$. By (\ref{Euler_Lagrange_identity}), $u$ is a (Lipschitz) extremal in the sense of \cite{Morrey1966} 1.10 (see also (1.4.7)). Then $j$ is regular and analytic on $G\times\mathbb{R}\times\mathbb{R}$. By \cite{Morrey1966} Theorem 1.10.4, $u$ is analytic on $G$.
\qed

\medskip

\noindent Suppose the open set $E$ in $\mathbb{R}^2_0$ has $C^1$ boundary $M$ in $\mathbb{R}^2_0$. Denote by $n:M\rightarrow\mathbb{S}$ the inner unit normal vector field. Given $p\in M$ we choose a tangent vector $t(p)\in\mathbb{S}$ in such a way that the pair $\{t(p),n(p)\}$ forms a positively oriented basis for $\mathbb{R}^2$. For $p\in M$ let $\sigma(p)$ stand for the angle measured anti-clockwise from the position vector $p$ to the tangent vector $t(p)$; $\sigma(p)$ is uniquely determined up to integer multiples of $2\pi$. This angle will feature frequently in our considerations.
  
\smallskip
 
\noindent Let $E$ be an open set in $\mathbb{R}^2_0$ with $C^2$ boundary $M$ in $\mathbb{R}^2_0$. Let $p\in M$. There exists a local parametrisation $\gamma_1:I\rightarrow M$ where $I=(-\delta,\delta)$ for some $\delta>0$ of class $C^2$ with $\gamma_1(0)=p$. We always assume that $\gamma_1$ is parametrised by arc-length and that $\dot{\gamma_1}(0)=t(p)$ where the dot signifies differentiation with respect to arc-length.  The (geodesic) curvature $k_1$ is defined on $I$ via the relation
 \begin{equation}
\ddot{\gamma}_1=k_1n_1.
\end{equation}
The curvature $k$ of $M$ is defined on $M$ by
 \begin{equation}
 k(x):=k_1(s)
 \end{equation}
 whenever $x=\gamma_1(s)$ for some $s\in I$. 

\smallskip 
 
\noindent With $E$ as above define $\Lambda:=M\cap\{\cos\sigma=0\}$ and 
\begin{align}
\Lambda_1&:=\{x\in M:\mathscr{H}^1(\Lambda\cap B(x,\rho))=\mathscr{H}^1(M\cap B(x,\rho))\text{ for some }\rho>0\}.\label{definition_of_Lambda_1}
\end{align}
This set corresponds to the union of open centred circular arcs contained in the boundary $M$. Put $\Lambda_1^{\pm}:=\Lambda_1\cap\{x\in M:\pm\langle x,n\rangle>0\}$. The next important Theorem characterises the generalised (mean) curvature
\[
\tau^{\beta-\alpha}\Big\{k+\frac{\beta}{\tau}\sin\sigma\Big\}
\]
of the boundary of an isoperimetric minimiser in the spherical cap symmetric case. 

\smallskip

\begin{theorem}\label{constant_weighted_mean_curvature}
Let $(\alpha,\beta)\in\mathcal{Q}$. Given $v>0$ let $E$ be a minimiser of (\ref{isoperimetric_problem}). Assume that $E$ is a bounded open set with analytic boundary $M$ in $\mathbb{R}^2_0$ and suppose that $E=E^{sc}$. Then there exists $\lambda\in\mathbb{R}$ such that 
\[
k+\frac{\beta}{\tau}\sin\sigma+\lambda\tau^{\alpha-\beta}=0
\]
on $M$.
\end{theorem}

\smallskip

\noindent{\em Proof.} The result follows on the set $M\setminus\overline{\Lambda}_1$ by \cite{McGillivray2018_1} Theorem 6.5 {\em (i)} and follows on the set $\overline{\Lambda}_1$ by {\em (ii)}-{\em (iii)} because $\widetilde{\varrho}=(\zeta\psi)^\prime/(\zeta\psi)=g^\prime/g=\beta/\tau$ is continuous. The $\mathscr{H}^1$-a.e. statement in this last can be dropped as $M$ is analytic. 
\qed

\smallskip

\noindent Let $E$ be an open set with $C^1$ boundary $M$ in $\mathbb{R}^2_0$ and assume that $E=E^{sc}$. Introduce the projection $\pi:\mathbb{R}^2_0\rightarrow(0,\infty);x\mapsto|x|$. The set
\begin{equation}\label{definition_of_Omega}
\Omega:=\pi\Big[M\cap\{\cos\sigma\neq 0\}\Big]
\end{equation}
plays an important r\^{o}le in the sequel. 
Bearing in mind \cite{McGillivray2018} Lemma 5.4 we may define
\begin{align}
&\theta_2:\Omega\rightarrow(0,\pi);\tau\mapsto L(\tau)/2\tau;\label{definition_of_theta}\\
&\gamma:\Omega\rightarrow M;\tau\mapsto(\tau\cos\theta_2(\tau),\tau\sin\theta_2(\tau)).\label{definition_of_gamma}
\end{align}
The function 
\begin{align}
&u:\Omega\rightarrow[-1,1];\tau\mapsto\sin(\sigma(\gamma(\tau))).\label{definition_of_y}
\end{align}
is closely related to the geodesic curvature of $M$ (in fact $k=(1/\tau)(\tau u)^\prime$) and we may rephrase the last Theorem \ref{constant_weighted_mean_curvature} as follows.

\smallskip

\begin{theorem}\label{ode_for_y}
Let $(\alpha,\beta)\in\mathcal{Q}$. Given $v>0$ let $E$ be a minimiser of (\ref{isoperimetric_problem}). Assume that $E$ is a bounded open set with analytic boundary $M$ in $\mathbb{R}^2_0$ and suppose that $E=E^{sc}$. Suppose that $\Omega\neq\emptyset$. Then $u\in C^{1}(\Omega)$ and
\begin{align}
u^\prime+\frac{\beta+1}{\tau}u+\lambda\tau^{\alpha-\beta}&=0\label{ode_for_u}
\end{align}
on $\Omega$.
\end{theorem}

\smallskip

\noindent{\em Proof.} See \cite{McGillivray2018_1} Theorem 6.6.
\qed

\smallskip

\noindent Finally we recount \cite{McGillivray2018_1} Theorem 6.7 which gives an explicit expression for the derivative of the angular coordinate $\theta_2$ of a boundary point involving the function $u$.  

\smallskip

\begin{lemma}\label{derivative_of_theta}
Suppose that $E$ is a bounded open set with $C^1$ boundary $M$ in $\mathbb{R}^2_0$ and that $E=E^{sc}$.  Suppose that $\Omega\neq\emptyset$. Then 
\begin{itemize}
\item[(i)] $\theta_2\in C^1(\Omega)$;
\item[(ii)] $\theta_2^\prime=-\frac{1}{\tau}\frac{u}{\sqrt{1-u^2}}$ on $\Omega$.
\end{itemize}
\end{lemma}

\smallskip
\section{Some algebraic inequalities}

\smallskip

\noindent In this Section we obtain two algebraic inequalities. These are used in the next Section to derive some integral inequalities: Theorem \ref{main_algebraic_inequality} is used in the proof of Corollary \ref{integral_of_w} while Proposition \ref{inequality_for_hat_m} is used in the proof of Corollary \ref{integral_of_solution_of_first_order_ode}. Let us fix a pair of indices $(\alpha,\beta)\in\mathcal{P}$ for the time being. The definition of $\mathcal{P}$ entails that $\beta+1<\alpha+2\leq 2(\beta+1)$ and in particular $\beta+1>0$ and $\alpha+2>\beta+1>0$. Given positive numbers $a$ and $b$ with $a<b$ define
\begin{align}
m=m(a,b)=m_{\alpha,\beta}(a,b)
&:=(\alpha+2)\frac{b^{\beta+1}-a^{\beta+1}}{b^{\alpha+2}-a^{\alpha+2}}.\label{definition_of_m}
\end{align}
We shall make frequent use of the homogeneity property
\begin{align}
m(a,b)&=a^{-(\alpha-\beta+1)}m(t)\label{homogeneity}
\end{align}
where $t:=b/a>1$ and $m(t):=m(1,t)$. The bulk of this Section will be devoted to establishing the following algebraic inequality.

\smallskip

\begin{theorem}\label{main_algebraic_inequality}
Fix $(\alpha,\beta)\in\mathcal{P}$. Then for any $0<a<b<\infty$,
\begin{align}
&\frac{1}{1+\beta-ma^{\alpha-\beta+1}}-\frac{1}{1+\beta-mb^{\alpha-\beta+1}}\geq 2\frac{a+b}{b-a}.\label{basic_inequality_for_w_problem}
\end{align}
If $(\alpha,\beta)\in\mathcal{P}\setminus\{(0,0)\}$ then strict inequality holds in (\ref{basic_inequality_for_w_problem}) for each $0<a<b<\infty$. If $(\alpha,\beta)=(0,0)$ then equality holds in (\ref{basic_inequality_for_w_problem}) for each $0<a<b<\infty$.
\end{theorem}

\smallskip

\noindent Let us first mention some landmarks in the proof of this result. We first define the set
\begin{align}
\mathcal{P}^+&:=\Big\{(\alpha,\beta)\in\mathbb{R}^2:\alpha\geq 0,\beta\geq 0\text{ and }\alpha(\beta+1)=\beta^2\Big\}.\nonumber
\end{align}
A first step in the proof of Theorem \ref{main_algebraic_inequality} is to show that it holds for indices $(\alpha,\beta)$ belonging to $\mathcal{P}^+\setminus\{(0,0)\}$. With this in hand a monotonicity argument is then used to derive the result for indices in the set $\mathcal{P}\setminus\mathcal{P}^-$. A separate argument deals with the case $\mathcal{P}^-$. To begin we require some preliminary lemmas.

\smallskip

\begin{lemma}\label{Holder_inequality_lemma}
For $\alpha>\beta>0$ the function
\[
(1,\infty)\rightarrow\mathbb{R};t\mapsto\frac{t^{\beta}-1}{t^\alpha-1}
\]
is strictly decreasing.
\end{lemma}

\smallskip

\noindent{\em Proof.}
The derivative of this function is given by
\begin{align}
\frac{d}{dt}\frac{t^\beta-1}{t^\alpha-1}&=\frac{(\beta-\alpha)t^{\alpha+\beta-1}+\alpha t^{\alpha-1}-\beta t^{\beta-1}}{(t^\alpha-1)^2}\nonumber
\end{align}
for $t>1$. Put
\[
p:=\frac{\alpha}{\alpha-\beta}\text{ and }q:=\frac{\alpha}{\beta}
\]
so that $p>1$, $q>1$ and $1/p+1/q=1$. For $t>1$ set
\[
a:=t^{\frac{\alpha^2-\beta^2}{\alpha}}\text{ and }b:=t^{\frac{\beta^2}{\alpha}}
\]
and apply Young's inequality to obtain
\begin{align}
t^{\alpha}&=ab\leq\frac{a^p}{p}+\frac{b^q}{q}
=\frac{\alpha-\beta}{\alpha}t^{\alpha+\beta}+\frac{\beta}{\alpha}t^\beta
\nonumber
\end{align}
and in fact strict inequality holds as the equality condition in Young's inequality reads $a^p-b^q=t^{\alpha+\beta}-t^\beta=t^\beta(t^\alpha-1)>0$. This entails that the derivative is strictly negative.
\qed

\smallskip

\begin{lemma}\label{properties_of_m}
Let $(\alpha,\beta)\in\mathcal{P}$, $0<a<b<\infty$ and define $m$ as in (\ref{definition_of_m}). Then
\begin{itemize}
\item[(i)] $a^{\alpha-\beta+1}m<\beta+1$;
\item[(ii)] $b^{\alpha-\beta+1}m>\beta+1$.
\end{itemize}
Moreover for each fixed $a>0$,
\begin{itemize}
\item[(iii)] the mapping $b\mapsto m(a,b)$ on $(a,\infty)$ is strictly decreasing;
\item[(iv)] the mapping $b\mapsto b^{\alpha-\beta+1}m(a,b)$ on $(a,\infty)$ is strictly increasing.
\end{itemize}
\end{lemma}

\smallskip

\noindent{\em Proof.} By homogeneity it suffices to consider the case $a=1$ and $b=t>1$. For {\em (i)} we require
\[
\frac{\alpha+2}{\beta+1}\frac{t^{\beta+1}-1}{t^{\alpha+2}-1}
=\frac{\int_1^t\tau^{\beta}\,d\tau}{\int_1^t\tau^{\alpha+1}\,d\tau}
<1
\]
for each $t>1$ which is the case as $\alpha+1>\beta$. By Cauchy's mean-value theorem,
\begin{align}
t^{\alpha-\beta+1}\frac{\alpha+2}{\beta+1}\frac{t^{\beta+1}-1}{t^{\alpha+2}-1}&=\Big(\frac{t}{c}\Big)^{\alpha-\beta+1}>1\nonumber
\end{align}
for some $c\in(1,t)$ and hence {\em (ii)}. {\em (iii)} This follows using homogeneity and Lemma \ref{Holder_inequality_lemma}. {\em (iv)} Using homogeneity we write
\begin{align}
t^{\alpha-\beta+1}m(1,t)&=(\alpha+2)\frac{t^{\alpha+2}-t^{\alpha-\beta+1}}{t^{\alpha+2}-1}
=(\alpha+2)\Big\{1-\frac{t^{\alpha-\beta+1}-1}{t^{\alpha+2}-1}\Big\}
\nonumber
\end{align}
and the claim follows from Lemma \ref{Holder_inequality_lemma}.
\qed

\smallskip

\noindent We are now in a position to rephrase the statement of Theorem \ref{main_algebraic_inequality} for indices belonging to $\mathcal{P}^+$. The condition mentioned in the next proposition is verified in Lemma \ref{inequality_on_boundary_of_P} leading to Corollary \ref{positivity_on_P_+}.

\smallskip

\begin{proposition}\label{inequality_on_boundary_of_P}
Let $(\alpha,\beta)\in\mathcal{P}^+\setminus\{(0,0)\}$. The following are equivalent:
\begin{itemize}
\item[(i)] for any $0<a<b<\infty$,
\begin{align}
&\frac{1}{1+\beta-ma^{\alpha-\beta+1}}-\frac{1}{1+\beta-mb^{\alpha-\beta+1}}>2\frac{a+b}{b-a};
\nonumber
\end{align}
\item[(ii)] for each $\lambda>0$,
\begin{align}
(1/x)^2\coth(\lambda/2x)
+x^2\coth(\lambda x/2)
-y^2\coth(\lambda y/2)
+y\tanh(\lambda/2)&>0
\label{hyperbolic_inequality_on_boundary_of_P}
\end{align}
where $x:=\beta+1$ and $y:=\alpha+2$. 
\end{itemize}
\end{proposition}

\smallskip

\noindent{\em Proof.}
By homogeneity it suffices to consider the case $a=1$ and $b=t>1$ in which case the inequality (\ref{basic_inequality_for_w_problem}) becomes
\begin{align}
&\frac{1}{1+\beta-m}-\frac{1}{1+\beta-mt^{\alpha-\beta+1}}\geq 2\frac{t+1}{t-1}\label{modified_inequality_for_w_problem}
\end{align}
with now
\begin{align}
m=m(t)&:=(\alpha+2)\frac{t^{\beta+1}-1}{t^{\alpha+2}-1}\label{formula_for_m_in_terms_of_t}
\end{align}
as in (\ref{definition_of_m}). Bearing in mind Lemma \ref{properties_of_m} {\em (i)} and {\em (ii)} the inequality in (\ref{modified_inequality_for_w_problem}) may be rewritten
\begin{align}
&2(t+1)t^{\alpha-\beta+1}m^2+\Big\{(t-1)(t^{\alpha-\beta+1}-1)-2(\beta+1)(t+1)(t^{\alpha-\beta+1}+1)\Big\}m\nonumber\\
&+2(\beta+1)^2(t+1)\geq 0\nonumber
\end{align}
as a quadratic in $m$. Write $\gamma:=\alpha-\beta+1$ and replace $m$ with the expression in (\ref{formula_for_m_in_terms_of_t}) to obtain
\begin{align}
&2(\alpha+2)^2(t+1)t^{\gamma}(t^{\beta+1}-1)^2
+2(\beta+1)^2(t+1)(t^{\alpha+2}-1)^2\nonumber\\
&+(\alpha+2)(t^{\alpha+2}-1)(t^{\beta+1}-1)\Big\{(t-1)(t^{\gamma}-1)-2(\beta+1)(t+1)(t^{\gamma}+1)\Big\}
\geq 0.\label{intermediate_inequality_1}
\end{align}

\smallskip 

\noindent Let us introduce some further notation. Put $x:=\beta+1$. Because of the constraint $\alpha(\beta+1)=\beta^2$,
\[
y:=\alpha+2=\frac{\beta^2}{\beta+1}+2=\beta+1+\frac{1}{\beta+1}=x+\frac{1}{x}.
\]
In this case $\gamma=y-x=1/x$. In this notation the inequality (\ref{intermediate_inequality_1}) can be rewritten
\begin{align}
&2y^2(t+1)t^{1/x}(t^x-1)^2+2x^2(t+1)(t^{y}-1)^2\nonumber\\
&+y(t^{y}-1)(t^x-1)\Big\{(t-1)(t^{1/x}-1)-2x(t+1)(t^{1/x}+1)\Big\}\geq 0.\nonumber
\end{align}
As $t>1$ we may put $t=e^\lambda$ for some $\lambda>0$. With this further substitution we obtain
\begin{align}
&2y^2(e^\lambda+1)e^{\lambda/x}(e^{\lambda x}-1)^2
+2x^2(e^\lambda+1)(e^{\lambda y}-1)^2\nonumber\\
&+y(e^{\lambda y}-1)(e^{\lambda x}-1)\Big\{(e^\lambda-1)(e^{\lambda/x}-1)-2x(e^\lambda+1)(e^{\lambda/x}+1)\Big\}\geq 0.\nonumber
\end{align}
Multiplication by $e^{-\lambda x}$ and $e^{-\lambda/x}$ in turn (or $e^{-\lambda y}$) results in the inequality
\begin{align}
&2y^2(e^\lambda+1)(e^{\lambda x/2}-e^{-\lambda x/2})^2
+2x^2(e^\lambda+1)(e^{\lambda y/2}-e^{-\lambda y/2})^2\nonumber\\
&+y(e^{\lambda y/2}-e^{-\lambda y/2})(e^{\lambda x/2}-e^{-\lambda x/2})\Big\{(e^\lambda-1)(e^{\lambda/2x}-e^{-\lambda/2x})\nonumber\\
&-2x(e^\lambda+1)(e^{\lambda/2x}+e^{-\lambda/2x})\Big\}\geq 0\nonumber
\end{align}
or more briefly
\begin{align}
&y^2\cosh(\lambda/2)\sinh^2(\lambda x/2)
+x^2\cosh(\lambda/2)\sinh^2(\lambda y/2)\nonumber\\
&+y\sinh(\lambda y/2)\sinh(\lambda x/2)
\Big\{
\sinh(\lambda/2)\sinh(\lambda/2x)
-2x\cosh(\lambda/2)\cosh(\lambda/2x)\Big\}\geq 0.\nonumber
\end{align}

\smallskip

\noindent To continue let us introduce the function 
\[
\phi(x):=\frac{\sinh(\lambda x/2)}{x}
\]
defined for positive $x$. The above expression can be written in terms of this function as
\[
\phi(x)^2+\phi(y)^2-2\cosh(\lambda/2x)\phi(x)\phi(y)+\phi(x)\phi(y)\tanh(\lambda/2)\frac{\sinh(\lambda/2x)}{x}\geq 0
\]
after dividing by $(xy)^2\cosh(\lambda/2)$. Or equivalently,
\begin{align}
\frac{\phi(x)}{\phi(y)}+\frac{\phi(y)}{\phi(x)}-2\cosh(\lambda/2x)+\tanh(\lambda/2)\frac{\sinh(\lambda/2x)}{x}&\geq0.\label{intermediate_inequality}
\end{align}
By the addition formula (cf. \cite{AbramowitzStegun1964} 4.5.24 for example),
\[
\frac{\phi(y)}{\phi(x)}
=\frac{\sinh(\lambda y/2)}{\sinh(\lambda x/2)}\frac{x}{y}
=\Big\{\cosh(\lambda/2x)+\coth(\lambda x/2)\sinh(\lambda/2x)\Big\}
\frac{x}{y}
\]
and its reciprocal is given by
\[
\frac{\phi(x)}{\phi(y)}
=\frac{\sinh(\lambda x/2)}{\sinh(\lambda y/2)}\frac{y}{x}
=\Big\{\cosh(\lambda/2x)-\coth(\lambda y/2)\sinh(\lambda/2x)\Big\}
\frac{y}{x}.
\]
Inserting these expressions into (\ref{intermediate_inequality}) leads to
\begin{align}
&\cosh(\lambda/2x)\Big(\frac{x}{y}+\frac{y}{x}-2\Big)
+\sinh(\lambda/2x)\Big(\coth(\lambda x/2)\frac{x}{y}-\coth(\lambda y/2)\frac{y}{x}
+\frac{\tanh(\lambda/2)}{x}\Big)\geq 0;\nonumber
\end{align}
or upon dividing by $\sinh(\lambda/2x)$ and multiplying by $xy$,
\begin{align}
&(1/x)^2\coth(\lambda/2x)
+x^2\coth(\lambda x/2)
-y^2\coth(\lambda y/2)
+y\tanh(\lambda/2)\geq 0.
\nonumber
\end{align}
\qed

\smallskip

\begin{lemma}\label{inequality_on_boundary_of_P}
Let $x$ and $y$ be positive numbers. Then the inequality (\ref{hyperbolic_inequality_on_boundary_of_P}) holds for each $\lambda>0$.
\end{lemma}

\smallskip

\noindent{\em Proof.} We first write the left-hand side in (\ref{hyperbolic_inequality_on_boundary_of_P}) as
\begin{align}
&x^2\Big\{\coth(\lambda x/2)-\coth(\lambda y/2)\Big\}+(1/x)^2\Big\{\coth(\lambda/2x)-\coth(\lambda y/2)\Big\}\nonumber\\
&-2\coth(\lambda y/2)+y\tanh(\lambda/2)\label{alternative_left_hand_side}
\end{align}
after expanding $y^2$. By the addition formula (cf. \cite{AbramowitzStegun1964} 4.5.27 for example),
\begin{align}
\coth(\lambda x/2)-\coth(\lambda y/2)&=\coth(\lambda x/2)-\frac{1+\coth(\lambda x/2)\coth(\lambda/2x)}{\coth(\lambda x/2)+\coth(\lambda/2x)}\nonumber\\
&=\frac{\coth(\lambda x/2)[\coth(\lambda x/2)+\coth(\lambda/2x)]-1-\coth(\lambda x/2)\coth(\lambda/2x)}{\coth(\lambda x/2)+\coth(\lambda/2x)}\nonumber\\
&=\frac{\coth^2(\lambda x/2)-1}{\coth(\lambda x/2)+\coth(\lambda/2x)}
=\frac{\sinh^{-2}(\lambda x/2)}{\coth(\lambda x/2)+\coth(\lambda/2x)}\nonumber
\end{align}
and likewise
\begin{align}
\coth(\lambda/2x)-\coth(\lambda y/2)&
=\frac{\sinh^{-2}(\lambda/2x)}{\coth(\lambda x/2)+\coth(\lambda/2x)}.\nonumber
\end{align}
The expression (\ref{alternative_left_hand_side})  may be rewritten by completing the square using the addition formula (cf. \cite{AbramowitzStegun1964} 4.5.24 for example),
\begin{align}
&\frac{(x/\sinh(\lambda x/2))^2+(1/x\sinh(\lambda/2x))^2}{\coth(\lambda x/2)+\coth(\lambda/2x)}
-2\coth(\lambda y/2)+y\tanh(\lambda/2)\nonumber
\end{align}
\[
=
\frac{(x/\sinh(\lambda x/2)-1/x\sinh(\lambda/2x))^2}{\coth(\lambda x/2)+\coth(\lambda/2x)}
+\frac{2[\sinh(\lambda x/2)\sinh(\lambda/2x)]^{-1}}{\coth(\lambda x/2)+\coth(\lambda/2x)}
-2\coth(\lambda y/2)+y\tanh(\lambda/2)
\]
The final three terms in the expression above can be expressed
\[
\frac{2}{\sinh(\lambda y/2)}
-2\coth(\lambda y/2)+y\tanh(\lambda/2)
=2\frac{1-\cosh(\lambda y/2)}{\sinh(\lambda y/2)}+y\tanh(\lambda/2)
\]
\[
=-2\tanh(\lambda y/4)+y\tanh(\lambda/2)
=
2\Big\{-\tanh(\lambda y/4)+(y/2)\tanh(\lambda/2)\Big\}
\]
where we used the half-angle formula (cf. \cite{AbramowitzStegun1964} 4.5.30 for example) in the penultimate line. This means that (\ref{alternative_left_hand_side}) can be written
\[
\frac{(x/\sinh(\lambda x/2)-1/x\sinh(\lambda/2x))^2}{\coth(\lambda x/2)+\coth(\lambda/2x)}
+2\Big\{-\tanh(\lambda y/4)+(y/2)\tanh(\lambda/2)\Big\}.
\]

\smallskip

\noindent We now claim that
\[
-\tanh(\lambda y/4)+(y/2)\tanh(\lambda/2)\geq 0\text{ for any }y\geq 2.
\]
At $y=2$ equality holds. Differentiating the left-hand side gives
\begin{align}
-\frac{\lambda}{4}\frac{1}{\cosh^2(\lambda y/4)}+(1/2)\tanh(\lambda/2)
&\geq-\frac{\lambda}{4}\frac{1}{\cosh^2(\lambda/2)}+(1/2)\tanh(\lambda/2)\nonumber\\
&=\frac{1}{4\cosh^2(\lambda/2)}\Big\{-\lambda+2\sinh(\lambda/2)\cosh(\lambda/2)\Big\}\nonumber\\
&=\frac{1}{4\cosh^2(\lambda/2)}\Big\{-\lambda+\sinh(\lambda)\Big\}\nonumber
\end{align}
by the double-angle formula (cf. \cite{AbramowitzStegun1964} 4.5.31 for example) and
this last is positive for $\lambda>0$.
\qed

\smallskip

\begin{corollary}\label{positivity_on_P_+}
Let $(\alpha,\beta)\in\mathcal{P}^+\setminus\{(0,0)\}$. Then for any $0<a<b<\infty$,
\begin{align}
&\frac{1}{1+\beta-ma^{\alpha-\beta+1}}-\frac{1}{1+\beta-mb^{\alpha-\beta+1}}>2\frac{a+b}{b-a}.
\nonumber
\end{align}
\end{corollary}

\medskip

\noindent Given $(\alpha,\beta)\in\mathcal{P}$ set
\begin{align}
\zeta&:=\frac{\beta+1}{\alpha-\beta+1}\label{zeta}
\end{align}
and note that $\zeta\geq 1$ in virtue of the property that $\alpha\leq 2\beta$. As before we put $\gamma:=\alpha-\beta+1$ and consider this as a function defined on $\mathcal{P}$. For fixed $\zeta\neq 0$ we introduce the line $\ell_\zeta$ in the $(\alpha,\beta)$-plane given by
\begin{align}
&\ell_\zeta:\mathbb{R}\rightarrow\mathbb{R};\beta\mapsto\ell_\zeta(\beta):=(1+1/\zeta)\beta+1/\zeta-1.
\label{line}
\end{align}
The graph of this line passes through the point $(-2,-1)$ as well as $(\alpha,\beta)$. 

\smallskip

\noindent We have seen that strict inequality in (\ref{basic_inequality_for_w_problem}) holds on $\mathcal{P}^+$. To obtain the inequality in $\mathcal{P}$ we consider the line passing through $(\alpha,\beta)$ in $\mathcal{P}$ and the point $(-2,-1)$. This meets $\mathcal{P}^+$ in a unique point $(\overline{\alpha},\overline{\beta})$. A monotonicity argument allows us to obtain the inequality at $(\alpha,\beta)$ by making use of the inequality at $(\overline{\alpha},\overline{\beta})$. See Figure \ref{fig:Figure_lines}. We now describe the details of this argument. Let us start by fixing some properties of the lines $\ell_\zeta$.

\smallskip

\begin{figure}[h]
\centering
\begin{tikzpicture}[>=stealth]
    \begin{axis}[
        xmin=-2.5,xmax=3.5,
        ymin=-1.5,ymax=3.5,
        axis x line=middle,
        axis y line=middle,
        axis line style=->,
        xlabel={$\alpha$},
        ylabel={$\beta$},
        ]
        \addplot[black, thick, domain=0:3, smooth]{(x+sqrt(x^2+4*x))/2};
        \addplot[black, domain=-2:3.5, smooth]{x/2};
        \addplot[black, domain=-2:3.5, smooth]{x+1};
        \addplot[black, dashed, domain=-2:3.5, smooth]{0.75*x+0.5};
        \addplot[smooth,mark=*,black] plot coordinates {(-0.5,0.125)};
    \end{axis}
\end{tikzpicture}
\caption{The region $\mathcal{P}$ is enclosed between the lines $\alpha-\beta+1=0$ and $\alpha-2\beta=0$ and the curve $\alpha(\beta+1)-\beta^2=0$. This latter curve is marked in bold and corresponds to $\mathcal{P}^+$. Suppose that $(\alpha,\beta)\in\mathcal{P}$. The line $\ell_\zeta$ that passes through $(-2,-1)$ and $(\alpha,\beta)$ intersects $\mathcal{P}^+$ in a unique point $(\overline{\alpha},\overline{\beta})$. The gradient of the line $\ell_\zeta$ when considered as a function of $\beta$ is given by $\frac{\alpha+2}{\beta+1}$ and has a value in the interval $(1,2]$.}
\label{fig:Figure_lines}
\end{figure}
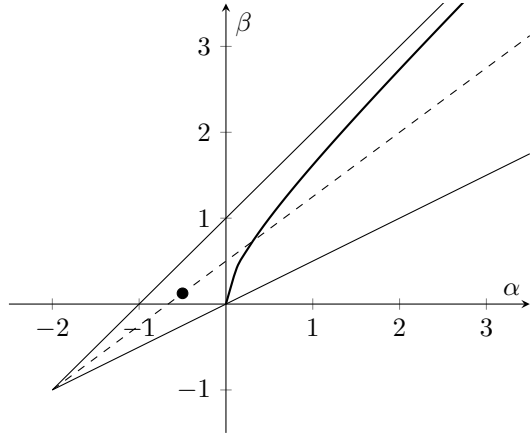

\smallskip

\begin{lemma}\label{lines_and_set_P}
For each point $(\alpha,\beta)\in\mathcal{P}$,
\begin{itemize}
\item[(i)] there exists a unique $\zeta\geq 1$ with the property that $(\alpha,\beta)$ lies on the graph of $\ell_\zeta$;
\item[(ii)] $\zeta$ is given by (\ref{zeta});
\end{itemize}
and for each $\zeta\in[1,\infty)$,
\begin{itemize}
\item[(iii)] the graph of the line $\ell_\zeta$ meets $\mathcal{P}^+$ exactly once;
\item[(iv)] the function $\beta\mapsto\gamma(\ell_\zeta(\beta),\beta)$ is increasing on $\mathbb{R}$.
\end{itemize}
\end{lemma}

\smallskip

\noindent{\em Proof.}
{\em (i)} and {\em (ii)}. For $(\alpha,\beta)\in\mathcal{P}$ define $\zeta$ as in (\ref{zeta}) and observe that
\begin{align}
1+\frac{1}{\zeta}&=\frac{\alpha+2}{\beta+1}\in(1,2];\label{zeta_alpha_and_beta}
\end{align}
this is the gradient of the line $\ell_\zeta$ when considered as a function of $\beta$. So $\ell_\zeta(\beta)=\alpha$ and $(\ell_\zeta(\beta),\beta)=(\alpha,\beta)$; that is, the point $(\alpha,\beta)$ lies on the graph of $\ell_\zeta$. Conversely, if $\ell_\zeta(\beta)=\alpha$ then $\zeta$ is given by (\ref{zeta}) and $\zeta\geq 1$ as $(\alpha,\beta)\in\mathcal{P}$.

\smallskip

\noindent{\em (iii)} It may be helpful to refer to Figure \ref{fig:Figure_hyperbola}. Define an affine transformation $T$ mapping the $(\alpha,\beta)$-plane into the $(X,Y)$-plane via
\begin{align}
T&:\mathbb{R}^2\rightarrow\mathbb{R}^2;(\alpha,\beta)\mapsto(X,Y)=
\Big(\frac{\alpha+2}{2},\frac{2\beta-\alpha}{2}\Big).\nonumber
\end{align}
This has inverse $\alpha=2X-2$ and $\beta=Y+X-1$. Note that
\[
X^2-Y^2=\alpha(\beta+1)-\beta^2+1.
\]
The image of $\mathcal{P}^+$ under $T$ is the branch of the hyperbola
\begin{align}
X^2-Y^2&=1\label{hyperbola}
\end{align}
in the first quadrant of the $(X,Y)$-plane noting that $X=(\alpha/2)+1$ has range $[1,\infty)$ when restricted to $\mathcal{P}^+$. The centre $(0,0)$ of the hyperbola (\ref{hyperbola}) in the $(X,Y)$-plane corresponds to the point $(-2,-1)$ in the $(\alpha,\beta)$-plane; its vertices correspond to the points $(0,0)$ and $(-4,-2)$; and its asymptotes $X\pm Y=0$ to the lines
\[
\alpha-\beta+1=0\text{ and }\beta+1=0.
\]
Note also that the line $X=1$ corresponds to the line $\alpha=0$ and the line $\alpha-2\beta=0$ corresponds to the line $Y=0$. The graph of the line $\ell_\zeta$ in (\ref{line}) corresponds to the line
\begin{align}
Y&=\frac{\zeta-1}{\zeta+1}X\nonumber
\end{align}
in the $(X,Y)$-plane. Notice that the gradient varies between $0$ and $1$ as $\zeta$ varies between $1$ and $\infty$. This makes evident the claim in {\em (iii)}. 

\begin{figure}[h]
\centering
\begin{tikzpicture}[>=stealth]
    \begin{axis}[
        xmin=-0.5,xmax=4.5,
        ymin=-0.5,ymax=4.5,
        axis x line=middle,
        axis y line=middle,
        axis line style=->,
        xlabel={$X$},
        ylabel={$Y$},
        ]
        \addplot[black, thick, domain=1:4.5, smooth]{sqrt(x^2-1)};
        \addplot[black, domain=0:4.5, smooth]{x};
    \end{axis}
\end{tikzpicture}
\caption{The image of $\mathcal{P}$ under the mapping $T$ is enclosed between the lines $Y=0$ and $Y=X$ and the arc of the hyperbola $X^2-Y^2-1=0$ marked in bold. This latter is the image of $\mathcal{P}^+$ under $T$. }
\label{fig:Figure_hyperbola}
\end{figure}
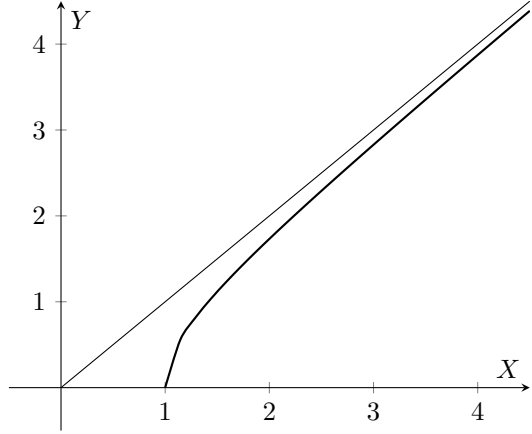

\smallskip

\noindent{\em (iv)} For $\beta\in\mathbb{R}$,
\begin{align}
\gamma(\ell_\zeta(\beta),\beta)&=\ell_\zeta(\beta)-\beta+1=
(1+1/\zeta)\beta+1/\zeta-1-\beta+1=(1/\zeta)(\beta+1)\nonumber
\end{align}
which increases as $\beta$ increases.
\qed

\smallskip

\noindent For $\zeta\geq 1$ define
\begin{align}
M:[1,\infty)\rightarrow\mathbb{R};\tau\mapsto
\left\{
\begin{array}{lc}
1& \text{ for }\tau=1;\\
\frac{\zeta+1}{\zeta}\frac{\tau^\zeta-1}{\tau^{\zeta+1}-1} & \text{ for }\tau>1;\\
\end{array}
\right.
\label{definition_of_M}
\end{align}
so that $M$ is continuous and strictly increasing on $[1,\infty)$ with limit $0$ at $\infty$ by Lemma \ref{Holder_inequality_lemma}. Define
\begin{align}
\Lambda:(1,\infty)\rightarrow\mathbb{R};\tau\mapsto\frac{1}{1-M(\tau)}-\frac{1}{1-\tau M(\tau)}.\label{definition_of_Lambda}
\end{align}

\smallskip

\begin{lemma}\label{Lambda_is_decreasing}
Let $(\alpha,\beta)\in\mathcal{P}$. Define $\Lambda$ as in (\ref{definition_of_Lambda}). Then $\Lambda$ is strictly decreasing on $(1,\infty)$.
\end{lemma}

\smallskip

\noindent{\em Proof.}
For $\tau>1$,
\[
\tau M(\tau)=\frac{\zeta+1}{\zeta}\Big\{1-\frac{\tau-1}{\tau^{\zeta+1}-1}\Big\}
\]
This is strictly increasing by Lemma \ref{Holder_inequality_lemma} while $M$ is strictly decreasing as already remarked. Note that $1-\tau M(\tau)<0$ for $\tau>1$ as this is equivalent to $\tau^{\zeta+1}-(\zeta+1)\tau+\zeta>0$ and this holds by Young's inequality because $a^p-b^q=\tau^{\zeta+1}-1>0$.
\qed

\smallskip

\noindent{\em Proof of Theorem \ref{main_algebraic_inequality}.}
Let $(\alpha,\beta)\in\mathcal{P}\setminus\mathcal{P}_-$. Define $\zeta$ as in (\ref{zeta}). It suffices to take $a=1$ and $b=t>1$ by homogeneity. Put $\tau:=t^{\alpha-\beta+1}=t^\gamma$. Then $t^{\beta+1}=\tau^\zeta$ and $t^{\alpha+2}=\tau^{\zeta+1}$. So we may write
\begin{align}
m(t)&=(\alpha+2)\frac{t^{\beta+1}-1}{t^{\alpha+2}-1}=(\beta+1)\frac{\zeta+1}{\zeta}\frac{\tau^\zeta-1}{\tau^{\zeta+1}-1}
=(\beta+1)M(\tau)\nonumber
\end{align}
making use of (\ref{zeta_alpha_and_beta}). This leads to in turn to the identity
\begin{align}
\frac{1}{\beta+1-m}-\frac{1}{\beta+1-t^{\alpha-\beta+1}m}
&=\frac{1}{\beta+1}\Big\{\frac{1}{1-M(\tau)}-\frac{1}{1-\tau M(\tau)}\Big\}
=\frac{1}{\beta+1}\Lambda(\tau).\label{identity_for_Lambda}
\end{align}
According to Lemma \ref{lines_and_set_P} the graph of the line $\ell_\zeta$ meets the curve $\mathcal{P}^+\setminus\{(0,0)\}$ at a unique point $(\overline{\alpha},\overline{\beta})$ with $\overline{\alpha}=\ell_\zeta(\overline{\beta})$ where $\beta\leq\overline{\beta}$. By (\ref{identity_for_Lambda}) and Corollary \ref{positivity_on_P_+},
\[
\frac{1}{\overline{\beta}+1}\Lambda(t^{\overline{\gamma}})=
\frac{1}{\overline{\beta}+1-m_{\overline{\alpha},\overline{\beta}}(t)}
-\frac{1}{\overline{\beta}+1-t^{\overline{\gamma}}m_{\overline{\alpha},\overline{\beta}}(t)}
>2\frac{1+t}{t-1}
\] 
for any $t>1$. Again by (\ref{identity_for_Lambda}) and Lemma \ref{Lambda_is_decreasing},
\begin{align}
\frac{1}{\beta+1-m}-\frac{1}{\beta+1-t^{\alpha-\beta+1}m}
&=\frac{1}{\beta+1}\Lambda(t^\gamma)
\geq\frac{1}{\overline{\beta}+1}\Lambda(t^{
\overline{\gamma}
}
)>2\frac{1+t}{t-1}
.\nonumber
\end{align}
where $\overline{\gamma}=\gamma(\ell_\zeta(\overline{\beta}),\overline{\beta})$ because $\gamma\leq\overline{\gamma}$ by Lemma \ref{lines_and_set_P} and $t>1$. 

\smallskip

\noindent Now suppose that $(\alpha,\beta)\in\mathcal{P}^-\setminus\{(0,0)\}$. In this case $-1<\beta<0$, $\alpha=2\beta$ and
\[
m(t)=\frac{2(\beta+1)}{1+t^{\beta+1}}
\]
for $t>1$. Moreover,
\begin{align}
\frac{1}{1+\beta-m(t)}-\frac{1}{1+\beta-m(t)t^{\alpha-\beta+1}}
&=\frac{2}{\beta+1}\frac{t^{\beta+1}+1}{t^{\beta+1}-1}.
\end{align}
Strict inequality holds in (\ref{basic_inequality_for_w_problem}) if and only if
\[
\beta[t^{\beta+2}-1]+(\beta+2)t[t^\beta-1]<0
\]
and this holds as $-1<\beta<0$ and $t>1$. Equality holds if $(\alpha,\beta)=(0,0)$. 
\qed

\bigskip

\noindent For $(\alpha,\beta)\in\mathcal{P}$ and $0<a<b<\infty$ define
\begin{align}
\hat{m}&=\hat{m}(a,b)=\hat{m}_{\alpha,\beta}(a,b):=(\alpha+2)\frac{b^{\beta+1}+a^{\beta+1}}{b^{\alpha+2}-a^{\alpha+2}}.\label{definition_of_hat_m}
\end{align}
The corresponding homogeneity property reads
\begin{align}
\hat{m}(a,b)&=a^{-(\alpha-\beta+1)}\hat{m}(t)\label{homogeneity}
\end{align}
where $t:=b/a>1$ and $\hat{m}(t):=\hat{m}(1,t)$. 

\smallskip

\begin{proposition}\label{inequality_for_hat_m}
Let $(\alpha,\beta)\in\mathcal{P}\setminus\mathcal{P}^-$ and $0<a<b<\infty$ and define $\hat{m}$ as in (\ref{definition_of_hat_m}). Then
\begin{align}
\hat{m}&<2\frac{\beta+1}{b^{\alpha-\beta+1}-a^{\alpha-\beta+1}}.\nonumber
\end{align}
\end{proposition}

\smallskip

\noindent{\em Proof.}
By homogeneity it suffices to show that
\[
(\alpha+2)\frac{t^{\beta+1}+1}{t^{\alpha+2}-1}<2\frac{\beta+1}{t^{\alpha-\beta+1}-1}
\]
for each $t>1$. This is equivalent to the inequality
\[
\frac{\alpha+2}{\beta+1}\Big\{t^{\alpha+2}-1+t^{\alpha-\beta+1}-t^{\beta+1}\Big\}
<2(t^{\alpha+2}-1).
\]
Now $\alpha<2\beta$ as $(\alpha,\beta)\not\in\mathcal{P}^-$ so $\frac{\alpha+2}{\beta+1}<2$ and $0<\alpha-\beta+1<\beta+1$. This proves the result. \qed

\section{Some integral inequalities}

\smallskip

\noindent Recall that the generalised curvature of the boundary of a spherical cap isoperimetric minimiser can be expressed in terms of the function $u$ as in Theorem \ref{ode_for_y}. In this Section we study solutions to the differential equation (\ref{ode_for_u}) subject to two types of Dirichlet boundary condition on a bounded open interval in $(0,\infty)$. The function $u$ determines the angular coordinate as in Lemma \ref{derivative_of_theta}. We obtain estimates for these angular integrals. The main results of this Section are Corollary \ref{integral_of_solution_of_first_order_ode}, Corollary \ref{integral_of_w} and Corollary \ref{comparision_result_for_u_in_case_a_is_zero}.

\smallskip

\noindent Let $\mathscr{L}$ stand for the collection of Lebesgue measurable sets in $[0,\infty)$. Define a measure $\mu$ on $([0,\infty),\mathscr{L})$ by $\mu(dx):=(1/x)\,dx$. Let $0\leq a<b<\infty$. Suppose that $u:[a,b]\rightarrow\mathbb{R}$ is an $\mathscr{L}^1$-measurable function with the property that
\begin{align}
\mu(\{u>t\})&<\infty\text{ for each }t>0.
\end{align}
 The distribution function  $\mu_u:(0,\infty)\rightarrow[0,\infty)$ of $u$ with respect to $\mu$ is given by
\[
\mu_u(t):=\mu(\{u>t\})\text{ for }t>0.
\]
Note that $\mu_u$ is right-continuous and non-increasing on $(0,\infty)$ and $\mu_u(t)\rightarrow 0$ as $t\rightarrow\infty$. 

\smallskip

\noindent Let $0<a<b<\infty$ and $(\alpha,\beta)\in\mathcal{P}$. Let $\eta\in\{\pm 1\}^2$. We study solutions $(u,\lambda)$ to the first-order linear ordinary differential equation
\begin{align}
u^\prime+(\beta+1)\tau^{-1}u+\lambda\tau^{\alpha-\beta}&=0
\text{ on }(a,b)\text{ with }u(a)=\eta_1\text{ and }u(b)=\eta_2\label{ode_for_u_with_plus_minus_1_at_endpoints}
\end{align}
where $u\in C^1([a,b])$ and $\lambda\in\mathbb{R}$. In case $(\alpha,\beta)=(0,0)$ we use the notation $u_0$. This latter is given explicitly by
\begin{align}
u_0(\tau)&=\frac{1}{b-a}\Big\{\tau-\frac{ab}{\tau}\Big\}\label{formula_for_u_0}
\end{align}
for $\tau\in[a,b]$ and $u_0$ is strictly increasing on $[a,b]$.

\smallskip

\begin{lemma}\label{first_order_ode}
Let $(\alpha,\beta)\in\mathcal{P}$ and $0<a<b<\infty$. Let $\eta\in\{\pm 1\}^2$. Then
\begin{itemize}
\item[(i)] there exists a solution $(u,\lambda)$ of (\ref{ode_for_u_with_plus_minus_1_at_endpoints}) with $u$ real analytic and $\lambda=\lambda_\eta\in\mathbb{R}$;
\item[(ii)] the pair $(u,\lambda)$ in (i) is unique;
\item[(iii)] $\lambda_\eta$ is given by
\begin{align}
&-\lambda_{(1,1)}=\lambda_{(-1,-1)}=m;\,\lambda_{(1,-1)}=-\lambda_{(-1,1)}=\hat{m};\nonumber
\end{align}
with $m$ resp. $\hat{m}$ as in (\ref{definition_of_m}) resp. (\ref{definition_of_hat_m});
\item[(iv)] if $\eta=(1,1)$ or $\eta=(-1,-1)$ then $u$ is uniformly bounded away from zero on $[a,b]$.
\end{itemize}
\end{lemma}

\smallskip

\noindent{\em Proof.}
{\em (i)} Multiplying by an integrating factor we may rewrite the ordinary differential equation in (\ref{ode_for_u_with_plus_minus_1_at_endpoints}) in the form $(\tau^{\beta+1}u)^\prime+\lambda\tau^{\alpha+1}=0$ and hence
\begin{align}
&u=\tau^{-\beta-1}\Big\{c-\frac{\lambda}{\alpha+2}\tau^{\alpha+2}\Big\}\label{general_equation_for_u}
\end{align}
for real constants $c$ and $\lambda$. Fitting the boundary conditions we derive
\begin{align}
\lambda&=-m=-(\alpha+2)\frac{b^{\beta+1}-a^{\beta+1}}{b^{\alpha+2}-a^{\alpha+2}};\label{lambda_for_eta_equal_to_11}\\
c&=-\frac{m}{\alpha+2}a^{\alpha+2}+a^{\beta+1};\label{c_for_eta_equal_to_11}
\end{align}
in case $\eta=(1,1)$ while
\begin{align}
\lambda&=-\hat{m}=-(\alpha+2)\frac{a^{\beta+1}+b^{\beta+1}}{b^{\alpha+2}-a^{\alpha+2}};\label{lambda_for_eta_equal_to_-11}\\
c&=-\frac{\hat{m}}{\alpha+2}b^{\alpha+2}+b^{\beta+1};;\label{c_for_eta_equal_to_-11}
\end{align}
in case $\eta=(-1,1)$. Similar expressions hold in cases $\eta=(-1,-1)$ and $\eta=(1,-1)$ by linearity. {\em (ii)} We consider the case $\eta=(1,1)$. Suppose that $(u_1,\lambda_1)$ resp. $(u_2,\lambda_2)$ solve  (\ref{ode_for_u_with_plus_minus_1_at_endpoints}). By linearity $u:=u_1-u_2$ solves
\[
u^\prime+(\beta+1)\tau^{-1}u+\lambda\tau^{\alpha-\beta}=0
\text{ on }(a,b)
\text{ with }u(a)=u(b)=0
\]
where $\lambda=\lambda_1-\lambda_2$. So $u$ takes the form (\ref{general_equation_for_u}) for some real constants $c$ and $\lambda$ and the boundary conditions entail that these vanish. The other cases are similar. {\em (iii)} follows as in {\em (i)}. {\em (iv)} If $\eta=(1,1)$ then $u>0$ on $[a,b]$ as
$c+\frac{m}{\alpha+2}a^{\alpha+2}=a^{\beta+1}>0$.
\qed

\medskip

\noindent{\em The boundary condition $\eta_1\eta_2=-1$.}

\smallskip

\begin{lemma}\label{derivative_of_distribution_function_for_y_resp_z}
Let $(\alpha,\beta)\in\mathcal{P}\setminus\mathcal{P}^-$ and $0<a<b<\infty$. Let $(u,\lambda)$ solve (\ref{ode_for_u_with_plus_minus_1_at_endpoints}) with $\eta=(-1,1)$. Put $v:=-u$. Then
\begin{itemize}
\item[(i)] $u^\prime>0$ on $[a,b]$;
\item[(ii)]
$\int_{\{v=1\}}\frac{1}{|v^\prime|}\frac{d\mathscr{H}^0}{\tau}
<
\int_{\{u=1\}}\frac{1}{|u^\prime|}\frac{d\mathscr{H}^0}{\tau}$.
\end{itemize}
\end{lemma}

\smallskip

\noindent{\em Proof.} {\em (i)} Write the derivative as $u^\prime=-(\beta+1)\tau^{-1}u+\hat{m}\tau^{\alpha-\beta}$ and replace $u$ with the expression in (\ref{general_equation_for_u}). Strict positivity of the derivative is equivalent to the condition
\[
-b^{\alpha+2}\frac{\hat{m}}{\alpha+2}+b^{\beta+1}<\frac{\alpha-\beta+1}{(\alpha+2)(\beta+1)}\tau^{\alpha+2}\hat{m}
\]
on substituting the expression for $u$ in (\ref{general_equation_for_u}), (\ref{lambda_for_eta_equal_to_-11}) and (\ref{c_for_eta_equal_to_-11}). After inserting the definition (\ref{definition_of_hat_m}) for $\hat{m}$ the left-hand side above becomes
\[
-\frac{(ab)^{\beta+1}}{b^{\alpha+2}-a^{\alpha+2}}\Big(a^{\alpha-\beta+1}+b^{\alpha-\beta+1}\Big)
\]
which is plainly negative. This leads to the statement.

\smallskip

\noindent{\em (ii)} From (\ref{ode_for_u_with_plus_minus_1_at_endpoints}),
\[
\tau u^\prime=-(\beta+1)u+\hat{m}\tau^{\alpha-\beta+1}
\]
and in particular,
\begin{align}
bu^\prime(b)&=-(\beta+1)+\hat{m}b^{\alpha-\beta+1}
>-(\beta+1)+m b^{\alpha-\beta+1}
>0\nonumber
\end{align}
as $\hat{m}>m$ and by Lemma \ref{properties_of_m}. On the other hand,
\begin{align}
av^\prime(a)&=-au^\prime(a)=-\beta-1-a^{\alpha-\beta+1}\hat{m}<0.\nonumber
\end{align}
Then
\begin{align}
\int_{\{v=1\}}\frac{1}{|v^\prime|}\frac{d\mathscr{H}^0}{\tau}
-
\int_{\{u=1\}}\frac{1}{|u^\prime|}\frac{d\mathscr{H}^0}{\tau}
&=
\frac{1}{au^\prime(a)}-\frac{1}{bu^\prime(b)}\nonumber\\
&=\frac{1}{\beta+1+a^{\alpha-\beta+1}\hat{m}}
-\frac{1}{-\beta-1+\hat{m}b^{\alpha-\beta+1}}\nonumber\\
&=\frac{-2(\beta+1)+\hat{m}(b^{\alpha-\beta+1}-a^{\alpha-\beta+1})}{
(\beta+1+\hat{m}a^{\alpha-\beta+1})(-\beta-1+\hat{m}b^{\alpha-\beta+1})
}<0\nonumber
\end{align}
by Proposition \ref{inequality_for_hat_m}.
\qed

\smallskip

\begin{theorem}\label{comparision_of_derivatives_of_distribution_functions}
Let $(\alpha,\beta)\in\mathcal{P}\setminus\mathcal{P}^-$ and $0<a<b<\infty$. Let $(u,\lambda)$ solve (\ref{ode_for_u_with_plus_minus_1_at_endpoints}) with $\eta=(-1,1)$. Put $v:=-u$.
Then
\begin{itemize}
\item[(i)]  $-\mu_u^\prime(t)=\int_{\{u=t\}}\frac{1}{|u^\prime|}\frac{d\mathscr{H}^0}{\tau}$ for each $t\in(0,1)$;
\item[(ii)] $-\mu_v^\prime<-\mu_u^\prime$ on $(0,1)$.
\end{itemize}
\end{theorem}

\smallskip

\noindent{\em Proof.}
{\em (i)} As $u^\prime>0$ on $[a,b]$ according to Lemma \ref{derivative_of_distribution_function_for_y_resp_z} the function $u$ possesses a $C^1$ inverse $u^{-1}:[-1,1]\rightarrow[a,b]$. This allows us to write
\begin{align}
\mu_u(t)&=\mu(\{u>t\})=\mu((u^{-1}(t),b])=\log\Big(\frac{b}{u^{-1}(t)}\Big)\nonumber
\end{align}
for $t>1$ and further
\begin{align}
\mu_u^\prime(t)&=-\frac{1}{(\tau u^\prime)(u^{-1}(t))}=-\int_{\{u=t\}}\frac{1}{|u^\prime|}\frac{d\mathscr{H}^0}{\tau}.
\nonumber
\end{align}
{\em (ii)} Let $t\in(0,1)$. Put $c:=\min\{v\leq t\}=\min\{u\geq-t\}$ and $d:=\min\{u\geq t\}$. By continuity, $a<c<d<b$. Put $\widetilde{u}:=u/t$ and $\widetilde{v}:=v/t$ on $[c,d]$. Then
\begin{align}
\widetilde{u}^\prime+(\beta+1)\tau^{-1}\widetilde{u}-(\hat{m}/t)\tau^{\alpha-\beta}&=0\text{ on }(c,d)\text{ and }-\widetilde{u}(c)=\widetilde{u}(d)=1;\nonumber\\
\widetilde{v}^\prime+(\beta+1)\tau^{-1}\widetilde{v}+(\hat{m}/t)\tau^{\alpha-\beta}&=0\text{ on }(c,d)\text{ and }\widetilde{v}(c)=-\widetilde{v}(d)=1.\nonumber
\end{align}
By Lemma \ref{derivative_of_distribution_function_for_y_resp_z},
\begin{align}
\int_{\{v=t\}}\frac{1}{|v^\prime|}\frac{d\mathscr{H}^0}{\tau}
&=\int_{[c,d]\cap\{v=t\}}\frac{1}{|v^\prime|}\frac{d\mathscr{H}^0}{\tau}
=(1/t)\int_{[c,d]\cap\{\widetilde{v}=1\}}\frac{1}{|\widetilde{v}^\prime|}\frac{d\mathscr{H}^0}{\tau}\nonumber\\
&<(1/t)\int_{[c,d]\cap\{\widetilde{u}=1\}}\frac{1}{|\widetilde{u}^\prime|}\frac{d\mathscr{H}^0}{\tau}
=\int_{\{u=t\}}\frac{1}{|u^\prime|}\frac{d\mathscr{H}^0}{\tau}.\nonumber
\end{align}
The claim now follows with the help of {\em (i)}.
\qed

\smallskip

\begin{corollary}\label{inequality_for_mu_u_vs_mu_v}
Let $(\alpha,\beta)\in\mathcal{P}\setminus\mathcal{P}^-$ and $0<a<b<\infty$. Let $(u,\lambda)$ solve (\ref{ode_for_u_with_plus_minus_1_at_endpoints}) with $\eta=(-1,1)$. Put $v:=-u$. Then 
\[
\mu_u(t)>\mu_v(t)$ for each $t\in(0,1).
\]
\end{corollary}

\smallskip

\noindent{\em Proof.}
By Theorem \ref{comparision_of_derivatives_of_distribution_functions},
\[
\mu_u(t)=\mu_u(t)-\mu_u(1)=-\int_{(t,1]}\mu_u^\prime\,ds
\]
for each $t\in(0,1)$ as $\mu_u(1)=0$. On the other hand,
\[
\mu_v(t)=\mu_v(1)+(\mu_v(t)-\mu_v(1))
=\mu_v(1)-\int_{(t,1]}\mu_v^\prime\,ds=-\int_{(t,1]}\mu_v^\prime\,ds
\]
for each $t\in(0,1)$. The claim follows from Theorem \ref{comparision_of_derivatives_of_distribution_functions}.
\qed

\smallskip

\begin{corollary}\label{integral_of_solution_of_first_order_ode}
Let $(\alpha,\beta)\in\mathcal{P}\setminus\mathcal{P}^-$ and $0<a<b<\infty$. Let $(u,\lambda)$ solve (\ref{ode_for_u_with_plus_minus_1_at_endpoints}) with $\eta=(-1,1)$. Let $\varphi\in C^1((-1,1))$ be an odd strictly increasing function with $\varphi\in L^1((-1,1))$. Then
\begin{itemize}
\item[(i)] $\int_{\{u>0\}}\varphi(u)\,d\mu<\infty$;
\item[(ii)] $\int_a^b\varphi(u)\,d\mu>0$;
\end{itemize}
and in particular, 
\begin{itemize}
\item[(iii)] $\int_a^b\frac{u}{\sqrt{1-u^2}}\,d\mu>0$.
\end{itemize}
\end{corollary}

\smallskip

\noindent{\em Proof.}
{\em (i)} The function $u$ is $C^1$ and $u^\prime$ is bounded away from $0$ on $[a,b]$ by Lemma \ref{derivative_of_distribution_function_for_y_resp_z}. Put $I:=\{0<u<1\}$. The restriction of $u$ to the interval $I$ has $C^1$ inverse $v:(0,1)\rightarrow I$ with derivative $v^\prime=1/(u^\prime\circ v)$ and $|v^\prime|$ is uniformly bounded on $(0,1)$. By a change of variables, 
\begin{align}
\int_{\{u>0\}}\varphi(u)\,d\mu&=-\int_0^1\varphi(v^\prime/v)\,dt\nonumber
\end{align}
from which the claim is apparent.
{\em (ii)} The integral is well-defined because $\varphi(u)^+=\varphi(u)\chi_{\{u>0\}}\in L^1((a,b),\mu)$ by {\em (i)}. The set $\{u=0\}$ consists of a singleton and has $\mu$-measure zero. So
\begin{align}
\int_a^b\varphi(u)\,d\mu&=\int_{\{u>0\}}\varphi(u)\,d\mu+\int_{\{u<0\}}\varphi(u)\,d\mu
=\int_{\{u>0\}}\varphi(u)\,d\mu-\int_{\{v>0\}}\varphi(v)\,d\mu\nonumber
\end{align}
where $v:=-u$ as $\varphi$ is an odd function. We remark that in a similar way as above,
\[
\int_0^1\varphi^\prime\mu_u\,dt=\int_{\{u>0\}}\Big\{\varphi(u)-\varphi(0)\Big\}\,d\mu
=\int_{\{u>0\}}\varphi(u)\,d\mu
\]
using oddness of $\varphi$ and an analogous formula holds with $v$ in place of $u$. Thus we may write
\begin{align}
\int_a^b\varphi(u)\,d\mu&=\int_0^1\varphi^\prime\mu_u\,dt-\int_0^1\varphi^\prime\mu_v\,dt
=\int_0^1\varphi^\prime\Big\{\mu_u-\mu_v\Big\}\,dt>0\nonumber
\end{align}
by Corollary \ref{inequality_for_mu_u_vs_mu_v} as $\varphi^\prime>0$ on $(0,1)$. {\em (iii)} follows from {\em (ii)} with the particular choice $\varphi:(-1,1)\rightarrow\mathbb{R};t\mapsto t/\sqrt{1-t^2}$.
\qed

\bigskip

\noindent{\em The boundary condition $\eta_1\eta_2=1$.} Let $(\alpha,\beta)\in\mathcal{P}$ and $0<a<b<\infty$. In this subsection study solutions $(w,\lambda)$ of the auxilliary Riccati equation 
\begin{align}
w^\prime+\lambda\tau^{\alpha-\beta}w^2&=\frac{\beta+1}{\tau}w
\text{ on }(a,b)\text{ with }w(a)=w(b)=1;
\label{Riccati_equation}
\end{align}
with $w$ a $C^1$ function on $[a,b]$ and $\lambda\in\mathbb{R}$. If $(\alpha,\beta)=(0,0)$ then we write $w_0$ instead of $w$. Suppose $(u,\lambda)$ solves (\ref{ode_for_u_with_plus_minus_1_at_endpoints}) with $\eta=(1,1)$. Then $u>0$ on $[a,b]$ by Lemma \ref{first_order_ode} and we may set $w:=1/u$. Then $(w,-\lambda)$ satisfies (\ref{Riccati_equation}).

\smallskip

\begin{lemma}\label{solution_of_Riccati_equation}
Let $0<a<b<\infty$ and $(\alpha,\beta)\in\mathcal{P}$. Then
\begin{itemize}
\item[(i)] there exists a solution $(w,\lambda)$ of (\ref{Riccati_equation}) with $w$ real analytic and
$\lambda\in\mathbb{R}$;
\item[(ii)] the pair $(w,\lambda)$ in (i) is unique;
\item[(iii)] $\lambda=m$ with $m$ as in (\ref{definition_of_m}).
\end{itemize}
\end{lemma}

\smallskip

\noindent{\em Proof.}
{\em (i)} With $c$ as in (\ref{c_for_eta_equal_to_11}),
\[
c+\frac{m}{\alpha+2}\tau^{\alpha+2}\geq c+\frac{m}{\alpha+2}a^{\alpha+2}=a^{\beta+1}>0
\]
for $\tau\in[a,b]$. It therefore makes sense to define $w:[a,b]\rightarrow\mathbb{R}$ by
\begin{align}
w&:=\frac{\tau^{\beta+1}}{c+\frac{m}{\alpha+2}\tau^{\alpha+2}}.\label{formula_for_w}
\end{align}
Then $w$ is real analytic on $[a,b]$ and $(w,m)$ satisfies (\ref{Riccati_equation}). {\em (ii)} We claim that $w>0$ on $[a,b]$ for any solution $(w,\lambda)$ of (\ref{Riccati_equation}). For otherwise, $c:=\min\{w=0\}\in(a,b)$. Then $u:=1/w$ on $[a,c)$ satisfies
\[
u^\prime+\frac{\beta+1}{\tau}u+\lambda\tau^{\alpha-\beta}=0
\text{ on }(a,c)
\text{ and }u(a)=1, u(c-)=\infty.
\]
Integrating, we obtain
\[
\tau^{\beta+1}u-a^{\beta+1}+\frac{\lambda}{\alpha+2}\Big[\tau^{\alpha+2}-a^{\alpha+2}\Big]=0\text{ on }[a,c)
\]
and this entails the contradiction that $u(c-)<\infty$. We may now use the uniqueness statement in Lemma \ref{first_order_ode}. {\em (iii)} follows from {\em (ii)} and the particular solution given in {\em (i)}. 
\qed

\smallskip

\begin{lemma}\label{lemma_on_derivative_of_mu_v}
Let $(\alpha,\beta)\in\mathcal{P}\setminus\mathcal{P}^-$ and $0<a<b<\infty$. Let $(w,\lambda)$ solve (\ref{Riccati_equation}). Then\begin{itemize}
\item[(i)] $w^\prime(a)>0$ and $w^\prime(b)<0$;
\item[(ii)] $w>1$ on $(a,b)$;
\item[(iii)] $w$ is unimodal;
\item[(iv)] $\int_{\{w=1\}}\frac{1}{|w^\prime|}\frac{d\mathscr{H}^0}{\tau}>2\coth(\mu((a,b))/2).
$
\end{itemize}
\end{lemma}

\smallskip

\noindent{\em Proof.}
{\em (i)} follows from Lemma \ref{properties_of_m}. {\em (ii)} The inequality $w(\tau)>1$ may be rewritten in the form $m_{\alpha,\beta}(a,\tau)>m_{\alpha,\beta}(a,b)=m$ upon making use of the expression for $w$ given in (\ref{formula_for_w}). This last inequality is a consequence of Lemma \ref{properties_of_m}. {\em (iii)} The reciprocal $u:=1/w$ satisfies (\ref{ode_for_u_with_plus_minus_1_at_endpoints}) with $\eta=(1,1)$. The derivative of $u$ vanishes at a point $\tau$ characterised by the condition
\[
\tau^{\alpha+2}=\frac{(\alpha+2)(\beta+1)}{\alpha-\beta+1}\frac{c}{m}.
\]
It can be shown using Lemma \ref{properties_of_m} that $\tau\in(a,b)$. It follows that $w$ is unimodal.
{\em (iv)} By Theorem \ref{main_algebraic_inequality},
\begin{align}
\int_{\{w=1\}}\frac{1}{|w^\prime|}\frac{d\mathscr{H}^0}{\tau}
&=\frac{1}{1+\beta-ma^{\alpha-\beta+1}}-\frac{1}{1+\beta-mb^{\alpha-\beta+1}}
>2\frac{a+b}{b-a}=2\coth(\mu((a,b))/2).\nonumber
\end{align}
\qed

\smallskip

\begin{theorem}\label{derivative_of_Ricatti_distribution_function}
Let $(\alpha,\beta)\in\mathcal{P}\setminus\mathcal{P}^-$ and $0<a<b<\infty$ . Let $(w,\lambda)$ solve (\ref{Riccati_equation}). Then for each $t\in(1,\|w\|_\infty)$,
\begin{itemize}
\item[(i)] $-\mu_w^\prime=\int_{\{w=t\}}\frac{1}{|w^\prime|}\frac{d\mathscr{H}^0}{\tau}$;
\item[(ii)] $-\mu_w^\prime>(2/t)\coth((1/2)\mu_w)$.
\end{itemize}
\end{theorem}

\smallskip

\noindent{\em Proof.}
{\em (i)} This follows from the unimodality property of $w$ in Lemma \ref{lemma_on_derivative_of_mu_v}. Let $c$ be the unique point in $(a,b)$ such that $w^\prime(c)=0$. The $C^1$ functions $w_1:[a,c]\rightarrow[1,\|w\|_\infty]$ resp. $w_2:[c,b]\rightarrow[1,\|w\|_\infty]$ have $C^1$ inverses $v_1:[1,\|w\|_\infty]\rightarrow[a,c]$ resp. $v_2:[1,\|w\|_\infty]\rightarrow[c,b]$. Then
\[
\mu_w(t)=\mu((v_1(t),v_2(t))=\log\Big(\frac{v_2(t)}{v_1(t)}\Big)
\]
and
\[
\mu_w^\prime(t)=\frac{v_2^\prime(t)}{v_2(t)}-\frac{v_1^\prime(t)}{v_1(t)}
=\frac{1}{(\tau w)^\prime(v_2(t))}-\frac{1}{(\tau w)^\prime(v_1(t))}
=-\int_{\{w=t\}}\frac{1}{|w^\prime|}\frac{d\mathscr{H}^0}{\tau}
\]
for each $t\in(1,\|w\|_\infty)$.
{\em (ii)} Let $t\in(1,\|w\|_\infty)$. Put $\widetilde{w}:=w/t$ on $\overline{\{w>t\}}$ so
\begin{align}
\widetilde{w}^\prime+(mt)\tau^{\alpha-\beta}\widetilde{w}^2&=\frac{\beta+1}{\tau}\widetilde{w}
\text{ on }\{w>t\}\text{ and }\widetilde{w}=1\text{ on }\{w=t\}.\nonumber
\end{align}
By Lemma \ref{lemma_on_derivative_of_mu_v},
\begin{align}
(0,\infty)\ni\int_{\{w=t\}}\frac{1}{|w^\prime|}\frac{d\mathscr{H}^0}{\tau}
&=(1/t)\int_{\{\widetilde{w}=1\}}\frac{1}{|\widetilde{w}^\prime|}\frac{d\mathscr{H}^0}{\tau}\nonumber\\
&>(2/t)\coth((1/2)\mu(\{\widetilde{w}>1\}))
=(2/t)\coth((1/2)\mu_w(t)).\nonumber
\end{align}
The statement now follows from {\em (i)}.
\qed

\smallskip

\noindent We introduce the mapping
\[
\omega:(0,\infty)\times(0,\infty)\rightarrow\mathbb{R};(t,x)\mapsto-(2/t)\coth(x/2).
\]
For $\xi>0$,
\begin{align}
|\omega(t,x)-\omega(t,y)|&\leq\mathrm{cosech}^2[\xi/2](1/t)|x-y|\label{Lipschitz_estimate_for_omega}
\end{align}
for $(t,x),(t,y)\in(0,\infty)\times(\xi,\infty)$ and $\omega$ is locally Lipschitzian in $x$ on $(0,\infty)\times(0,\infty)$ in the sense of \cite{Hale1969} I.3.
Let $0<a<b<\infty$ and set $\lambda:=A/G>1$. Here, $A=A(a,b)$ stands for the arithmetic mean of $a,b$ as introduced in the previous Section while $G=G(a,b):=\sqrt{|ab|}$ stands for their geometric mean. We refer to the inital value problem
\begin{equation}\label{auxilliary_ode_1}
z^\prime=\omega(t,z)\text{ on }(0,\lambda)
\text{ and }z(1)=\mu((a,b)).
\end{equation}
Define
\[
z_0:(0,\lambda)\rightarrow\mathbb{R};t\mapsto
2\log\Big\{\frac{\lambda+\sqrt{\lambda^2-t^2}}{t}\Big\}.
\]

\smallskip

\begin{lemma}\label{properties_of_w_0}
Let $0<a<b<\infty$.  Then
\begin{itemize}
\item[(i)] $w_0(\tau)=\frac{2A\tau}{G^2+\tau^2}$ for $\tau\in[a,b]$;
\item[(ii)] $\|w_0\|_\infty=\lambda$;
\item[(iii)] $\mu_{w_0}=z_0$ on $[1,\lambda)$;
\item[(iv)] $z_0$ satisfies (\ref{auxilliary_ode_1}) and this solution is unique;
\item[(v)] $\int_{\{w_0=1\}}\frac{1}{|w_0^\prime|}\,\frac{d\mathscr{H}^0}{\tau}=2\coth(\mu((a,b))/2)$;
\item[(vi)] $\int_a^b\frac{1}{\sqrt{w_0^2-1}}\,\frac{dx}{x}=\pi$. 
\end{itemize}
\end{lemma}

\smallskip

\noindent{\em Proof.} {\em (i)} Note that 
\[
m_0=\frac{1+ab}{a+b}.
\]
{\em (i)} follows from the representation of $w_0$ in Lemma \ref{solution_of_Riccati_equation} by direct computation. {\em (ii)} follows by calculus. {\em (iii)} follows by solving the quadratic equation $t\tau^2-2A\tau+G^2t=0$ for   $\tau$ with $t\in(0,\lambda)$. Uniqueness in {\em (iv)} follows from \cite{Hale1969} Theorem 3.1 as $\omega$ is locally Lipschitzian with respect to $x$ in $(0,\infty)\times(0,\infty)$. For {\em (v)} note that $|aw_0^\prime(a)|=1-a/A$ and $|bw_0^\prime(b)|=b/A-1$ and 
\[
2\coth(\mu((a,b))/2)=2(a+b)/(b-a).
\]
{\em (vi)} We may write
\begin{align}
\int_a^b\frac{1}{\sqrt{w_0^2-1}}\,\frac{d\tau}{\tau}
&=\int_a^b\frac{ab+\tau^2}{\sqrt{(a+b)^2\tau^2-(ab+\tau^2)^2}}
\frac{d\tau}{\tau}=\int_a^b\frac{ab+\tau^2}{\sqrt{(\tau^2-a^2)(b^2-\tau^2)}}
\frac{d\tau}{\tau}.\nonumber
\end{align}
The substitution $s=\tau^2$ followed by the Euler substitution (cf. \cite{Gradshteynetal1965} 2.251)
$\sqrt{(s-a^2)(b^2-s)}=t(s-a^2)$ gives
\begin{align}
\int_a^b\frac{1}{\sqrt{w_0^2-1}}\,\frac{d\tau}{\tau}
&=\int_0^\infty\frac{1}{1+t^2}+\frac{ab}{b^2+a^2t^2}\,dt=\pi.\nonumber
\end{align}
\qed

\smallskip

\begin{theorem}\label{distribution_function_inequality_for_Ricatti_equation}
Let $(\alpha,\beta)\in\mathcal{P}\setminus\mathcal{P}^-$ and $0<a<b<\infty$. Let $(w,\lambda)$ solve (\ref{Riccati_equation}). Then 
\begin{itemize}
\item[(i)] $\|w\|_\infty\leq\|w_0\|_\infty$;
\item[(ii)] $\mu_w(t)<\mu_{w_0}(t)$ for each $t\in(1,\|w\|_\infty)$.
\end{itemize}
\end{theorem}

\smallskip

\noindent{\em Proof.} {\em (i)} From the inequality in Theorem \ref{derivative_of_Ricatti_distribution_function} we infer that $\mu_w^\prime\leq\omega(t,\mu_w)$ on $[1,\|w\|_\infty]$ and $\mu_w(1)=\mu([a,b])$. In virtue of Lemma \ref{properties_of_w_0}, $\mu_{w_0}^\prime=\omega(t,\mu_{w_0})$ on $[1,\|w_0\|_\infty]$ and $\mu_{w_0}(1)=\mu([a,b])$. Put $T:=\min\{\|w_0\|_\infty,\|w\|_\infty\}$. By \cite{Hale1969} Theorem 6.1 we deduce that $\mu_w\leq\mu_{w_0}$ on $[1,T]$. This implies that $\|w\|_\infty\leq\|w_0\|_\infty$. {\em (ii)} From Theorem \ref{derivative_of_Ricatti_distribution_function} we see that $\mu_w^\prime<\omega(t,\mu_w)\leq\omega(t,\mu_{w_0})=\mu_{w_0}^\prime$ on $(1,\|w\|_\infty)$ using the observation that the function $\omega(t,z)$ is strictly increasing in its second argument. This establishes the second claim. \qed

\smallskip

\begin{corollary}\label{integral_of_w}
Let $(\alpha,\beta)\in\mathcal{P}\setminus\mathcal{P}^-$ and $0<a<b<\infty$. Let $(w,\lambda)$ solve (\ref{Riccati_equation}). Let $0\leq\varphi\in C^1((1,\infty))$ be strictly decreasing with $\int_a^b\varphi(w_0)\,d\mu<\infty$. Then
\begin{itemize}
\item[(i)] $\int_a^b\varphi(w)\,d\mu>\int_a^b\varphi(w_0)\,d\mu$;
\end{itemize}
and in particular, 
\begin{itemize}
\item[(ii)] $\int_a^b\frac{1}{\sqrt{w^2-1}}\,d\mu>\pi$.
\end{itemize}
\end{corollary}

\smallskip

\noindent{\em Proof.}
{\em (i)} Let $\varphi\geq 0$ be a decreasing function on $(1,\infty)$ which is piecewise $C^1$. Suppose that $\varphi(1+)<\infty$. By Tonelli's Theorem,
\begin{align}
\int_{(1,\infty)}\varphi^\prime\mu_w\,ds
&=\int_{(1,\infty)}\varphi^\prime\Big\{\int_{(a,b)}\chi_{\{w>s\}}\,d\mu\Big\}\,ds\nonumber\\
&=\int_{(a,b)}\Big\{\int_{[1,\infty)}\varphi^\prime\chi_{\{w>s\}}\,ds\Big\}\,d\mu\nonumber\\
&=\int_{(a,b)}\Big\{\varphi(w)-\varphi(1)\Big\}\,d\mu
=\int_{(a,b)}\varphi(w)\,d\mu-\varphi(1)\mu((a,b))\nonumber
\end{align}
and a similar identity holds for $\mu_{w_0}$. By Theorem \ref{distribution_function_inequality_for_Ricatti_equation}, $\int_a^b\varphi(w)\,d\mu\geq\int_a^b\varphi(w_0)\,d\mu$. Strict inequality holds if $\varphi$ is nonconstant on $(1,\|w\|_\infty)$.  Now suppose that $\varphi\geq 0$ is a decreasing piecewise $C^1$ function on $(1,\infty)$ with $\int_a^b\varphi(w_0)\,d\mu<\infty$. The inequality holds for the truncated function $\varphi\wedge n$ for each $2\leq n\in\mathbb{N}$. An application of the monotone convergence theorem establishes the result result that
\[
\int_a^b\varphi(w)\,d\mu\geq\int_a^b\varphi(w_0)\,d\mu
\]
for such $\varphi$.

\smallskip

\noindent Let $0\leq\varphi\in C^1((1,\infty))$ be strictly decreasing with $\int_a^b\varphi(w_0)\,d\mu<\infty$. Suppose that equality holds for this $\varphi$. For $c\in(1,\infty)$ put $\varphi_1:=\varphi\vee\varphi(c)-\varphi(c)$ and $\varphi_2:=\varphi\wedge\varphi(c)$. Choose $c\in(1,\|w\|_\infty)$. By the results of the last paragraph,
\[
\int_a^b\varphi_1(w)\,d\mu\geq\int_a^b\varphi_1(w_0)\,d\mu
\text{ while }
\int_a^b\varphi_2(w)\,d\mu>\int_a^b\varphi_2(w_0)\,d\mu.
\]
By linearity,
\begin{align}
\int_a^b\varphi_2(w)\,d\mu&=\int_a^b\varphi_1(w)\,d\mu+\int_a^b\varphi_2(w)\,d\mu\nonumber\\
&>\int_a^b\varphi_1(w_0)\,d\mu+\int_a^b\varphi_2(w_0)\,d\mu=\int_a^b\varphi_2(w_0)\,d\mu.\nonumber
\end{align}
\qed

\smallskip

\noindent{\em The case $a=0$.} Let $0<b<\infty$ and $(\alpha,\beta)\in\mathcal{P}$. We study solutions to the first-order linear ordinary differential equation
\begin{align}
u^\prime+(\beta+1)\tau^{-1}u+\lambda\tau^{\alpha-\beta}&=0
\text{ on }(0,b)\text{ with }u(0)=0\text{ and }u(b)=1\label{first_order_ode_with_bc_and_a_zero}
\end{align}
where $u\in C^2([0,b])$ and $\lambda\in\mathbb{R}$. In case $(\alpha,\beta)=(0,0)$ we use the notation $u_0$.

\smallskip

\begin{lemma}\label{first_order_ode_with_a_zero}
Let $0<b<\infty$ and $(\alpha,\beta)\in\mathcal{P}$. The problem (\ref{first_order_ode_with_bc_and_a_zero}) admits the unique solution
\[
u(\tau):=(\tau/b)^{\alpha-\beta+1}\text{ for }\tau\in[0,b]
\]
and $\lambda:=-\frac{\alpha+2}{b^{\alpha-\beta+1}}$.
\end{lemma}

\smallskip

\noindent{\em Proof.}
It is straightforward to check that the pair $(u,\lambda)$ is a solution to (\ref{first_order_ode_with_bc_and_a_zero}). Let us check uniqueness. Suppose that $(u_1,\lambda_1)$ resp. $(u_2,\lambda_2)$ solve  (\ref{first_order_ode_with_bc_and_a_zero}). By linearity $u:=u_1-u_2$ solves
\[
u^\prime+(\beta+1)\tau^{-1}u+\lambda\tau^{\alpha-\beta}=0
\text{ on }(0,b)\text{ with }u(0)=0\text{ and }u(b)=0
\]
where $\lambda=\lambda_1-\lambda_2$. Making use of an integrating factor the zero function solves this problem uniquely.
\qed

\smallskip

\begin{lemma}\label{integral_of_solution_of_first_order_ode_with_a_equal_to_0}
Let $0<b<\infty$. Then $\int_0^b\frac{u_0}{\sqrt{1-u_0^2}}\,d\mu=\pi/2$. 
\end{lemma}

\smallskip

\noindent{\em Proof.} The integral is elementary as $u_0(t)=t/b$ for $t\in[0,b]$.
\qed

\smallskip

\begin{corollary}\label{comparision_result_for_u_in_case_a_is_zero}
Let $0<b<\infty$ and $(\alpha,\beta)\in\mathcal{P}$. Let $(u,\lambda)$ satisfy (\ref{first_order_ode_with_bc_and_a_zero}). Then 
\[
\int_0^b\frac{u}{\sqrt{1-u^2}}\,d\mu=\frac{\pi}{2\gamma}.
\]
\end{corollary}

\smallskip

\noindent{\em Proof.}
This follows directly from Lemma \ref{first_order_ode_with_a_zero} and Lemma \ref{integral_of_solution_of_first_order_ode_with_a_equal_to_0}. \qed

\section{A beta function inequality}

\smallskip

\noindent Recall from \cite{Gradshteynetal1965} 8.380 that the beta function is defined by 
\begin{align}
B(x,y)&:=\int_0^1 t^{x-1}(1-t)^{y-1}\,dt\label{beta_function}
\end{align}
for $x,y>0$. The goal of the present Section is to prove the following theorem. This is used in (\ref{relation_between_isoperimetric_ratios}) to show that a natural competitor domain (spherical cap symmetric whose boundary has a singularity at the origin) is not an isoperimetric minimiser.

\smallskip

\begin{theorem}\label{beta_function_inequality}
For $x\geq 1/2$,
\[
W(x):=\sqrt{2x}\Big(1+\frac{1}{2x}\Big)^{2x}B(x,1/2)\geq 2\pi
\]
and equality holds only at $x=1/2$.
\end{theorem}

\smallskip

\noindent Define a function
\begin{align}
\phi(x)&:=\frac{1}{x}-\log x+\psi(x).\label{definition_of_phi}
\end{align}
for $x>0$ where $\psi$ stands for the digamma function (cf. \cite{AbramowitzStegun1964} 6.3) and the logarithm refers to the natural logarithm. According to Binet's first formula for the logarithm of the gamma function (cf. \cite{WhittakerWatson1935} 12 31),
\[
\log\Gamma(x)=(x-\frac{1}{2})\log x-x+\frac{1}{2}\log(2\pi)+\int_0^\infty\Big(\frac{1}{2}-\frac{1}{t}+\frac{1}{e^t-1}\Big)\frac{e^{-tx}}{t}\,dt
\]
for any $x>0$. The partial derivative of the integrand with respect to $x$ is integrable on $(0,\infty)$ so differentiation yields the identity
\begin{align}
\phi(x)&=\int_0^\infty\Big(\frac{1}{t}-\frac{1}{e^t-1}\Big)e^{-tx}\,dt\label{integral_representation_of_phi}
\end{align}
for each $x>0$. From this representation it can be seen that the function $\phi$ is positive and completely monotonic; and in particular $\phi$ is convex and strictly decreasing. 

 \smallskip
 
 \noindent Define
\begin{align}
w(x)&:=\phi(x)-\phi(x+1/2)+\log\Big(1+\frac{1}{2x}\Big)-\frac{1}{2x}\label{definition_of_w}
\end{align}
for $x>0$ with $\phi$ as above. A computation reveals that
\[
\frac{W^\prime}{W}=(\log W)^\prime=w
\]
on $(0,\infty)$ so that
\[
(\log W)(x)-(\log W)(1/2)=\int_{1/2}^x w\,dy
\]
and
\[
\frac{W(x)}{W(1/2)}=e^{\int_{1/2}^x w\,dy}
\]
for each $x\geq 1/2$. There are two steps in the proof of Theorem \ref{beta_function_inequality}. We first show  in Proposition \ref{nonnegativity_of_integral_of_w} that the integral in the exponent above is positive for $x\geq 1$. In the second step we show that the function $w$ is positive on $(1/2,1)$ in Proposition \ref{positivity_of_w_on_half_one} (see Figure \ref{fig:graph_of_w}). This gives the conclusion of the Theorem.

\begin{figure}[h]
\centering
\includegraphics[width=80mm]{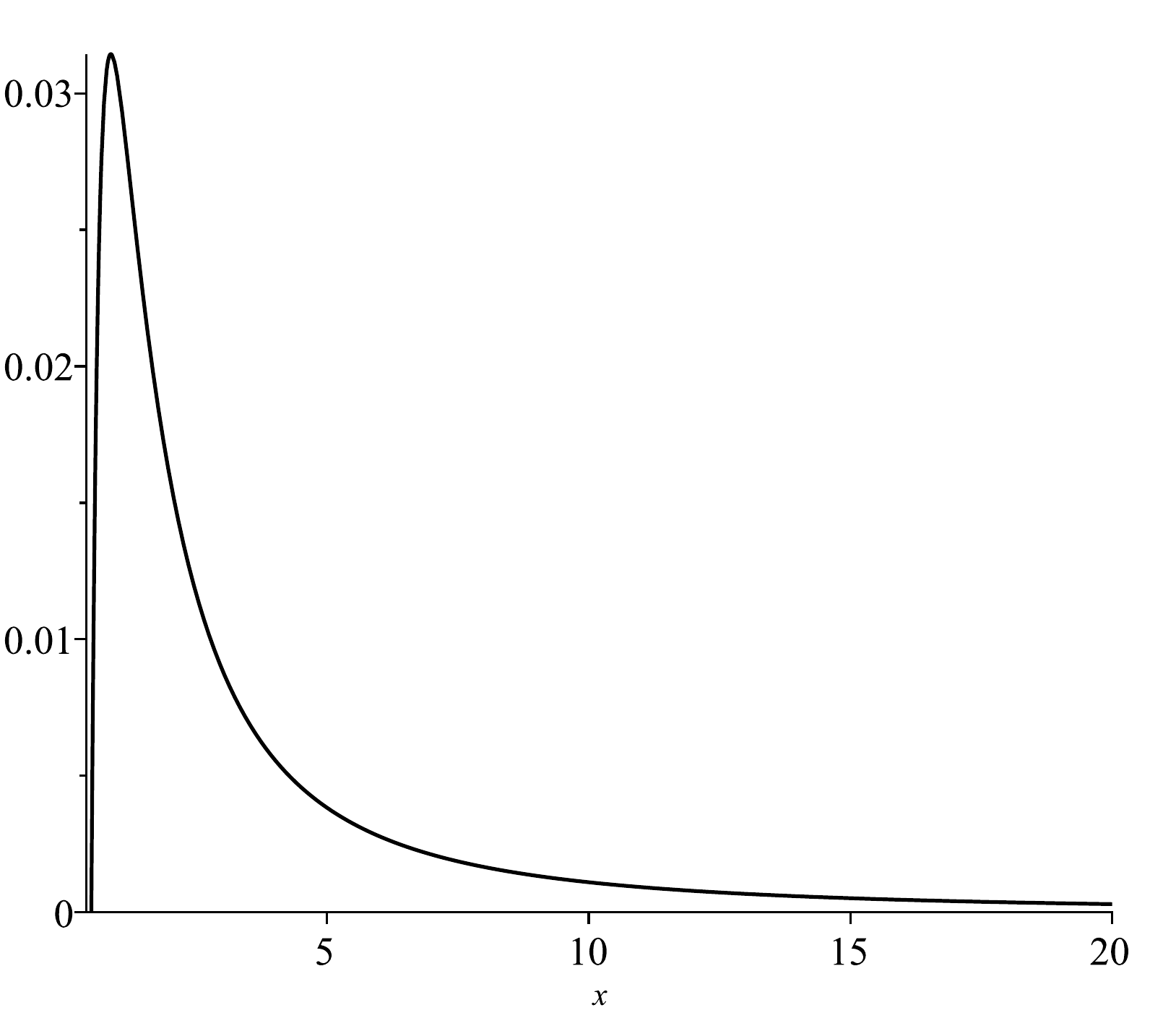}
\caption{Graph of the function $w$ on the interval $[1/2,20]$}
\label{fig:graph_of_w}
\end{figure}

 \smallskip

\begin{lemma}\label{representation_of_integral_of_w}
For $x\geq 1/2$,
\begin{align}
\int_{1/2}^x w\,dy&=\int_0^\infty\frac{1}{t}\Big(e^{-t/2}-e^{-tx}\Big)\Big[\frac{1}{t}-\frac{1}{e^t-1}
+\frac{1}{2}\Big\{\frac{1}{t/2}-\frac{1}{e^{t/2}-1}\Big\}-\frac{1}{2}\Big]\Big(1-e^{-t/2}\Big)\,dt.\nonumber
\end{align}
\end{lemma}

\smallskip

\noindent{\em Proof.} 
Note that
\begin{align}
\int_{1/2}^x e^{-ty}\,dy&=\frac{e^{-t/2}-e^{-tx}}{t}
\nonumber
\end{align}
and
\begin{align}
\int_{1/2}^x e^{-t(y+1/2)}\,dy&=e^{-t/2}\frac{e^{-t/2}-e^{-tx}}{t}
\nonumber
\end{align}
for $1/2\leq x$. From (\ref{integral_representation_of_phi}) and the above,
\begin{align}
\int_{1/2}^x\phi(y)-\phi(y+1/2)\,dy&=\int_0^\infty\Big(\int_{1/2}^x e^{-ty}-e^{-t(y+1/2)}\,dy\Big)\Big(\frac{1}{t}-\frac{1}{e^t-1}\Big)\,dt.\nonumber\\
&=\int_0^\infty\frac{1}{t}\Big(e^{-t/2}-e^{-tx}\Big)\Big(\frac{1}{t}-\frac{1}{e^t-1}\Big)\Big(1-e^{-t/2}\Big)\,dt.\label{integral_involving_phi}
\end{align}
On the other hand for $y>0$,
\begin{align}
\log(y+1/2)-\log y&=\int_y^{y+1/2}\frac{1}{t}\,dt=\int_y^{y+1/2}\int_0^\infty e^{-\tau t}\,d\tau\,dt\nonumber\\
&=\int_0^\infty \int_y^{y+1/2}e^{-\tau t}\,dt\,d\tau
=\int_0^\infty\frac{e^{-\tau y}-e^{-\tau(y+1/2)}}{\tau}\,d\tau\nonumber
\end{align}
We continue
\begin{align}
\int_{1/2}^x\log\Big(1+\frac{1}{2y}\Big)\,dy&=\int_{1/2}^x\log(y+1/2)-\log y\,dy\nonumber\\
&=\int_{1/2}^x
\int_0^\infty\frac{e^{-\tau y}-e^{-\tau(y+1/2)}}{\tau}\,d\tau
\,dy\nonumber\\
&=\int_0^\infty\int_{1/2}^x\frac{e^{-\tau y}-e^{-\tau(y+1/2)}}{\tau}\,dy\,d\tau
\nonumber\\
&=\int_0^\infty\frac{1}{t^2}\Big(e^{-t/2}-e^{-tx}\Big)\Big(1-e^{-t/2}\Big)\,dt.
\label{integral_involving_log}
\end{align}
Finally,
\begin{align}
\int_{1/2}^x\frac{1}{2y}\,dy&=\frac{1}{2}\int_{1/2}^x\int_0^\infty e^{-ty}\,dt\,dy
=\frac{1}{2}\int_0^\infty\int_{1/2}^xe^{-ty}\,dy\,dt
=\frac{1}{2}\int_0^\infty\frac{1}{t}\Big(e^{-t/2}-e^{-tx}\Big),dt
\label{integral_involving_reciprocal_of_2y}
\end{align}
From (\ref{integral_involving_phi}), (\ref{integral_involving_log}) and (\ref{integral_involving_reciprocal_of_2y}),
\begin{align}
\int_{1/2}^x w\,dy&=\int_{1/2}^x\phi(y)-\phi(y+1/2)+\log\Big(1+\frac{1}{2y}\Big)-\frac{1}{2y}\,dy\nonumber\\
&=\int_0^\infty\frac{1}{t}\Big(e^{-t/2}-e^{-tx}\Big)\Big(\frac{2}{t}-\frac{1}{e^t-1}\Big)\Big(1-e^{-t/2}\Big)\,dt
-\frac{1}{2}\int_0^\infty\frac{e^{-t/2}-e^{-tx}}{t}\,dt\nonumber\\
&=\int_0^\infty\frac{1}{t}\Big(e^{-t/2}-e^{-tx}\Big)\Big[\Big(\frac{2}{t}-\frac{1}{e^t-1}\Big)\Big(1-e^{-t/2}\Big)-\frac{1}{2}\Big]\,dt\nonumber\\
&=\int_0^\infty\frac{1}{t}\Big(e^{-t/2}-e^{-tx}\Big)\Big[\frac{1}{t}-\frac{1}{e^t-1}
+\frac{1}{t}-\frac{1}{2(1-e^{-t/2})}\Big]\Big(1-e^{-t/2}\Big)\,dt\nonumber\\
&=\int_0^\infty\frac{1}{t}\Big(e^{-t/2}-e^{-tx}\Big)\Big[\frac{1}{t}-\frac{1}{e^t-1}
+\frac{1}{2}\Big\{\frac{1}{t/2}-\frac{1}{e^{t/2}-1}\Big\}-\frac{1}{2}\Big]\Big(1-e^{-t/2}\Big)\,dt\nonumber
\end{align}
as in the statement of the Theorem.
\qed

\smallskip

\noindent Define
\[
h:[0,\infty)\rightarrow\mathbb{R};t\mapsto
\left\{
\begin{array}{cl}
\frac{1}{t}-\frac{1}{e^t-1} & \text{ for }t>0;\\
\frac{1}{2} & \text{ for }t=0;
\end{array}
\right.
\]
and
\begin{align}
\rho(t)&:=h(t)+\frac{1}{2}h(t/2)-\frac{1}{2}\label{definition_of_rho}
\end{align}
for $t\geq 0$.

\smallskip

\begin{lemma}\label{convexity_of_rho}
It holds that
\begin{itemize}
\item[(i)] $\cosh t\leq\Big(\frac{\sinh t}{t}\Big)^3$ for each $t>0$;
\item[(ii)] $h$ is continuous and convex on $[0,\infty)$;
\item[(iii)] $\rho$ is continuous and convex on $[0,\infty)$ with exactly one root;
\item[(iv)] $\rho>0$ on $[0,5/2]$.
\end{itemize}
\end{lemma}

\smallskip

\noindent{\em Proof.} 
{\em (i)} This reduces to the inequality
\[
\sum_{p=0}^\infty\frac{t^{2p}}{(2p)!}\leq\sum_{p=0}^\infty\Big\{\sum_{k+l+m=p}\frac{1}{(2k+1)!(2l+1)!(2m+1)!}\Big\}t^{2p}
\]
for each $t\geq 0$ after writing in power series. The inequality can be seen after comparing like coefficients.
{\em (ii)} Continuity at the origin is straightforward. A computation leads to the identity
\[
h^{\prime\prime}(t)=\frac{2}{t^3}-\frac{1}{4}\frac{\cosh(t/2)}{\sinh^3(t/2)}
\]
for each $t>0$. The right-hand side is nonnegative by {\em (i)}. {\em (iii)} The last item entails that the function $\rho$ is convex. Note that $\rho(0)=1/4$ and $\rho(t)\rightarrow -1/2$ as $t\rightarrow\infty$. It follows that $\rho$ has exactly one root. {\em (iv)} A calculation shows that $\rho(5/2)>0$.
\qed

\smallskip

\noindent We now record two integral identities for later use. First,
\begin{align}
\int_0^\infty\frac{1}{t}\Big(e^{-t/2}-e^{-tx}\Big)\Big(1-e^{-t/2}\Big)\,dt
&=
\int_0^\infty
\Big\{
\int_0^\infty e^{-t\tau}\,d\tau
\Big\}
\Big(e^{-t/2}-e^{-tx}\Big)\Big(1-e^{-t/2}\Big)\,dt\nonumber\\
&=
\int_0^\infty\Big\{
\int_0^\infty
 e^{-t\tau}
\Big(e^{-t/2}-e^{-tx}\Big)\Big(1-e^{-t/2}\Big)\,dt
\Big\}\,d\tau
\nonumber\\
&=
\int_0^\infty
\frac{1}{\tau+1/2}-\frac{1}{\tau+1}-\frac{1}{\tau+x}+\frac{1}{\tau+x+1/2}\,d\tau
\nonumber\\
&=-\log\Big(\frac{x+1/2}{2x}\Big)
=\log\Big(1+\frac{x-1/2}{x+1/2}\Big).\label{first_integral_identity}
\end{align}
Secondly,
\begin{align}
\int_0^\infty\Big(e^{-t/2}-e^{-tx}\Big)\Big(1-e^{-t/2}\Big)\,dt&
=\frac{x^2+(1/2)x-1/2}{x(x+1/2)}.\label{second_integral_identity}
\end{align}

\smallskip

\begin{lemma}\label{bound_for_Y}
For $x\geq 1$,
\[
\log\Big(1+\frac{x-1/2}{x+1/2}\Big)\geq\frac{2x^2+x-1}{5x(x+1/2)}>0.
\]
\end{lemma}

\smallskip

\noindent{\em Proof.} We bound the logarithm from below using the estimate
\[
\log(1+y)\geq\frac{2y}{2+y}
\]
for $y\geq 0$ which follows by the Hermite-Hadamard inequality. The resulting algebraic inequality holds if $x^3-x+1/8\geq 0$ and this cubic is positive on $[1,\infty)$.
\qed

\smallskip

\begin{proposition}\label{nonnegativity_of_integral_of_w}
For $x\geq 1$,
\begin{align}
\int_{1/2}^x w\,dy&>0.\nonumber
\end{align}
\end{proposition}

\smallskip

\noindent{\em Proof.}
Fix $x\geq 1$. We introduce the probability measure $\mathbb{P}$ on $(0,\infty)$ given by
\[
\mathbb{P}_x(d\tau)=
\frac{\frac{1}{\tau}\Big(e^{-\tau/2}-e^{-\tau x}\Big)\Big(1-e^{-\tau/2}\Big)\,d\tau}{\int_0^\infty\frac{1}{t}\Big(e^{-t/2}-e^{-tx}\Big)\Big(1-e^{-t/2}\Big)\,dt}
\]
and define the random variable $X:(0,\infty)\rightarrow\mathbb{R}$ via $t\mapsto t$. From (\ref{first_integral_identity}) and (\ref{second_integral_identity}) the expectation of $X$ under $\mathbb{P}$ is given by
\begin{align}
Y:=\mathbb{E}X
&=
\frac{\int_0^\infty\Big(e^{-t/2}-e^{-tx}\Big)\Big(1-e^{-t/2}\Big)\,dt}{\int_0^\infty\frac{1}{t}\Big(e^{-t/2}-e^{-tx}\Big)\Big(1-e^{-t/2}\Big)\,dt}
=\frac{\frac{x^2+(1/2)x-1/2}{x(x+1/2)}}{\log\Big(1+\frac{x-1/2}{x+1/2}\Big)}.\nonumber
\end{align}
By Lemma \ref{representation_of_integral_of_w}, Lemma \ref{convexity_of_rho} and Jensen's inequality,
\begin{align}
\frac{\int_{1/2}^x w\,dy}{\log\Big(1+\frac{x-1/2}{x+1/2}\Big)}&=\mathbb{E}\Big[\rho(X)\Big]
\geq\rho(\mathbb{E}[X])=\rho(Y).\nonumber
\end{align}
By Lemma \ref{convexity_of_rho}, $\rho>0$ on the interval $[0,5/2]$; on the other hand, $Y\leq 5/2$ by Lemma \ref{bound_for_Y}. This leads to the result.
\qed

\bigskip

\noindent We now turn to the second step in the proof of Theorem \ref{beta_function_inequality}. Recall from \cite{AbramowitzStegun1964} 6.4.10,
\[
\frac{\psi^{(n)}(x)}{n!}=(-1)^{n+1}\zeta(x,n+1)
\]
for each $x>0$. For $a>0$ and $s>1$,
\[
\zeta(a,s):=\sum_{p=0}^\infty\frac{1}{(a+p)^s}
\]
is the Hurwitz zeta function. For $|h|<1/2$,
\begin{align}
\psi(1/2+h)&=\sum_{k=0}^\infty\frac{\psi^{(k)}(1/2)}{k!}h^k
=\psi(1/2)+\sum_{k=1}^\infty(-1)^{k+1}\zeta(1/2,k+1)h^k;\label{Taylor_series_for_psi_at_half}\\
\psi(1+h)&=\sum_{k=0}^\infty\frac{\psi^{(k)}(1)}{k!}h^k=\psi(1)+\sum_{k=1}^\infty(-1)^{k+1}\zeta(1,k+1)h^k
.\label{Taylor_series_for_psi_at_one}
\end{align}

\smallskip

\begin{lemma}\label{sum_of_difference}
Suppose that $f$ is a positive strictly decreasing function of class $C^1$ defined on $[0,\infty)$ and suppose that both $f$ and $f^\prime$ are integrable. Then
\begin{align}
\sum_{j=1}^\infty\Big\{f(j)-f(j+1/2)\Big\}
&=\frac{1}{2}f(1/2)+\frac{1}{2}\sum_{j=0}^\infty(-1)^j\Big\{f(j/2+1)-f(j/2+1/2)\Big\}.\nonumber
\end{align}
\end{lemma}

\smallskip

\noindent{\em Proof.}
By the Euler summation formula (see \cite{Apostol1976} Theorem 3.1 for example),
\begin{align}
\sum_{j=1}^\infty\Big\{f(j)-f(j+1/2)\Big\}&=\int_0^\infty\Big\{f(t)-f(t+1/2)\Big\}\,dt+\int_0^\infty \{t\}\Big[f^\prime(t)-f^\prime(t+1/2)\Big]\,dt\nonumber\\
&=\int_0^{1/2}f(t)\,dt+\int_0^\infty \{t\}\Big[f^\prime(t)-f^\prime(t+1/2)\Big]\,dt\nonumber
\end{align}
where we use the notation $\{t\}:=t-[t]$ for $t\geq 0$ and $[\cdot]$ refers to the integer part. Let us observe that for any $t\geq 0$,
\[
\{t+1/2\}-\{t\}=
\left\{
\begin{array}{rcl}
1/2 & \text{ if } & \frac{j}{2}\leq t\leq\frac{j+1}{2} \text{ and }j\in\mathbb{Z}\text{ is even};\\
-1/2 & \text{ if } & \frac{j}{2}< t\leq\frac{j+1}{2} \text{ and }j\in\mathbb{Z}\text{ is odd}.\\
\end{array}
\right.
\]
\noindent Concentrating on the last member in the expression for the sum,
\begin{align}
\int_0^\infty\{t\}\Big[f^\prime(t)-f^\prime(t+1/2)\Big]\,dt
&=\int_0^\infty\{t\}f^\prime(t)-\{t+1/2\}f^\prime(t+1/2)+\Big[\{t+1/2\}-\{t\}\Big]f^\prime(t+1/2)\,dt
\nonumber\\
&=\int_0^{1/2} tf^\prime(t)\,dt+\sum_{j=0}^\infty\int_{\frac{j}{2}}^{\frac{j+1}{2}}\Big[\{t+1/2\}-\{t\}\Big]f^\prime(t+1/2)\,dt\nonumber\\
&=\int_0^{1/2} tf^\prime(t)\,dt+\frac{1}{2}\sum_{j=0}^\infty(-1)^j\int_{\frac{j}{2}}^{\frac{j+1}{2}}f^\prime(t+1/2)\,dt\nonumber\\
&=\int_0^{1/2} tf^\prime(t)\,dt+\frac{1}{2}\sum_{j=0}^\infty(-1)^j
\Big\{f(j/2+1)-f(j/2+1/2)\Big\}.\nonumber
\end{align}
Lastly let us note that
\[
\int_0^{1/2}f(t)+tf^\prime(t)\,dt=\int_0^{1/2}(tf)^\prime\,dt=\frac{1}{2}f(1/2).
\]
Combining this with the earlier identities gives the result. 
\qed

\smallskip

\begin{corollary}\label{formula_for_difference_between_zeta_functions}
For $a>0$ and $s>1$,
\begin{align}
\zeta(a,s)-\zeta(a+1/2,s)&=
\frac{1}{a^s}-\frac{1}{2}\frac{1}{(a+1/2)^s}
+\frac{1}{2}\sum_{j=0}^\infty(-1)^j
\Big\{\frac{1}{(a+j/2+1)^s}-\frac{1}{(a+j/2+1/2)^s}\Big\}.\nonumber
\end{align}
\end{corollary}

\smallskip

\noindent{\em Proof.}
First note that
\begin{align}
\zeta(a,s)-\zeta(a+1/2,s)&=
\frac{1}{a^s}-\frac{1}{(a+1/2)^s}
+\sum_{j=1}^\infty\Big[\frac{1}{(a+j)^s}-\frac{1}{(a+1/2+j)^s}\Big].\nonumber
\end{align}
Now apply Lemma \ref{sum_of_difference} with $f(t):=(a+t)^{-s}$ to the summation.
\qed

\smallskip

\begin{lemma}
For $|h|<1/2$,
\begin{align}
w(1/2+h)&=2\log\Big(\frac{1+h}{1+2h}\Big)+\frac{h(6h+5)}{2(1+h)(1+2h)}\nonumber\\
&-2h\sum_{j=0}^\infty(-1)^j\frac{2j+5+2h}{(j+2)(j+3)(j+2+2h)(j+3+2h)}\label{representation_of_w}
\end{align}
\end{lemma}

\smallskip

\noindent{\em Proof.} Making use of the identities (\ref{Taylor_series_for_psi_at_half}) and (\ref{Taylor_series_for_psi_at_one}), an application of Corollary \ref{formula_for_difference_between_zeta_functions} yields the identity
\begin{align}
\psi(1/2+h)-\psi(1+h)&=\psi(1/2)-\psi(1)+\frac{7h+6h^2}{(1+2h)(2+2h)}\nonumber\\
&-2h\sum_{j=0}^\infty(-1)^j\frac{2j+5+2h}{(j+2)(j+3)(j+2+2h)(j+3+2h)}\nonumber
\end{align}
after some careful calculation where $|h|<1/2$. We remark that 
\begin{align}
w(1/2)&=\psi(1/2)-\psi(1)+2\log 2=0\label{w_at_half}
\end{align}
after substituting the $\psi$-values in \cite{AbramowitzStegun1964} 6.3.2 and 6.3.3. Upon writing $x=1/2+h$ for $|h|<1/2$ we compute
\begin{align}
w(x)&=\phi(x)-\phi(x+1/2)+\log\Big(1+\frac{1}{2x}\Big)-\frac{1}{2x}\nonumber\\
&=2\log\Big(1+\frac{1}{2x}\Big)+\frac{1}{2x}-\frac{1}{x+1/2}+\psi(x)-\psi(x+1/2)\nonumber\\
&=2\log\Big(\frac{2+2h}{1+2h}\Big)-\frac{h}{(1+h)(1+2h)}+\psi(1/2+h)-\psi(1+h)\nonumber\\
&=2\log\Big(\frac{2+2h}{1+2h}\Big)-\frac{h}{(1+h)(1+2h)}+\psi(1/2)-\psi(1)+\frac{7h+6h^2}{(1+2h)(2+2h)}\nonumber\\
&-2h\sum_{j=0}^\infty(-1)^j\frac{2j+5+2h}{(j+2)(j+3)(j+2+2h)(j+3+2h)}\nonumber\\
&=2\log\Big(\frac{1+h}{1+2h}\Big)+\frac{h(6h+5)}{2(1+h)(1+2h)}
-2h\sum_{j=0}^\infty(-1)^j\frac{2j+5+2h}{(j+2)(j+3)(j+2+2h)(j+3+2h)}\nonumber
\end{align}
where we made use of (\ref{w_at_half}) to obtain the last line.
\qed

\smallskip

\begin{proposition}\label{positivity_of_w_on_half_one}
The function $w$ is positive on the interval $(1/2,1)$.
\end{proposition}

\smallskip

\noindent{\em Proof.}
For fixed $h$ in $(0,1/2)$ the sequence $(a_j)_{j=0}^\infty$ with
\[
a_j:=\frac{2j+5+2h}{(j+2)(j+3)(j+2+2h)(j+3+2h)}
\]
is a strictly decreasing null sequence noting that $2j+5+2h=2(j+2+2h)+1-2h$. The alternating series in (\ref{representation_of_w}) is bounded above by its first term. We thus obtain
\begin{align}
w(1/2+h)&\geq2\log\Big(\frac{1+h}{1+2h}\Big)+\frac{h(6h+5)}{2(1+h)(1+2h)}-\frac{h(5+2h)}{6(1+h)(3+2h)}.\nonumber
\end{align}
We bound the logarithm from below using the estimate
\[
\log(1+y)\geq\frac{y(2+y)}{2(1+y)}\text{ for }-1<y\leq 0.
\]
This follows from the Hermite-Hadamard inequality. Continuing, 
\begin{align}
w(1/2+h)&\geq
-\frac{h(2+3h)}{(1+h)(1+2h)}+\frac{h(6h+5)}{2(1+h)(1+2h)}-\frac{h(5+2h)}{6(1+h)(3+2h)}\nonumber\\
&=\frac{2h(2+h)(1/2-h)}{3(1+h)(1+2h)(3+2h)}>0\nonumber
\end{align}
for $0<h<1/2$.
\qed
\section{An isoperimetric competitor}\label{isoperimetric_competitor}

\smallskip

\noindent There is a natural candidate for isoperimetric minimiser in problem (\ref{isoperimetric_problem}) which is not a centred ball. This candidate set is spherical cap symmetric and its boundary contains the origin but this fails to be of class $C^1$ there. Our task here is to show that this candidate is not isoperimetric. 

\smallskip

\noindent Let $(\alpha,\beta)\in\mathcal{P}$ and $b>0$. As in Lemma \ref{first_order_ode_with_a_zero} the constant generalised curvature equation $(\tau^{\beta+1}u)^\prime+\lambda\tau^{\alpha+1}=0$ on $(0,b)$ with boundary conditions $u(0)=0$ and $u(b)=1$ has the unique solution $u(\tau)=(\tau/b)^\gamma$ for $\tau\in[0,b]$ where we write $\gamma:=\alpha-\beta+1$. Put
\begin{align}
\theta(\tau)&:=\int_\tau^b\frac{u}{\sqrt{1-u^2}}\,d\mu=
\int_{\tau/b}^1\frac{t^{\gamma-1}}{\sqrt{1-t^{2\gamma}}}\,dt\label{theta_for_competitor}
\end{align}
for $\tau\in[0,b]$. This angular coordinate may be expressed in more explicit form. Let us explain this briefly - though we do not make further use of this fact. On substituting $s=t^\gamma$ we obtain
\[
\theta(\tau)=\frac{1}{2\gamma}\int_0^{1-(\tau/b)^{2\gamma}}\frac{ds}{\sqrt{s(1-s)}}.
\]
At this point let us note the integral identity
\begin{align}
\int_x^{1/2}\frac{ds}{\sqrt{s(1-s)}}&=\frac{\pi}{2}-\mathrm{arcsin}\sqrt{4x(1-x)}\label{integral_identity_involving_arcsin}
\end{align}
which holds for $0\leq x\leq 1/2$. This can be seen using the substitution $\sqrt{s(1-s)}=z$ as in \cite{Gradshteynetal1965} 2.292. Making use of this last identity (\ref{integral_identity_involving_arcsin}) we alight on a more explicit formula
\begin{align}
\theta(\tau&)=
\left\{
\begin{array}{rl}
\frac{1}{2\gamma}\Big(\pi-\mathrm{arcsin}\Big[2(\tau/b)^\gamma\sqrt{1-(\tau/b)^{2\gamma}}\Big]\Big) & \text{ if }0\leq(\tau/b)^{2\gamma}\leq 1/2;\\
\frac{1}{2\gamma}\mathrm{arcsin}\Big[2(\tau/b)^\gamma\sqrt{1-(\tau/b)^{2\gamma}}\Big] & \text{ if }1/2<(\tau/b)^{2\gamma}\leq 1.\label{explicit_formula_for_theta}
\end{array}
\right.
\end{align}
for $\theta(\tau)$. It can be seen that $\theta(0)=\frac{\pi}{2\gamma}$ (though this also follows from the first equality in (\ref{theta_for_competitor}) using the substitution $t^\gamma=\sin\theta$). 

\begin{figure}[h]
\centering
\includegraphics[width=80mm]{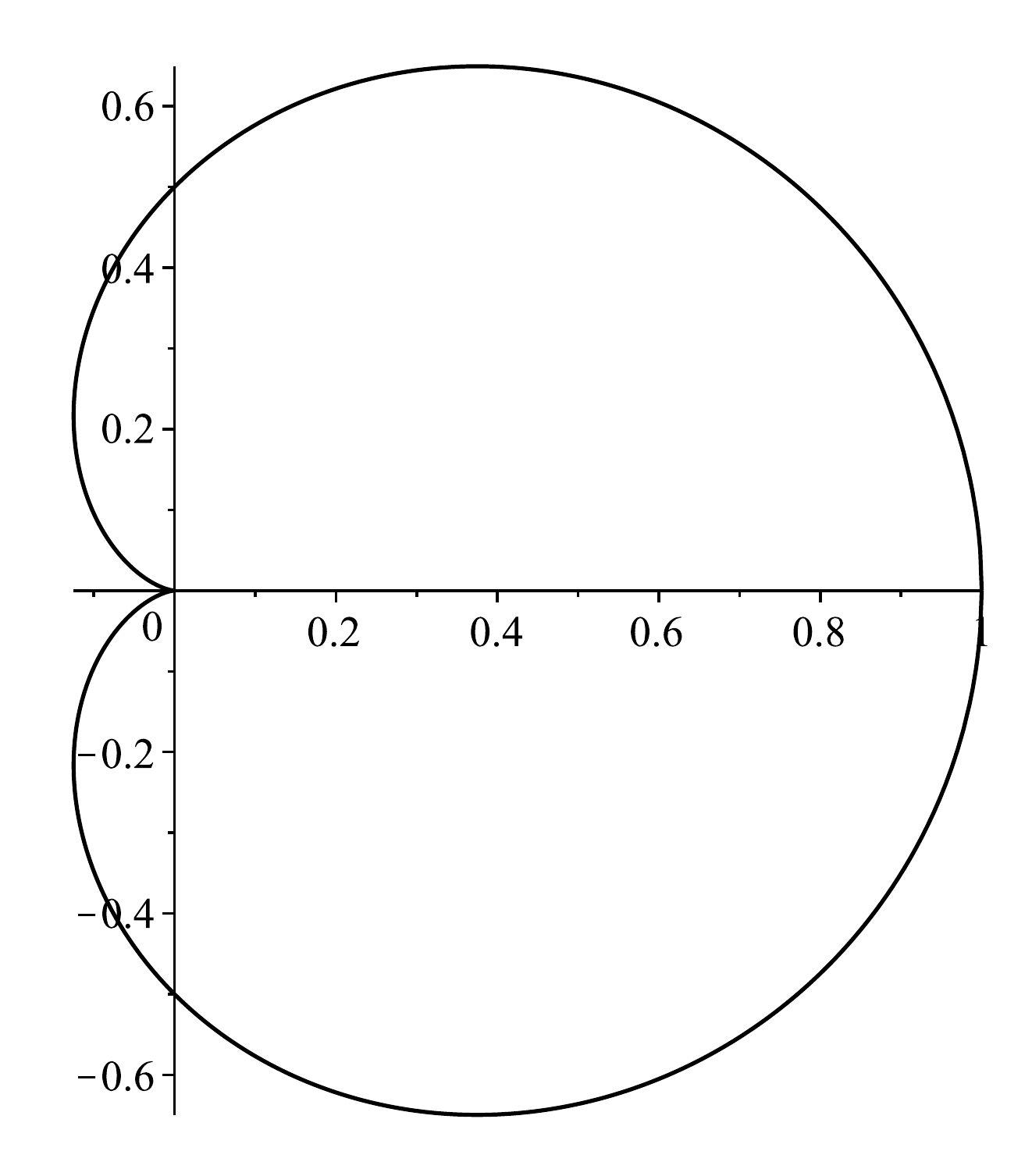}
\caption{Plot of the set $E$ in (\ref{definition_of_the_set_E}) with $\gamma=1/2$ and $b=1$. Its boundary has a cusp at $0$. In case $1/2<\gamma<1$ there is a Lipschitz singularity at $0$. If $\gamma=1$ and $b=1$ then $E$ is a ball with centre $(1/2,0)$ and tangent to the origin.}
\label{fig:graph_of_the_set_E}
\end{figure}

\smallskip

\noindent In the remainder of this subsection we shall assume that $1/2\leq\gamma<1$ so that $\theta(0)\in(\pi/2,\pi]$. Define the open set
\begin{align}
E&:=\{(\tau\cos\phi,\tau\sin\phi):0<\tau<b\text{ and }|\phi|<\theta(\tau)\}\subset\mathbb{R}^2_0\label{definition_of_the_set_E}
\end{align}
with $\theta$ defined as in (\ref{theta_for_competitor}). For $1/2\leq\gamma<1$ the boundary of $E$ fails to be $C^1$ at the origin as $\theta(0)\in(\pi/2,\pi]$. See Figure \ref{fig:graph_of_the_set_E}. Incidentally, the regularity result for an isoperimetric minimiser in \cite{Morgan2003} 3.10 fails to apply in this setting because the metric $\zeta=r^{\gamma-1}$ has a singularity at the origin and the density $\psi=r^{2\beta-\alpha}$ vanishes at the origin.

\smallskip

\noindent Let us now compute the weighted perimeter and volume of the set $E$. From Lemma \ref{derivative_of_theta} and Lemma \ref{first_order_ode_with_a_zero},
\begin{align}
P_\beta(E)&=2\int_0^b\sqrt{1+t^2(\theta^\prime)^2}t^\beta\,dt
=2b^{\beta+1}\int_0^1\frac{t^\beta}{\sqrt{1-t^{2\gamma}}}\,dt
=\frac{b^{\beta+1}}{\gamma}B(x,1/2)\nonumber
\end{align}
where the beta function is defined as in (\ref{beta_function}) and $x:=\frac{\beta+1}{2\gamma}$. Turning to the weighted volume,
\begin{align}
V_\alpha(E)&=2\int_0^b\theta t^{\alpha+1}\,dt=-2\int_0^b\theta^\prime\frac{t^{\alpha+2}}{\alpha+2}\,dt
=2\int_0^b\frac{1}{t}\frac{(t/b)^\gamma}{\sqrt{1-(t/b)^{2\gamma}}}\frac{t^{\alpha+2}}{\alpha+2}\,dt\nonumber\\
&=\frac{b^{\alpha+2}}{(\alpha+2)\gamma}B(x+1,1/2)
=\frac{b^{\alpha+2}}{(2x+1)\gamma^2}B(x+1,1/2)
\nonumber
\end{align}
with $x$ as before using the relation $\alpha+2=(2x+1)\gamma$. By the recurrence relation for the beta function $B(x+1,1/2)=\frac{x}{x+1/2}B(x,1/2)$ we may write
\begin{align}
V_\alpha(E)&=\frac{b^{\alpha+2}}{\gamma^2}\frac{2x}{(2x+1)^2}B(x,1/2).
\nonumber
\end{align}
Let us notice here that the ratio
\begin{align}
\frac{P_\beta(E)^{\alpha+2}}{V_\alpha(E)^{\beta+1}}&
=\gamma^{2\beta-\alpha}\Big[\frac{(2x+1)^2}{2x}\Big]^{\beta+1}B(x,1/2)^\gamma
\label{isoperimetric_ratio_for_E}
\end{align}
does not depend on the choice of $b$.

\smallskip

\noindent Now let $B$ stand for the open centred ball in $\mathbb{R}^2_0$ with radius $r>0$. Then the isoperimetric ratio
\begin{align}
\frac{P_\beta(B)^{\alpha+2}}{V_\alpha(B)^{\beta+1}}&=(2\pi)^\gamma(\alpha+2)^{\beta+1}
=(2\pi)^\gamma\gamma^{\beta+1}(2x+1)^{\beta+1}
\label{isoperimetric_ratio_for_ball}
\end{align}
is likewise scale invariant. The main claim of this subsection is that
\begin{align}
\frac{P_\beta(E)^{\alpha+2}}{V_\alpha(E)^{\beta+1}}
&>\frac{P_\beta(B)^{\alpha+2}}{V_\alpha(B)^{\beta+1}}.\label{relation_between_isoperimetric_ratios}
\end{align}
From (\ref{isoperimetric_ratio_for_E}) and (\ref{isoperimetric_ratio_for_ball}) this holds if and only if
\begin{align}
\gamma^{-1}\Big(1+\frac{1}{2x}\Big)^{2x}B(x,1/2)&>2\pi\label{desired_inequality}
\end{align}

\smallskip

\noindent Comparing the definition of $\zeta$ in (\ref{zeta}) with the expression for $x$ above we see that $\zeta=2x$. In other words, the expression $x$ which depends on $(\alpha,\beta)$ in $\mathscr{P}$ is constant on the graph of the line $\ell_{2x}$. From Lemma \ref{lines_and_set_P}, $\gamma$ attains its maximum where the graph of $\ell_{2x}$ meets $\mathcal{P}^+$.  Denote this point by $(\overline{\alpha},\overline{\beta})$ as before. The pair $(\overline{\alpha},\overline{\beta})$ satisfies the simultaneous equations
\[
\left\{
\begin{array}{rcl}
\alpha(\beta+1)&=&\beta^2;\\
\alpha+2&=&[1+\frac{1}{2x}](\beta+1).\\
\end{array}
\right.
\]
Eliminating $\alpha$ we see that $\overline{\beta}+1=\sqrt{2x}$ using positivity of $\beta+1$ and $\overline{\gamma}=\frac{1}{\sqrt{2x}}$ appealing to the identity $\beta+1=(2x)\gamma$.

\smallskip

\noindent Now take $(\alpha,\beta)\in\mathcal{P}\setminus\mathcal{P}^-$ with $1/2\leq\gamma<1$. Then
\begin{align}
\gamma^{-1}\Big(1+\frac{1}{2x}\Big)^{2x}B(x,1/2)
&\geq\overline{\gamma}^{-1}\Big(1+\frac{1}{2x}\Big)^{2x}B(x,1/2)
=\sqrt{2x}\Big(1+\frac{1}{2x}\Big)^{2x}B(x,1/2)=W(x)\nonumber
\end{align}
in the notation of Theorem \ref{beta_function_inequality}. By Lemma \ref{lines_and_set_P}, $\zeta=2x\geq 1$ so that $x\geq 1/2$. The case $x=1/2$ corresponds to $\zeta=1$ and $\alpha=2\beta$ or in other words $(\alpha,\beta)\in\mathcal{P}^-$ which we are here excluding so $x>1/2$. By Theorem \ref{beta_function_inequality} the right-hand side above is positive. This establishes the claim in (\ref{relation_between_isoperimetric_ratios}).

\section{A necessary condition for spherical cap symmetry}

\smallskip

\noindent In this short Section we record two results which will be used in the proof of Theorem \ref{Omega_is_empty}; namely, Lemma \ref{curvature_and_spherial_cap_symmetry} and Proposition \ref{positivity_of_lambda}.

\smallskip

\begin{lemma}\label{lemma_on_lb_for_curvature}
Let $E$ be a bounded open set with $C^2$ boundary $M$ in $\mathbb{R}^2_0$ and assume that $E=E^{sc}$. Put $R:=\sup\{|x|:x\in M\}>0$ and $x:=(R,0)$. Then the curvature of $M$ at $x$ satisfies the estimate $k(x)\geq 1/R$.
\end{lemma}

\smallskip

\noindent{\em Proof.} Let $\gamma_1:I\rightarrow M$ be a $C^2$ local parametrisation of $\gamma$ with $\gamma_1(0)=x$ as in Section \ref{Section_on_regularity_etc}. Put $r_1:=|\gamma_1|$. By \cite{McGillivray2018} (2.9) (for example), $\dot{r}_1=\cos\sigma_1$ and differentiating once more $\ddot{r}_1=-\sin\sigma_1\Big(k_1-\frac{\sin\sigma_1}{r_1}\Big)$ on $I$ using the decomposition $\sigma_1=\alpha_1-\theta_1$ and \cite{McGillivray2018} (2.10). In particular, $\ddot{r}_1(0)=-k(x)+1/R$. As $r_1$ has a global maximum on $I$ at $s=0$ it follows by the second derivative test that $-k(x)+1/R\leq 0$.
\qed

\smallskip

\noindent Let
\[
H:=\{x=(x_1,x_2)\in\mathbb{R}^2:x_2>0\}
\]
stand for the open upper half-plane in $\mathbb{R}^2$.

\smallskip

\begin{lemma}\label{curvature_and_spherial_cap_symmetry}
Let $E$ be an open set with $C^2$ boundary $M:=\partial E\cap\mathbb{R}^2_0$ in $\mathbb{R}^2_0$. Assume that $E=E^{sc}$. Suppose
\begin{itemize}
\item[(i)] $x\in M\cap H$;
\item[(ii)] $\sin(\sigma(x))=-1$.
\end{itemize}
Then $k(x)+1/a=0$ where $a:=|x|$.
\end{lemma}

\smallskip

\noindent{\em Proof.} Let $\gamma_1:I\rightarrow M$ be a $C^2$ parametrisation of $M$ in a neighbourhood of $x$ with $\gamma_1(0)=x$ as above.  We use the same notation as before. At $s=0$, $\dot{r}_1(0)=0$ and $\ddot{r}_1(0)=k_1(0)+1/a$. By Taylor's Theorem,
\[
r_1(s)=r_1(0)+\dot{r}_1(0)s+\frac{\ddot{r}_1(c)}{2}s^2
=a-(1/2)s^2\Big(\sin\sigma_1\Big[k_1-\frac{\sin\sigma_1}{r_1}\Big]\Big)(c)
\]
for $s\in I$ and some $c$ between $0$ and $s$. At $s=0$ the expression in round brackets takes the value $k(x)+1/a$. If $k(x)+1/a>0$ then $r_1>a$ in a neighbourhood of $s=0$. This contradicts spherical cap symmetry. If the expression $k(x)+1/a$ happens to be strictly negative on the other hand then the above expansion implies that $r_1<a$ near $0$. This contradicts the fact that $r_1\geq a$ on $I\cap\{s<0\}$. This is because $\cos\sigma_1\leq 0$ on $I$ by \cite{McGillivray2018} Lemma 5.4 as well as the above expression for the derivative of $r_1$. See Figure \ref{fig:cuvature_at_sin_sigma}. 
\qed

\begin{figure}[h]
\centering
\begin{tikzpicture}[>=stealth][scale=0.33]

\draw[->] (-1,0) -- (4,0);
\draw[->] (0,-0.5) -- (0,4);
\draw[dotted] (2,0) arc (0 : 105 : 2);

\draw[thick] (3,3.5) to[out=165,in=145,looseness=1.2] (1.414,1.414);
\draw[thick] (1,0.5) to[out=30,in=-45,looseness=1.2] (1.414,1.414);
\draw[->,dashed] (1.414,1.414) to (2.414,2.414) node[right] {$n$};

\coordinate (A) at   (2,0);
\coordinate (B) at   (1.414,1.414);
\coordinate (C) at   (2.414,2.414);

\filldraw (A) circle (2pt) node[below] {$a$};
\filldraw (B) circle (2pt) node[left] {$x$};

\end{tikzpicture}
\label{fig:cuvature_at_sin_sigma}
\caption{At the point $x\in\partial E$ the sine of $\sigma$ is $-1$. The inward unit normal to $E$ at $x$ is the vector $n$. The curvature at $x$ is $-1/a$. } 
\end{figure}
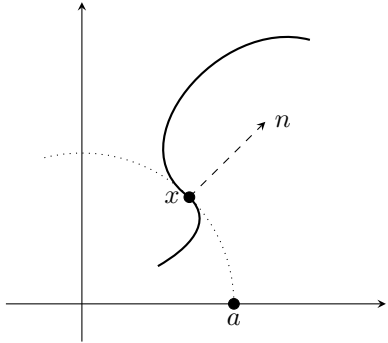

\smallskip

\begin{proposition}\label{positivity_of_lambda}
Let $(\alpha,\beta)\in\mathcal{Q}$. Given $v>0$ let $E$ be a minimiser of (\ref{isoperimetric_problem}). Assume that $E$ is a bounded open set with analytic boundary $M$ in $\mathbb{R}^2_0$ and suppose that $E=E^{sc}$. Let $\lambda$ be as in Theorem \ref{constant_weighted_mean_curvature}. Then $\lambda<0$.
\end{proposition}

\smallskip

\noindent{\em Proof.}
Let $R$ be as in Lemma \ref{lemma_on_lb_for_curvature} and put $x:=(R,0)$ as before. By Lemma \ref{lemma_on_lb_for_curvature}, 
$k(x)\geq 1/R$. By Theorem \ref{constant_weighted_mean_curvature}, $k(x)=-\frac{\beta}{R}-\lambda R^{\alpha-\beta}\geq 1/R$; that is, $\lambda\leq-(\beta+1)R^{-\alpha+\beta-1}<0$.
\qed
\section{Proof of main result}

\smallskip

\noindent This Section contains the proof of Theorem \ref{main_theorem}. Recall the set $\Omega$ defined in (\ref{definition_of_Omega}). This set comprises all those radii where the boundary curve crosses the foliation of centred circles transversally. The first step in the proof is to show that this set is empty.

\noindent 

\begin{theorem}\label{Omega_is_empty}
Let $(\alpha,\beta)\in\mathcal{P}\setminus\mathcal{P}^-$. Given $v>0$ let $E$ be a bounded minimiser of (\ref{isoperimetric_problem}). Assume that $E$ is open in $\mathbb{R}^2_0$, $M:=\partial E\cap\mathbb{R}^2_0$ is an analytic hypersurface in $\mathbb{R}^2_0$ and $E=E^{\mathrm{sc}}$. Then $\Omega=\emptyset$.
\end{theorem}

\smallskip

\noindent{\em Proof.}
Suppose that $\Omega\neq\emptyset$. The set $\Omega$ is open in $(0,\infty)$ by \cite{McGillivray2018} Lemma 5.6. So we may write $\Omega$ as a countable union of disjoint open intervals in $(0,\infty)$. By a suitable choice of one of these intervals we may assume that the interval $(a,b)$ is a connected component of $\Omega$ for some $0\leq a<b<\infty$. Let us assume for the time being that $a>0$. As $M$ is analytic, $[a,b]\subset\pi(M)$ and $\cos\sigma$ vanishes on $M_a\cup M_b$. Let $\lambda\in\mathbb{R}$ be as in Theorem \ref{constant_weighted_mean_curvature}. 

\smallskip

\noindent Let $u:\Omega\rightarrow[-1,1]$ be as in (\ref{definition_of_y}). Then $u$ has a continuous extension to $[a,b]$ and $u=\pm 1$ at $\tau=a,b$. This may be seen as follows. For $\tau\in(a,b)$ the set $M_\tau\cap H$ consists of a singleton by \cite{McGillivray2018} Lemma 5.4. The limit $x:=\lim_{\tau\downarrow a}M_\tau\cap H\in\mathbb{S}_a\cap\overline{H}$ exists as the boundary $M$ is analytic in $\mathbb{R}^2_0$. In other words, $\gamma(\tau)\rightarrow x$ as $\tau\downarrow a$. Observe that $\cos\sigma(x)=0$ by openness of $\Omega$ and hence $\sin\sigma(x)=\pm 1$. The tangent vector varies continuously along $M$. So $\sigma(\gamma(\tau))\rightarrow\sigma(x)$ and hence $\sin(\sigma(\gamma(\tau)))\rightarrow\sin(\sigma(x))$ as $\tau\downarrow a$. This shows that $u\rightarrow\pm 1$ as $\tau\downarrow a$. The argument at $b$ is similar.

\smallskip

\noindent Put $\eta_1:=u(a)$ and $\eta_2:=u(b)$. Let us consider the case $\eta=(\eta_1,\eta_2)=(1,1)$. Note that $u<1$ on $(a,b)$ for otherwise $\cos(\sigma\circ\gamma)$ vanishes at some point in $(a,b)$ contradicting the definition of $\Omega$. By Theorem \ref{ode_for_y} the pair $(u,\lambda)$ satisfies (\ref{ode_for_u_with_plus_minus_1_at_endpoints}) with $\eta=(1,1)$. By Lemma \ref{first_order_ode}, $u>0$ on $[a,b]$. Put $w:=1/u$. Then $(w,-\lambda)$ satisfies (\ref{Riccati_equation}) and $w>1$ on $(a,b)$. By Lemma \ref{derivative_of_theta},
\[
\theta_2(b)-\theta_2(a)=\int_a^b\theta_2^\prime\,d\tau
=-\int_a^b\frac{u}{\sqrt{1-u^2}}\,\frac{d\tau}{\tau}
=-\int_a^b\frac{1}{\sqrt{w^2-1}}\,\frac{d\tau}{\tau}.
\]
By Corollary \ref{integral_of_w}, $|\theta_2(b)-\theta_2(a)|>\pi$. But this contradicts the definition of $\theta_2$ in (\ref{definition_of_theta}) as $\theta_2$ takes values in $(0,\pi)$ on $(a,b)$ according to \cite{McGillivray2018} Lemma 5.4. If $\eta=(-1,-1)$ then $\lambda>0$ by Lemma \ref{first_order_ode}; this contradicts Proposition \ref{positivity_of_lambda}.

\smallskip

\noindent Now let us consider the case $\eta=(-1,1)$. As before the limit $x:=\lim_{\tau\downarrow a}M_\tau\cap H\in\mathbb{S}_a\cap\overline{H}$ exists as $M$ is $C^1$ in $\mathbb{R}^2_0$ and $u(a)=\sin\sigma(x)=-1$. Using the same formula as above, $\theta_2(b)-\theta_2(a)<0$ by Corollary \ref{integral_of_solution_of_first_order_ode}. This means that $\theta_2(a)\in(0,\pi]$. We may assume that $\theta_2(a)\in(0,\pi)$. For if $\theta_2(a)=\pi$ then $x=(-a,0)$ and the tangent vector to $M$ at $x$ is pointing into the upper half-plane. This implies that $\theta_2$ exceeds $\pi$ near $a$ and contradicts the definition of $\theta_2$ in (\ref{definition_of_theta}). We proceed on the basis of this assumption. By Theorem \ref{constant_weighted_mean_curvature} the curvature at $a$ is given by 
\[
k=\frac{\beta}{a}-\lambda a^{\alpha-\beta}=\frac{1}{a}\Big\{\beta-\lambda a^{\alpha-\beta+1}\Big\}
=\frac{1}{a}\Big\{\beta+\hat{m}a^{\alpha-\beta+1}\Big\}
\]
and
\[
k+1/a=\frac{1}{a}\Big\{\beta+1+\hat{m}a^{\alpha-\beta+1}\Big\}>0.
\]
This contradicts Lemma \ref{curvature_and_spherial_cap_symmetry}. If $\eta=(1,-1)$ then $\lambda>0$ by Lemma \ref{first_order_ode}. This contradicts Proposition \ref{positivity_of_lambda}.
\smallskip

\noindent Suppose finally that $a=0$. Then $u$ satisfies the ordinary differential equation (\ref{ode_for_u}) on $(0,b)$ for some $\lambda\in\mathbb{R}$. This may be rewritten $(\tau^{\beta+1}u)^\prime+\lambda\tau^{\alpha+1}=0$ or integrating
\[
u+\frac{\lambda}{\alpha+2}\tau^{\alpha-\beta+1}+c\tau^{-(\beta+1)}=0\text{ on }(0,b)
\]
for some constant $c\in\mathbb{R}$. The first two terms are bounded on $(0,b)$ while the third is unbounded if $c\neq 0$. This entails that $c=0$ and hence that $u(0+)=0$. 

\smallskip

\noindent Let us assume that $u(b)=1$. From the earlier part of the proof the set $\Omega$ coincides with $(0,b)$. Notice that $\theta_2$ takes values in $(0,\pi)$ on $(0,b)$. As in Section \ref{isoperimetric_competitor}, $\theta(0+)=\frac{\pi}{2\gamma}$. If $\gamma\in(0,1/2)$ this last exceeds $\pi$ strictly. This contradicts the definition of $\theta_2$ because $\theta$ coincides with $\theta_2$. Suppose then that $\gamma\in[1/2,1)$. We claim that $\theta_2(b)=0$. Suppose for a contradiction that $\theta_2(b)\in(0,\pi)$. Then the boundary $M$ of $E$ in $\mathbb{R}^2_0$ includes a centred circular arc with radius $b$. From the parametric formulae \cite{McGillivray2018} (2.9) and (2.10) (for example),
\[
\frac{d\tau}{d\theta_2}=\frac{\dot{r}_1}{\dot{\theta}_1}=-\tau\frac{\sqrt{1-u^2}}{u}
\]
on the image of $(0,b)$ under $\theta_2$. By scale-invariance we may assume that $b=1$. We derive that
\[
\frac{d\tau}{d\theta}=
\left\{
\begin{array}{lcl}
0 & \text{ if } & 0<\theta<\theta_2(b);\\
-\tau^{1-\gamma}\sqrt{1-\tau^{2\gamma}} & \text{ if } &\theta_2(b)<\tau<\theta_2(0);\\
\end{array}
\right.
\]
using Lemma \ref{first_order_ode_with_a_zero}. It is then apparent that the boundary $M$ fails to be analytic at the point with polar coordinates $(b,\theta_2(b))$ which contradicts a hypothesis of the Theorem. This shows that $\theta_2(b)=0$. Let us continue with the argument. Let $B$ stand for the centred ball which has the same weighted area as $E$. According to the claim in (\ref{relation_between_isoperimetric_ratios}) the set $B$ has smaller weighted perimeter compared to the set $E$. This contradicts the fact that $E$ is an isoperimetric minimiser.

\smallskip

\noindent Now consider the case in which $u(b)=-1$. We may choose $x$ in $M$ in the closed upper half-plane; and in fact $x$ may be chosen in the open upper half-pane $H$ because of the condition $u(b)=-1$. By Theorem \ref{constant_weighted_mean_curvature}, 
\[
k-\frac{\beta}{b}+\lambda b^{\alpha-\beta}=0
\]
at $x$. By Lemma \ref{curvature_and_spherial_cap_symmetry}, $k(x)=-1/b$. Substituting into the expression above we obtain that $\lambda=(\beta+1)b^{-\gamma}>0$. This contradicts Proposition \ref{positivity_of_lambda}. 
\qed

\smallskip

\begin{lemma}\label{boundary_of_E_as_union_of_circles}
Let $(\alpha,\beta)\in\mathcal{P}\setminus\mathcal{P}^-$. Given $v>0$ let $E$ be a bounded minimiser of (\ref{isoperimetric_problem}). Assume that $E$ is open in $\mathbb{R}^2_0$, $M:=\partial E$ is an analytic hypersurface in $\mathbb{R}^2_0$ and $E=E^{sc}$. Then $M$ consists of an at most countable union of disjoint centred circles which accumulate at the origin if at all. 
\end{lemma}

\smallskip

\noindent{\em Proof.} 
The set $E$ is open and bounded so $M\neq\emptyset$ and hence $\emptyset\neq\pi(M)=\pi(M)\setminus\Omega$ by Theorem \ref{Omega_is_empty}. By definition of $\Omega$, $\cos\sigma=0$ on $M$. Let $0<\tau\in\pi(M)$. We claim that $M_\tau=\mathbb{S}_\tau$. Suppose for a contradiction that $M_\tau\neq\mathbb{S}_\tau$. By \cite{McGillivray2018} Lemma 5.2, $M_\tau$ is the union of two closed  spherical arcs in $\mathbb{S}_\tau$. Let $x$ be a point on the boundary of one of these spherical arcs relative to $\mathbb{S}_\tau$. There exists an analytic parametrisation $\gamma_1:I\rightarrow M$ of $M$ in a neighbourhood of $x$ with $\gamma_1(0)=x$ with conventions as mentioned in Section \ref{Section_on_regularity_etc}. By \cite{McGillivray2018} (2.9), $\dot{r}_1=0$ on $I$ as $\cos\sigma_1=0$ on $I$ because $\cos\sigma=0$ on $M$; that is, $r_1$ is constant on $I$. This means that $\gamma_1(I)\subset\mathbb{S}_\tau$. As the function $\sin\sigma_1$ is continuous on $I$ it takes the value $\pm 1$ there. By \cite{McGillivray2018} (2.10), $r_1\dot{\theta}_1=\sin\sigma_1=\pm 1$ on $I$. This means that $\theta_1$ is either strictly decreasing or strictly increasing on $I$. This entails that the point $x$ is not a boundary point of $M_\tau$ in $\mathbb{S}_\tau$ and this proves the claim.
\qed

\smallskip

\noindent{\em Proof of Theorem \ref{main_theorem}.}
The result follows from \cite{MorganRitore2002} Theorem 1.1 in case $(\alpha,\beta)\in\mathcal{P}^-\setminus\{(0,0)\}$. Let us assume then that $(\alpha,\beta)\in\mathcal{P}\setminus\mathcal{P}^-$.

\smallskip

\noindent{\em Existence.} 
By Corollary \ref{existence_and_boundedness_for_indices_in_curly_Q} the variational problem (\ref{isoperimetric_problem}) admits a bounded minimiser $E$. By Theorem \ref{C1_property_of_reduced_boundary} we may assume that $E$ is open in $\mathbb{R}^2_0$, $M:=\partial E$ is an analytic hypersurface in $\mathbb{R}^2_0$ and $E=E^{sc}$. By Lemma \ref{boundary_of_E_as_union_of_circles} $M$ consists of an at most countable union of disjoint centred circles which accumulate at the origin if at all. As such we may write
\[
E=\bigcup_{h=0}^NA((a_{2h+1},a_{2h}))
\]
where $N\in\mathbb{N}\cup\{\infty\}$ and $\infty>a_0>a_1>\cdots>a_{2N}>a_{2N+1}>\cdots>0$. Define
\begin{align}
&F:[0,\infty)\rightarrow\mathbb{R};t\mapsto V_\alpha(B(0,t))=2\pi\frac{t^{\alpha+2}}{\alpha+2}.\nonumber
\end{align}
Define the strictly increasing function
\begin{align}
&J:[0,\infty)\rightarrow\mathbb{R};t\mapsto P_\beta(B(0,F^{-1}(t))=\Big[(2\pi)^\gamma(\alpha+2)^{\beta+1}\Big]^{\frac{1}{\alpha+2}}t^{\frac{\beta+1}{\alpha+2}}.\nonumber
\end{align}
The function $J$ is continuous, strictly increasing and $J(0)=0$. Put $t_h:=F(a_h)$ for $h=0,1,\ldots$. Then $\infty>t_0>t_1>\cdots>0$. Put $B:=B(0,r)$ where $r:=F^{-1}(v)$ so that $V_f(B)=v$. Note that
\[
v=V_f(E)=\sum_{h=0}^N\Big\{F(a_{2h})-F(a_{2h+1})\Big\}=\sum_{h=0}^{2N+1}(-1)^{h}t_h.
\]
By \cite{McGillivray2018} Lemma 10.5,
\[
P_\beta(E)=2\pi\sum_{h=0}^{2N+1}a_h^{\beta+1}=\sum_{h=0}^{2N+1}J(t_h)
\geq J(\sum_{h=0}^{2N+1}(-1)^h t_h)=J(v)=P_\beta(B).
\]
This shows that the centred ball $B$ is a minimiser for (\ref{isoperimetric_problem}).

\smallskip

\noindent{\em Uniqueness.} Let $v>0$ and $E$ be any minimiser for (\ref{isoperimetric_problem}); this is essentially bounded by Corollary \ref{existence_and_boundedness_for_indices_in_curly_Q}. By Theorem \ref{C1_property_of_reduced_boundary} there exists an $\mathscr{L}^2$-measurable set $\widetilde{E}$ with the properties
\begin{itemize}
\item[(a)] $\widetilde{E}$ is a minimiser of (\ref{isoperimetric_problem});
\item[(b)] $L_{\widetilde{E}}=L_E$ a.e. on $(0,\infty)$;
\item[(c)] $\widetilde{E}$ is open with analytic boundary in $\mathbb{R}^2_0$;
\item[(d)] $\widetilde{E}=\widetilde{E}^{\mathrm{sc}}$.
\end{itemize}
Let $B$ be the open centred ball with the property that $V_\alpha(B)=V_\alpha(E)=v$ and denote its radius by $r>0$. Suppose that $\widetilde{E}\setminus B\neq\emptyset$. By Lemma \ref{boundary_of_E_as_union_of_circles} there exists $t>r$ such that $\mathbb{S}_t\subset\widetilde{M}$. So $P_\beta(E)=P_\beta(\widetilde{E})\geq 2\pi t^{\beta+1}>2\pi r^{\beta+1}=P_\beta(B)$. This contradicts the fact that $E$ is a minimiser for (\ref{isoperimetric_problem}). Consequently, $\widetilde{E}=B$. We infer that $L_E=L_B$ a.e. on $(0,\infty)$; in particular, $|E\Delta B|=0$. This entails that $E$ is equivalent to $B$.
\qed

\bigskip

\noindent{\em Acknowledgement.} I am grateful to Professor F. Morgan for bringing this problem to my attention.

\end{document}